\title{Numerical Stability on Local Integral Methods using RBF-QR}

\author{L. Ponzellini Marinelli$^{[1,2]}$\\
        ponzellini@cifasis-conicet.gov.ar\\
        \and
        N. Caruso$^{[1,2]}$\\
        caruso@cifasis-conicet.gov.ar\\
        \and
        M. Portapila$^{[1,2]}$\\
        portapila@cifasis-conicet.gov.ar\vspace{.5cm}\\
        $[1]$ Faculty of Exact Sciences, Engineering and Surveying,\\
        National University of Rosario,\\
        Rosario S2000BTP, Argentina.\\
        $[2]$ French Argentine International Center\\
        for Information and Systems Sciences,\\
        UAM (France)/UNR-CONICET,\\
        Rosario S2000BTP, Argentina.\\
}
\date{\today}

\documentclass[12pt]{article}

\usepackage{amssymb}
\usepackage{amsthm}

\usepackage{bm}
\usepackage{amsmath}

\usepackage{graphicx}

\begin{document}
\maketitle

\begin{abstract}
Many local integral methods are based on an integral formulation over small and heavilly overlapping stencils with local RBF interpolations. 
These functions have become an extremely effective tool for interpolation on scattered node sets, however the ill-conditioning of the interpolation matrix 
-when the RBF shape parameter tends to zero corresponding to best accuracy- is a serious task. Several stabilizing methods have been developed to deal with 
this near flat RBFs. The inclusion of the RBF-QR technique in the process of approximating in local integral methods makes possible to avoid this problem 
and stabilize the numerical error. In this paper we combine this technique in a local integral method and present accuracy results for Poisson, 
convection-difussion and thermal boundary layer PDEs.
\end{abstract}

\section{Introduction}\label{sec:intro}
The boundary element method (BEM) is now a well-established numerical technique in engineering.
The basis of this method is to transform the original partial differential equation (PDE), or 
system of PDEs that define a given physical problem, into an equivalent integral equation 
(or system) by means of the corresponding Green’s second identity and its fundamental solution, 
i.e. the Green's integral representation formula. In this way some or all of the field variables 
and their derivatives are only necessary to be defined at the boundary. 

Further increase in the number of applications of the BEM has been hampered by the need to operate 
with relatively complex fundamental solutions or by the difficulties encountered when these 
solutions cannot be expressed in a closed form. In the BEM formulation of this kind of problems, 
it is common to use an integral representation formula based upon a PDE with known closed-form 
fundamental solution, and express the remaining terms of the original equation as domain integrals. 
It is known that in these cases the BEM is in disadvantage in comparison with the classical domain 
schemes, such as the Control Volume (CV) and the Finite Element method (FEM). In the early BEM 
analysis the evaluation of domain integrals was done using cell integration, a technique which, 
while effective and general, made the approach too costly computationally due to the successive 
integration at each cell required for each of the surface collocation points. In order to deal with this, 
several methods have been developed in the literature to take domain integrals to the boundary in 
order to eliminate the need for internal cells (boundary-only BEM formulations). One of the most 
popular methods to date is the dual reciprocity method (DRM) introduced by Nardini and Brebbia 
\cite{nardini_brebbia_85}. In the DRM, the unknown densities of the corresponding domain
integrals are interpolated by a Radial Basis Function (RBF) scheme, and by applying the 
Green's second identity to a convolution integral of a particular solution and the fundamental solution, 
the domain integrals are converted into equivalent surface integrals. 
However, the DRM approach has the same computational limitations than the cell integration scheme, 
since very large fully populated matrix systems are obtained. It is important to mention that the 
DRM approximation is an alternative approach to evaluate domain integrals by defining global domain 
interpolations and only evaluating surface integrals, but still a domain integration scheme.  

When dealing with the BEM for large problems, with or without closed form fundamental solution, 
it is frequently used a domain decomposition technique, in which the original domain is divided 
into subdomains, and on each of them the full integral representation formulae are applied. At the 
interfaces of the adjacent subdomains the corresponding full-matching conditions are imposed 
(local matrix assembly), as is required in the CV and FEM methods, for which it is necessary to 
define subdomains or elements connectivity. However, in contrast with the CV and FEM methods, which 
integral representations of the original PDE are based on weighted residual approximation, 
in the BEM technique the Green's integral representation formula is an exact representation of 
the original PDE at each integration subdomain. The BEM matrices for subdomain formulation leads to 
block banded matrix systems with one block for each subregion and overlaps between blocks when 
subdomains have a common interface. In the limit of a very large number of subdomains, the resulting 
internal mesh pattern looks like a finite element grid. 

One of these approaches based on large number of subdomains but using the DRM to evaluate the domain 
integrals at each subdomain, instead of cell integration, has been referred by Popov and Power 
\cite{popov_power_99_IJNME} as the Dual Reciprocity Multi Domain approach (DRM-MD), for more details 
see Portapila and Power \cite{portapila_power_07}. As previously commented, the most attractive aspect of this 
type of local BEM approach at the subdomain level is the use of an exact integral representation formula 
of the original PDE instead of a weighted residual approximation. However, the numerical efficiency of this 
type of local BEM approaches is still behind of those classical domain numerical schemes. For this 
reason in recent years significant efforts have been given to the improvement of this type local BEM approaches.  

As has been the case in the FEM, see Atluri and Zhu \cite{atluri_zhu_98}, meshless formulations of 
local BEM approaches, see Zhu et al., \cite{zhu_zhang_atluri_98}, are attractive and efficient techniques 
to improve the performance of local BEM schemes. As in the meshless FEM, in the meshless BEM the integral 
representation formulae are applied at local internal integration subdomains embedded into interpolation 
stencils that are heavily overlapped. In this type of approach the continuity of the field variables are 
satisfied by the interpolation functions avoiding the local connectivity between subdomains or elements 
needed to enforce the matching conditions between them. Different interpolation schemes can be employed 
at the interpolation stencils, being the moving least squares shape functions and RBF interpolations the 
most popular approaches used in the literature. A major advantage of the meshless local BEM formulations 
in comparison with the classical BEM multi domain decomposition approaches, as the DRM-MD, is that the 
resulting integrands of the integral representation formulae are all regular, instead of singular, since 
the collocation points are always selected inside the integration subdomain. 

In the Local Boundary Integral Element Methods (LBEM or LBIEM) the solution domain is covered by a series 
of small and heavily overlapping local interpolation stencils, where a direct interpolation of the field 
variables is used to approximate the densities of the integral operator, and the boundary conditions 
of the problem are imposed at the integral representation formula; i.e. at the global system of equations, 
resulting in the evaluation of the corresponding weakly and singular surface integrals and if it is the 
case regular domain integrals, over each of the integration subdomains including those in contact with the 
problem boundary \cite{ zhu_zhang_atluri_98,sladek_sladek_zhang_04,sellountos_polyzos_atluri_12}. 
In this type of approach, the domains of integration usually are defined over several stencils, resulting 
in highly overlapping integration subdomains, in addition to the overlapping of interpolation stencils. 
Both polynomial moving least squares (MLS) approximation and direct RBF interpolations have been previously 
used in the LBEM as local interpolation algorithms. 

In Caruso et. al. \cite{caruso_portapila_power_15}, the Localized Regular Dual Reciprocity Method (LRDRM) 
is presented. The LRDRM is an integral domain descomposition method with two distinguishing features, the 
boundary conditions are imposed at the local interpolation (a local RBF interpolation) level and all the 
calculated integrals are regular. The "following" work \cite{caruso_portapila_power_16} is shown an 
enhancement of this method where the interpolation functions themselves satisfy the partial differential 
equation to be solve.

In recent years, the theory of RBFs has undergone intensive research and enjoyed considerable success as 
a technique for interpolating multivariable functions and for solving PDEs 
\cite{fasshauer_07_book,fasshauer_mccourt_15_book,fornber_flyer_15_book}. An RBF $\phi(r)$ depends only 
on the distance $r=\|\mathbf{x}-\mathbf{x}_k\|$ to a center node $\mathbf{x}_k$. The methods that use 
RBFs do not requiere a grid and it has been shown to be high-order accurate, flexible in nontrivial geometries, 
computationally efficient and easy to implement. 

When infinitely smooth RBFs $\phi(r,\varepsilon)$ are used, the spectral accuracy is often achieved 
when the shape parameter $\varepsilon$ tends to zero. This has been proven for some special cases 
\cite{madych_nelson_92,buhmann_dyn_93}, although numerical experiments suggestes that is also true in much general settings.
Nevertheless, in practice the interpolation error decreases to low levels until it breaks down due to the numerical ill-conditioning 
\cite{larsson_fornberg_03,larsson_fornberg_05}, i.e., when $\varepsilon \rightarrow 0$, the RBFs become relatively 
flat (named \emph{near-flat RBFs}) and the interpolation matrix increases the condition number. This was a -mistaken- trade-off between 
acuracy and numerical conditioning named as an \emph{uncertainty principle} due to R. Schaback \cite{schaback_95} which 
stablished that high accuracy and numerical stability cannot arrive simultaneously. 

This misconception about the uncertainty principle led to a negative impact on the development for RBFs approximation methods with scattered 
data. The reason was that the numerical solution denoted as \emph{RBF-Direct} amount to an ill-conditioned numerical procedure 
for a well-conditioned problem. So, many techniques for stabilizing the error has been developed in the last fiftteen year 
\cite{fornberg_wright_06,fornberg_piret_07,fornberg_larsson_flyer_11,fornberg_lehto_powell_13}. One of them, the \emph{RBF-QR algorithm} can 
stably compute interpolants in the case of near-flat RBFs using another basis that generates the same interpolation space.

The RBF-QR technique presented in \cite{fornberg_larsson_flyer_11} opened up new possibities for numerical methods based on local 
RBF approximations, such as LRDRM since it is possible to stabilize the shape parameter regime for small values of $\varepsilon$.

The following sections in this paper is structured as follows. In Section \ref{sec:formulation}, we describe different 
formulations of Local Integral Methods with local RBF interpolations. In Section \ref{sec:lim-rbf-qr} we describe the introduction of the 
RBF-QR technique into Local Integral approaches. And finally Section \ref{sec:num-results} contains a variety of numerical examples for 
Poisson's equations, convection-diffusion and thermal boundary layer equation. Section \ref{sec:conclusion} contains some concluding remarks.


\section{Local Integral Methods}\label{sec:formulation}

\subsection{Mathematical formulation and boundary integral represention formulae.}
\label{sec:math_and_intFormula}

Let us consider the following elliptic problem on a bounded open domain $\Omega$:
\begin{equation}
\left\{\begin{array}{cccl}
\mathcal L \left[u \left( \mathbf{x}\right)\right] & = & f\left(\mathbf{x}\right) & \hspace{1cm} \mathbf{x} \in \Omega,\\
\mathcal B \left[u \left( \mathbf{x}\right)\right] & = & g\left(\mathbf{x}\right) & \hspace{1cm} \mathbf{x} \in \Gamma=\partial \Omega,\\
\end{array}
\right.
\label{Eq:elliptic_Prob} 
\end{equation}
where $\mathcal L[.]$ is an elliptic operator and $\mathcal B[.]$ is a classical boundary operator related with different kind 
of boundary conditions (e.g. Dirichlet, Neumann or Robin Condition). We assume that the partial differential equation can be 
rewritten in the following way:
\begin{equation}
\Delta u \left( \mathbf{x}\right) = b\left(\mathbf{x},u\left(\mathbf{x}\right),\nabla u \left(\mathbf{x}\right) \right).
\label{Eq:gov_equation}
\end{equation}
The integral representation formula for the above PDE in terms of the Laplace's fundamental solution is 
obtained from the Green's second identity in terms of the  superposition of surfaces (single and double layers) 
and volume potentials is given by
\begin{equation}
c\left(\xi\right)u\left(\xi\right) = 
\int_{\Gamma}q^{*}\left(\mathbf{x},\xi\right)u\left(\mathbf{x}\right)d\Gamma_{\mathbf{x}} -
\int_{\Gamma}u^{*}\left(\mathbf{x},\xi\right)q\left(\mathbf{x}\right)d\Gamma_{\mathbf{x}} +
\int_{\Omega} u^{*}\left(\mathbf{x},\xi\right)\ b \ d\Omega_{\mathbf{x}},
\label{Eq:int_eq_conv-diff}
\end{equation}
with $\xi$ as the evaluation point, also referred as collocation point, and $u^{*}\left(\mathbf{x},\xi\right)$ 
as the fundamental solution of the Laplace problem, which in the case of two-dimensional problems is given by:
\begin{equation}
u^{*}\left(\mathbf{x},\xi\right) = 
\frac{1}{2\pi} \ln \left(\frac{1}{R\left(\mathbf{x},\xi\right)}\right),
\label{Eq:fund_sol_Laplace_2d}
\end{equation}
where $R\left(\mathbf{x},\xi\right)$ is the distance between the integration points $\mathbf{x}$ and collocation 
point $\xi$, i.e., $R\left(\mathbf{x},\xi\right)=\|\mathbf{x}-\xi\|$, and 
$q^{*}\left(\mathbf{x}\right)=\frac{\partial u^{*}}{\partial n}\left(\mathbf{x},\xi\right)$. The constant value 
$c\left(\xi\right)\in\left[0,1\right]$, being 1 if the point $\xi$ is inside the domain and $\frac{1}{2}$ if the point 
$\xi$ is on a smooth part of the domain boundary $\Gamma$ (for this work we always will consider $c\left(\xi\right)=1$).

The integral representation formula (\ref{Eq:int_eq_conv-diff}) is the basis of any meshless BEM approach, 
where the integration surface $\Gamma$ and domain $\Omega$ are chosen as integration subregions, $\Gamma_i$ and $\Omega_i$, 
embedded inside of a corresponding interpolation stencils, which are heavily overlapped. Despite of the above formulation, 
instead of using the fundamental solution, $u^{*}\left(\mathbf{x},\xi\right)$, and its normal derivative, 
$q^{*}\left(\mathbf{x},\xi\right)$, the Dirichlet Green's Function (DGF), $G\left(\mathbf{x},\xi\right)$ and its corresponding 
normal derivative, $Q\left(\mathbf{x},\xi\right)$, can be used, it leads that the Eq. (\ref{Eq:int_eq_conv-diff}) at each 
integration subregion is reduced to:
\begin{equation}
u\left(\xi\right) = 
\int_{\Gamma_i}Q\left(\mathbf{x},\xi\right)u\left(\mathbf{x}\right)d\Gamma_{\mathbf{x}} +
\int_{\Omega_i}  G\left(\mathbf{x},\xi\right) \ b \ d\Omega_{\mathbf{x}},
\label{Eq:int_eq_conv-diff(Green)}
\end{equation}
where by definition over the surfaces $\Gamma_i$ the value of $G$ is identically zero. 

\begin{figure}[!ht]
\centering
\begin{tabular}{c}
\includegraphics[width=120mm]{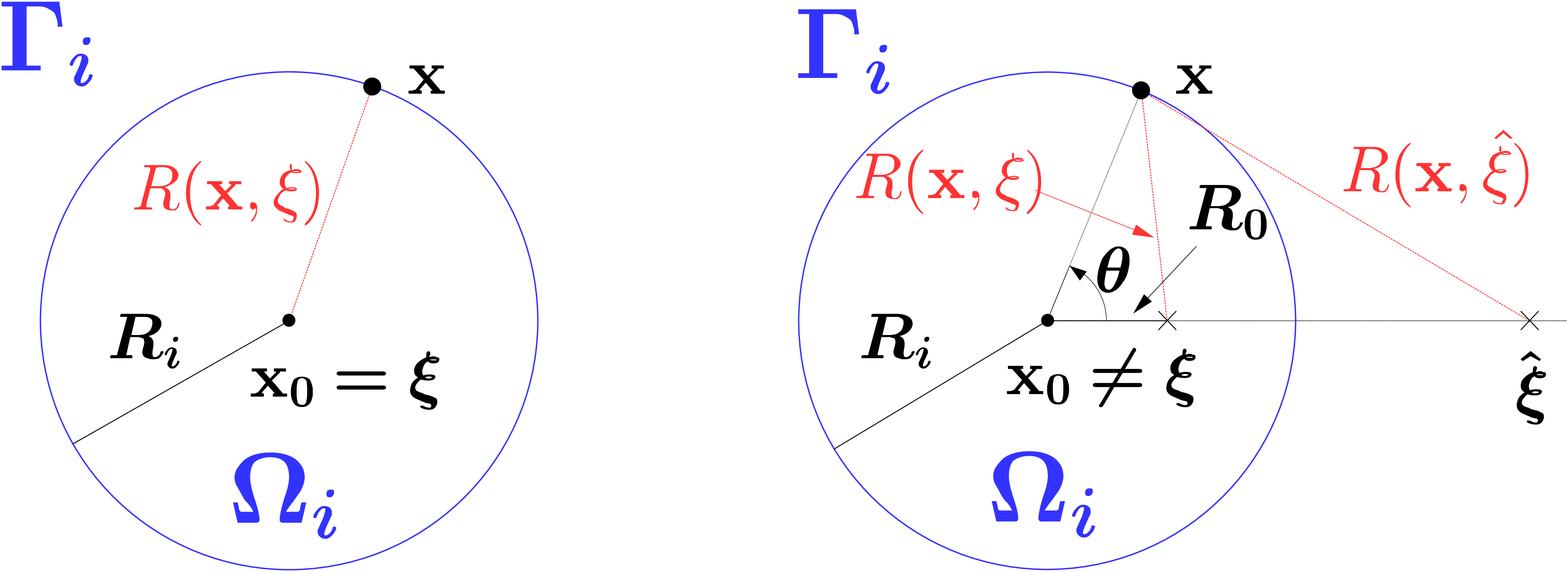}
\end{tabular}
\vspace{-.5cm}
\caption{Schematic representation of Dirichlet Green's Function in a circle for centered source $\mathbf{x}_0=\mathbb{\xi}$ (left) 
and outskirt source point $\mathbf{x}_0\neq\mathbb{\xi}$ (right).}
\label{F:GreenFunction}
\end{figure}

In the case a two dimensional problem and a circular subregion of integration $\Omega_i$ with radius $R_i$ and centre 
$\mathbf{x}_0$, the Dirichlet Green's Function for a source point, $\mathbf{\xi}$, inside the circle can be obtained from 
the circle theorem, and given by the Eq. (\ref{Eq:Green-function}) (see Figure \ref{F:GreenFunction} for an schematic 
representation about the elements of the DGF and the details about the circle theorem can be find at \cite{milne_thomson_68}):
\begin{equation}
G\left(\mathbf{x},\xi \right)
=\left\{
\begin{array}{ccc}
\dfrac{1}{2\pi}\ln \left(\dfrac{R_i}{R(\mathbf{x},\xi)}\right) & \hspace{1cm} if & \xi = \mathbf{x_0},\\
\\
\dfrac{1}{4\pi}\ln \left(\dfrac{R_0^2\;R(\mathbf{x},\hat{\xi})^2}{R_i^2\; R(\mathbf{x},\xi)^2}\right) & \hspace{1cm} if & \xi \neq \mathbf{x_0},\\
\end{array}
\right.
\label{Eq:Green-function}
\end{equation}
with the image or reflection point, $\mathbf{\hat{\xi}}$, located outside the circle along the same ray of the source point. 
In the above expression $R(\mathbf{x},\xi)$ is the distance between the field point $\mathbf{x}$ and the source point $\mathbf{\xi}$, 
similarly $R(\mathbf{x},\hat{\xi})$ is the distance between $\mathbf{x}$ and the image point $\hat\xi$ and $R_0$ is the distance 
between $\mathbf{x}_0$ and $\mathbf{\xi}$. In the work \cite{power_caruso_portapila_17} there are more details on the features and 
use of the DGF on the LRDRM. 

There are different approaches to evaluate the corresponding volume integrals in (\ref{Eq:int_eq_conv-diff}) and/or 
(\ref{Eq:int_eq_conv-diff(Green)}) in the literature of BEM or Local BEM, one of them consists in approximating 
the density $b$ of the volume integrals in terms of a interpolation function, i.e.  
\begin{equation}
 b\approx \sum_{k=1}^{N} \beta_k  \varphi_k\left(\mathbf{x}\right)
 \label{Eq:int_term_no_hom}
\end{equation}
with $N$ as the number of interpolation points and $\varphi_k\left(\mathbf{x}\right)$ usually is defined by a RBF. 
Then the integral representation formulae obtained are:
\begin{eqnarray}
u\left(\xi\right) & = &
\int_{\Gamma_i}q^*\left(\mathbf{x},\xi\right)u\left(\mathbf{x}\right)d\Gamma_{\mathbf{x}} -
\int_{\Gamma_i}u^*\left(\mathbf{x},\xi\right)\frac{\partial u^*}{\partial n}
\left(\mathbf{x}\right)d\Gamma_{\mathbf{x}} \nonumber \\
& + & \sum_{k=1}^{N}\beta_k 
\int_{\Omega_i} u^*\left(\mathbf{x},\xi\right)
\varphi_k\left(\mathbf{x}\right)d\Omega_{\mathbf{x}}
\label{Eq:Local_BEM-meshless-i}
\end{eqnarray}    
or
\begin{eqnarray}
u\left(\xi\right) = \int_{\Gamma_i}Q\left(\mathbf{x},\xi\right)
u\left(\mathbf{x}\right)d\Gamma_{\mathbf{x}} + 
\sum_{k=1}^{N}\beta_k 
  \int_{\Omega_i}G\left(\mathbf{x},\xi\right)
  \varphi_k\left(\mathbf{x}\right)d\Omega_{\mathbf{x}}.
 \label{Eq:DRM-meshless-i}
\end{eqnarray}   

Futhermore if it is possible find a particular solution $\widetilde{\varphi}_k$ such that, 
$\Delta \widetilde{\varphi}_k\left(\mathbf{x}\right)=\varphi_k\left(\mathbf{x}\right)$, then applying again the Green's second 
identity to the resulting volume integral with the particular solution as density and the fundamental solution as kernel or DGF, 
the formulae obtained for each subregion are:
\begin{eqnarray}
u\left(\xi\right) & = &
\int_{\Gamma_i}q^*\left(\mathbf{x},\xi\right)u\left(\mathbf{x}\right)d\Gamma_{\mathbf{x}} -
\int_{\Gamma_i}u^*\left(\mathbf{x},\xi\right)\frac{\partial u}{\partial n}
\left(\mathbf{x}\right)d\Gamma_{\mathbf{x}} \nonumber \\
& + & 
\sum_{k=1}^{N}\beta_k \left\{ \widetilde{\varphi}_k\left(\xi\right)  
 - \int_{\Gamma_i}q^*\left(\mathbf{x},\xi\right)\widetilde{\varphi}_k \left(\mathbf{x}\right)d\Gamma_{\mathbf{x}}
+  \int_{\Gamma_i}\frac{\partial u^*}{\partial n}\left(\mathbf{x},\xi\right)
\widetilde{\varphi}_k\left(\mathbf{x}\right)d\Gamma_{\mathbf{x}}\right\}
\label{Eq:Local_BEM-meshless}
\end{eqnarray}    
or
\begin{eqnarray}
u\left(\xi\right) = 
\int_{\Gamma_i}Q\left(\mathbf{x},\xi\right)u\left(\mathbf{x}\right)d\Gamma_{\mathbf{x}} 
+ \sum_{k=1}^{N}\beta_k \left\{  
\widetilde{\varphi}_k\left(\xi\right)  
 - \int_{\Gamma_i}Q\left(\mathbf{x},\xi\right)\widetilde{\varphi}_k\left(\mathbf{x}\right)
 d\Gamma_{\mathbf{x}}\right\}.
\label{Eq:DGF-DRM-meshless}
\end{eqnarray}    
  
\subsection{Local Integral approaches I: Localized Regular Dual Reciprocity Method (LRDRM)}

\begin{figure}[!ht]
\centering
\begin{tabular}{c}
\includegraphics[width=120mm]{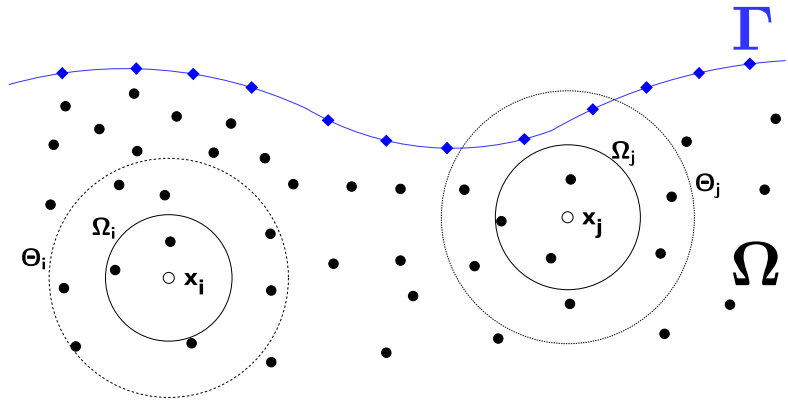}
\end{tabular}
\caption{Schematic representation of local stencils $\Theta_i$ and local subdomains or subregions $\Omega_i$ for different local meshless BEM }
\label{F:local-stencil-LBEM}
\end{figure}

In the local meshless BEM approaches the integral representation formulae are applied at local internal subdomain 
or subregion (as $\Omega_i$ and $\Omega_j$  in the Figure \ref{F:local-stencil-LBEM}) embedded into 
interpolation stencils that are heavily overlapped (as $\Theta_i$ and $\Theta_j$ in the Figure \ref{F:local-stencil-LBEM}).

In this article we consider the approach presented in Caruso et al. \cite{caruso_portapila_power_15} and Power et al. 
\cite{power_caruso_portapila_17}, i.e. the LRDRM, the computational domain is covered by a set of integration subregions 
(where the integral representation formula is applied) and a set of heavily overlapping interpolation stencils with the 
important features; the first one, the boundary conditions of the problem are imposed locally at the interpolation stencil; 
and the second one, the collocation point (where the Delta function is centered) is allways inside the subregion in order to 
obtain that every boundary integral in the local integral approach is regular.

Then the field variable $u\left(\mathbf{x}\right)$ is approximated by a RBF interpolation using the corresponding 
nodes of the interpolation stencil, plus some additional auxiliary boundary points (if the stencil is next to the global boundary). 
The set $\{\left( \mathbf{x}_j,u(\mathbf{x}_j) \right)\}_{j=1}^{n_i}$ is formed by the internal nodes $\mathbf{x}_j$ and 
the corresponding unknown nodal values $u(\mathbf{x}_j)$ for $j=1,\dots,n_i$. The other set 
$\{\left( \mathbf{x}_j,\mathcal B\left(u(\mathbf{x}_j)\right) \right)\}_{j=n_i+1}^{n_i+n_b}$ has the boundary nodes 
$\mathbf{x}_j$ and and the boundary data $\mathcal B\left(u(\mathbf{x}_j)\right)$ for $j=n_i+1,\dots,n_i+n_b$. Then the local approximation 
is presented as
\begin{equation}
 u\left(\mathbf{x}\right)=\sum^{n}_{j=1}\alpha_j\varphi_j\left(\mathbf{x}\right)
 \label{u_def_bound}
\end{equation}
being $\varphi_j\left(\mathbf{x}\right)=\phi\left(\| \mathbf{x} - \mathbf{x}_j\|\right)$ the RBFs interpolating functions and $n=n_i+n_b$ where $n_i$ is the number of internal points of the stencil and $n_b$ is the number of auxiliary points given by the boundary points belonging to 
an interpolation stencil next to the problem boundary. In this way, at interpolation stencils inside the problem domain $n_b=0$, and  at 
those in contact with the problem boundary, $n_b$ is equal to the number of boundary collocation points belonging to the given stencil.

Different RBFs have been considered as local interpolating functions in integral methods. In \cite{caruso_portapila_power_15} 
were used the MQ1, MQ2 and TPS. Other important and well stablished RBF are the Gaussians that depends on a shape parameter $\varepsilon$. All these RBFs 
corresponds to the two main groups of RBFs: piecewise smooth and infinitely smooth. Some examples are given in Table \ref{T:RBFs}. 
For these examples and many others, it had been demonstrated that the linear system of equations formed in the local 
interpolation is non-singular, in any dimension and for any number of nodes. This set of nodes must be different and \emph{unisolvent}. 
For more details, see \cite{fasshauer_07_book,fornber_flyer_15_book}.

\begin{table}[!ht]
\centering
\begin{small}
\begin{tabular}{|l|c|l|}
\hline
Infinitely smooth RBFs & $\phi(r,\varepsilon)$    & Polynomial \\
\hline
Gaussian (GA)          & $ e^{-(\varepsilon r)^2} $    &  no poly  \\ 
Multiquadric 1 (MQ1)   & $\sqrt{1 + (\varepsilon r)^2}$ &  $p(x)\in \mathbb{P}_0(x)$  \\ 
Multiquadric 2 (MQ2)   & $\left(1 + (\varepsilon r)^2\right)^{\frac{3}{2}}$ &  $p(x)\in \mathbb{P}_1(x)$  \\ 
\hline
Piecewise smooth RBFs  & $\phi(r)$ & Polynomial \\
\hline
Thin Plate Spline (TPS) & $r^4 log(r)$ &  $p(x)\in \mathbb{P}_2(x)$  \\ 
\hline
\end{tabular}
\caption{Piecewise (TPS) and infinitely smooth (GA, MQ1, MQ2) RBFs with the corresponding polynomials to ensure non-singularity 
of the interpolation matrix used in this work. $r=\|\mathbf{x}-\mathbf{x}_j\|$ denotes the distance to the centre of the RBF and $\varepsilon$ the shape parameter.}
\label{T:RBFs}
\end{small}
\end{table}

By using this interpolation scheme, the value of the unknown $u$ in (\ref{u_def_bound}), over an integration subregion is 
obtained from the interpolation reconstruction formula as:  

\begin{equation}
u\left(\mathbf{x}\right)=
{\bm \varphi\left(\mathbf{x}\right)}^T \mathbf{A}^{-1} \mathbf{d},
\label{Eq:u_phi}
\end{equation} 
with $\mathbf{A}^{-1}$ as the corresponding inverse interpolation matrix; and the vector $\mathbf{d}$ in terms of the unknown 
nodal values $\mathbf{u}=\left[u_1,\ldots,u_{n_i}\right]$ and the prescribed boundary condition values 
$g(\mathbf{u}_b)=\left[g(\mathbf{x}_{n_i+1}),\ldots,g(\mathbf{x}_{n})\right]$ with $g$ from Eq. 
(\ref{Eq:elliptic_Prob}) (i.e., $\mathbf{d}=\mathbf{u}^T$ for internal stencils and 
$\mathbf{d}=\left[\mathbf{u},g(\mathbf{u}_b)\right]^T$ for boundary stencils). After inversion of the interpolation matrix $A$, 
the interpolation coefficients are given by: 

\begin{equation}
\mathbf{\alpha}=\mathbf{A}^{-1} \mathbf{d}.
\label{alfa}
\end{equation}

In cases where the non-homogeneous term $b$ in (\ref{Eq:gov_equation}) is function of the derivative of the field variable $u$, 
we use the generalized finite different approximation where this value is approximated by the derivative of the interpolation 
reconstruction function in terms of the neighbouring values of $u$ at the interpolation stencils, i.e.,

\begin{equation} 
\frac{\partial u\left(\mathbf{x}\right)}{\partial x_{i}}=
{\frac{\partial {\bm \varphi}\left(\mathbf{x}\right)}{\partial x_i}}^T \mathbf{A}^{-1} \mathbf{d}.
\label{E:local.interp.w.matrix}
\end{equation}

Substituting (\ref{u_def_bound}) into the integral formula (\ref{Eq:DGF-DRM-meshless}) with $\xi=\mathbf{x}_i$ a trial point 
inside $\Omega_i$, the integration subregion, the discretized form for the unknown $u_i=u\left(\mathbf{x}_i\right)$ reduces to:

\begin{equation}
 u_i= \sum^n_{j=1}\alpha_j h_{ij} 
      + \sum^n_{j=1}\beta_j \widetilde{h}_{ij},
 \label{disc_form}  
\end{equation}

or in matrix notation:

\begin{equation}
   u_i = \mathbf{h}_i^T\mathbf{\alpha}  + \mathbf{\widetilde{h}}_i^T\mathbf{\beta},
   \label{matrix_form}
\end{equation}
where 
\begin{eqnarray}
      h_{ij} & = & \int_{\Gamma_{i}}Q\left(\mathbf{x},\mathbf{x}_i\right)
            \varphi_j\left(\mathbf{x}\right)d\Gamma_{\mathbf{x}} \\
      \widetilde{h}_{ij} & = & \widetilde{\varphi}_j\left(\mathbf{x}_i\right)- 
      \int_{\Gamma_{i}} Q\left(\mathbf{x},\mathbf{x}_i\right)
      \widetilde{\varphi}_j\left(\mathbf{x}\right)d\Gamma_{\mathbf{x}}
\end{eqnarray}
with the column vectors $\mathbf{h}_i=\left[\ldots,h_{ij},\ldots\right]^T$ and 
$\mathbf{\widetilde{h}}_i=[\ldots,\widetilde{h}_{ij},\ldots]^T$.

All the integrals in the above formulations are regulars, since the collocation points are located inside the integration subregion, 
and they are evaluated  through the Gauss-Legendre quadrature. Also notice that the interpolation coefficients $\alpha$ and $\beta$ 
in equation (\ref{matrix_form}) correspond to the interpolation of the field variable and the non-homogeneous part of the PDE, 
respectively, both of them given in terms of the stencil nodal values of the field variable.

When possible, in the DRM interpolation it is considered that

\begin{equation}
 b\left(\mathbf{x},u\left(\mathbf{x}\right),\nabla u\left(\mathbf{x}\right)\right)=
 f\left(\mathbf{x}\right) + \widetilde{b}\left(u\left(\mathbf{x}\right),\nabla u \left(\mathbf{x}\right)\right)
 \approx \sum^n_{j=1}\beta_j\varphi_j\left(\mathbf{x}\right)
 \label{b_def_int}
\end{equation}
similarly to equation (\ref{alfa}) beta is written as: $\beta = \mathbf{\widetilde{A}^{-1}} (\mathbf{f}_i+\mathbf{\widetilde{b}}_i)$, 
where the vector $\mathbf{f}_i$ is given by a data function evaluation and $\mathbf{\widetilde{b}}_i$ can be written as a function 
of $\mathbf{d}$ (i.e. in terms of the nodal values in $\mathbf{u}$ and boundary conditions values $g(\mathbf{u}_b)$), 
from the following expression for $\widetilde{b}\left(u\left(\mathbf{x}\right),\nabla u\left(\mathbf{x}\right) \right)$ and its linearity:

\begin{equation}
 \widetilde{b}\left(u\left(\mathbf{x}\right),\nabla u\left(\mathbf{x}\right) \right)=
 \widetilde{b}\left(\varphi_j\left(\mathbf{x}\right),
 \nabla  \varphi_j\left(\mathbf{x}\right) \right) \mathbf{A}^{-1} \mathbf{d},
\label{eq_b_2}
\end{equation}
therefore

\begin{equation}
 \beta = \mathbf{\widetilde{A}^{-1}}\left(\mathbf{f}_i + \mathbf{A}_b\mathbf{A}^{-1} \mathbf{d}\right)
\label{eq_beta_2}
\end{equation}
with matrix coefficients $(\mathbf{A}_b)_{kj}=\widetilde{b}\left(\varphi_j\left(\mathbf{x}_k\right),\nabla\varphi_j\left(\mathbf{x}_k\right)\right)$.

In the above expression the matrices $\mathbf{A}$ and $\mathbf{\tilde{A}}$ are identical at stencils in the interior 
of the problem domain, however, at boundary stencils the matrix $\mathbf{A}$ is defined by the corresponding interpolation 
matrix according to the boundary conditions of the problem, while the matrix $\mathbf{\tilde{A}}$ is the same direct 
interpolation matrix. We note $\mathbf{A_b}$ as the matrix corresponding to calculus of the vector $\mathbf{\widetilde{b}}_i$.

Equation (\ref{matrix_form}) can be written in terms of $\mathbf{d}$ by substituting into it the expression (\ref{alfa}), 
resulting the following equation:

\begin{equation}
 u_i = \mathbf{\widetilde{h}}_i^T \mathbf{\widetilde{A}}^{-1} \mathbf{f}_i +\left(\mathbf{h}_i^T \mathbf{A}^{-1} +
 \mathbf{\widetilde{h}}_i^T \mathbf{\widetilde{A}}^{-1} \mathbf{A_b} \mathbf{A}^{-1}\right)\mathbf{d}.
\label{local_int_eq}
\end{equation}

Finally equation (\ref{local_int_eq}) is collocated at each trial point of each stencil to form a global sparse matrix system. 
This equation is obtained in a way that is possibly to avoid calculating numerically the inverse matrix $\mathbf{A}^{-1}$ and 
$\mathbf{\widetilde{A}}^{-1}$. This method is called in the rest of this work as the \emph{Localized Regular Dual Reciprocity Method (LRDRM)}.

\subsection{Local Integral approaches II: Local Integral Method (LIM)}

The main aim in the using of DRM in a global integral method as BEM is try to avoid the cost of numerical calculus of domain 
integrals, because in this kind of method with $N$ degrees of freedom it is necessary $N^2$ of operations whereas for
boundary integrals is $N$, but in local method with a local $n$ (a fix low number) degrees of freedom it does not seem 
a big deal. 

From Eq. (\ref{b_def_int}) the non-homogeneous term was splitted up in a known data function $f$ plus a linear unknown term $\tilde{b}$, so in order to get more accuracy, the data function $f$ is integrated directly (with the corresponding $DGF$) instead of approximated it

\begin{eqnarray}
u\left(\mathbf{\xi}\right) = \int_{\Gamma_i}Q\left(\mathbf{x},\mathbf{\xi}\right)u\left(\mathbf{x}\right)d\Gamma_i + 
\int_{\Omega_i} f\left( \mathbf{x} \right) Q\left(\mathbf{x},\mathbf{\xi}\right) d\Omega_i + \int_{\Omega_i} \tilde{b}\left( u, \nabla u \right) Q\left(\mathbf{x},\mathbf{\xi}\right) d\Omega_i
\label{int_formulation_split}
\end{eqnarray}
where the linear term $\tilde{b}$ is locally interpolated with RBFs.

Then from the local approach in the above equation and using the same local RBF interpolation scheme $\Theta_i$ we can obtain the following equation:

\begin{equation}
 u_i= \sum^n_{j=1}\alpha_j h_{ij} + \sum^n_{j=1}\beta_j \widetilde{\widetilde{h}}_{ij} +f_i,
 \label{disc_form-ii}  
\end{equation}
where $\alpha_j,\beta_j$ and $h_{ij}$ are calculated as in Eq. (\ref{disc_form}) and 

\begin{eqnarray}
      \widetilde{\widetilde{h}}_{ij} & = & \int_{\Omega_i}G\left(\mathbf{x},\mathbf{x}_i\right)
      \varphi_j\left(\mathbf{x}\right)d\Omega_{\mathbf{x}}, \nonumber \\
      f_i & = & \int_{\Omega_i}G\left(\mathbf{x},\mathbf{x}_i\right)
      f\left(\mathbf{x}\right)d\Omega_{\mathbf{x}}.
\end{eqnarray}
The obtained the equation is

\begin{equation}
 u_i = f_i +\left(\mathbf{h}_i^T \mathbf{A}^{-1} +
 \mathbf{\widetilde{\widetilde{h}}}_i^T \mathbf{\widetilde{A}}^{-1} \mathbf{A_b}
 \mathbf{A}^{-1}\right)\mathbf{d}
\label{local_int_eq-ii}
\end{equation}
which is collocated at each trial point of each stencil to form a global sparse matrix system. 

As before, to avoid calculating numerically the inverse matrix $\mathbf{A}^{-1}$ and 
$\mathbf{\widetilde{A}}^{-1}$ we rewrite this expression as
\begin{equation}
 u_i = f_i +\mathbf{z}^T_i\mathbf{d}
\label{local_int_eq-ii_zi}
\end{equation}
where the algorithmic procedure to calculate this equation is the following:
\begin{enumerate}
 \item[] Step 1. Solve $\mathbf{\tilde{A}}\mathbf{\tilde{w}}_i=\mathbf{\tilde{h}}_i$ (since $\mathbf{\tilde{A}}$ simetric).
 \item[] Step 2. $\mathbf{w}_i^T = \mathbf{h}_i^T + \mathbf{\tilde{w}}_i^T \mathbf{A_b}$.
 \item[] Step 3. Solve $\mathbf{A}^T\mathbf{z}_i = \mathbf{w}_i$.
\end{enumerate}
In Steps 1 and 3 the ill-conditioning of the linear systems could significant. 
We call this method \emph{Local Integral Method (LIM)}.

All these equations are ensamble resulting a sparse linear system. In ths paper we used an iterative solver, as the restarted GMRES method that has computational cost of the order $O(\gamma N^\beta Nn)$ with $\gamma$ and $\beta$ depending on the structure of the matrix and numerical scheme employed (for detail about the computational cost of the GMRES scheme 
used in this work see Guttel and Pestana \cite{guttel_pestana_14}).

\section{Introducing RBF-QR into a Local Integral Method}
\label{sec:lim-rbf-qr}

\subsection{The Local Integral RBF-QR Method}

In this section, we introduce the RBF-QR method presented in \cite{fornberg_piret_07,fornberg_larsson_flyer_11,larsson_lehto_heryudono_fornberg_13} 
into the local RBF interpolations for the Local Integral Methods developed in Section \ref{sec:formulation}
 to get a new formulae 
for local meshless methods. The principal porpouse is to bypass the ill-conditioning of the RBF-Direct approach (present in the LRDRM) for near-flat RBFs.

As it is explain later we change the Dual Reciprocity formulation to introduce the RBF-QR 
to get a new formulation called \emph{Local Integral RBF-QR Method (LIM RBF-QR)}.

\begin{figure}[!ht]
\begin{center}
\includegraphics [height=8cm, width=9cm] {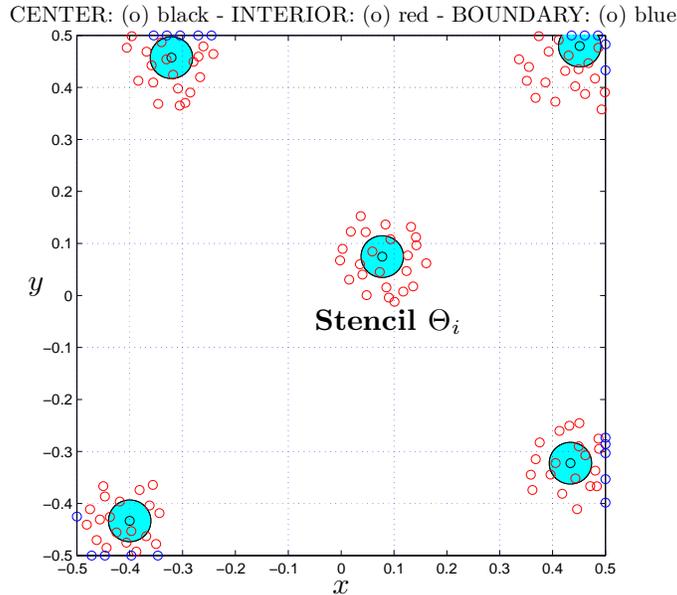}
\vspace{-.5cm}
\caption{Stencils $\Theta_i$ and integration regions $\Omega_i$}
\label{Fig:stencil_&_integration_regions}
\end{center}
\end{figure}

From Eq. (\ref{u_def_bound}), the unknown $u$ field is interpolated locally at each stencil $\Theta_i$ with RBF 
interpolating functions $\varphi_j\left(\mathbf{x},\varepsilon\right)=\phi\left(\| \mathbf{x} - \mathbf{x}_j\|,\varepsilon\right)$ that depends 
on a shape parameter $\varepsilon > 0$ where $\mathbf{\xi}_j \in \Theta_i$ the collocation point and $\mathbf{x}_j \in \Theta_i$ for $j=1,\dots,n_i+n_b$. 
See Fig. \ref{Fig:stencil_&_integration_regions}.

If there is no boundary node in $\Theta_i$ (as the interior stencil in that Figure), the interpolation matrix of the linear system for the 
local approximation (\ref{u_def_bound}) takes the form
\begin{equation}
\mathbf{A} = 
\left[
\begin{array}{cccc}
\phi(\|\mathbf{x}_1-\mathbf{x}_1\|,\varepsilon) & \phi(\|\mathbf{x}_1-\mathbf{x}_2\|,\varepsilon) & ... & \phi(\|\mathbf{x}_1-\mathbf{x}_{n_i}\|,\varepsilon)\\
\phi(\|\mathbf{x}_2-\mathbf{x}_1\|,\varepsilon) & \phi(\|\mathbf{x}_2-\mathbf{x}_2\|,\varepsilon) & ... & \phi(\|\mathbf{x}_2-\mathbf{x}_{n_i}\|,\varepsilon)\\
\vdots & \vdots & \ddots & \vdots \\
\phi(\|\mathbf{x}_{n_i}-\mathbf{x}_1\|,\varepsilon) & \phi(\|\mathbf{x}_{n_i}-\mathbf{x}_2\|,\varepsilon) & ... & \phi(\|\mathbf{x}_{n_i}-\mathbf{x}_{n_i}\|,\varepsilon)
\end{array}
\right].
\label{interpolation_matrix}
\end{equation}

As it is known, the RBFs constitute an ill-conditioned basis in a good approximation space. When the shape parameter 
$\varepsilon$ tends to zero, the interpolation error often decreases to low levels until it break downs \cite{larsson_fornberg_03,schaback_95} 
when solving the linear system of the interpolation with a direct method. This is because the interpolation matrix (\ref{interpolation_matrix}) 
becomes increasingly ill-conditioned and the expansion coefficients $\alpha_j$ in (\ref{u_def_bound}) becomes large magnitude and 
oscillatory causing numerical cancellations when using the reconstruction formula in (\ref{Eq:int_term_no_hom}) and (\ref{u_def_bound}).

\begin{figure}[!ht]
\centering
  \begin{center}
    \includegraphics[scale=0.225]{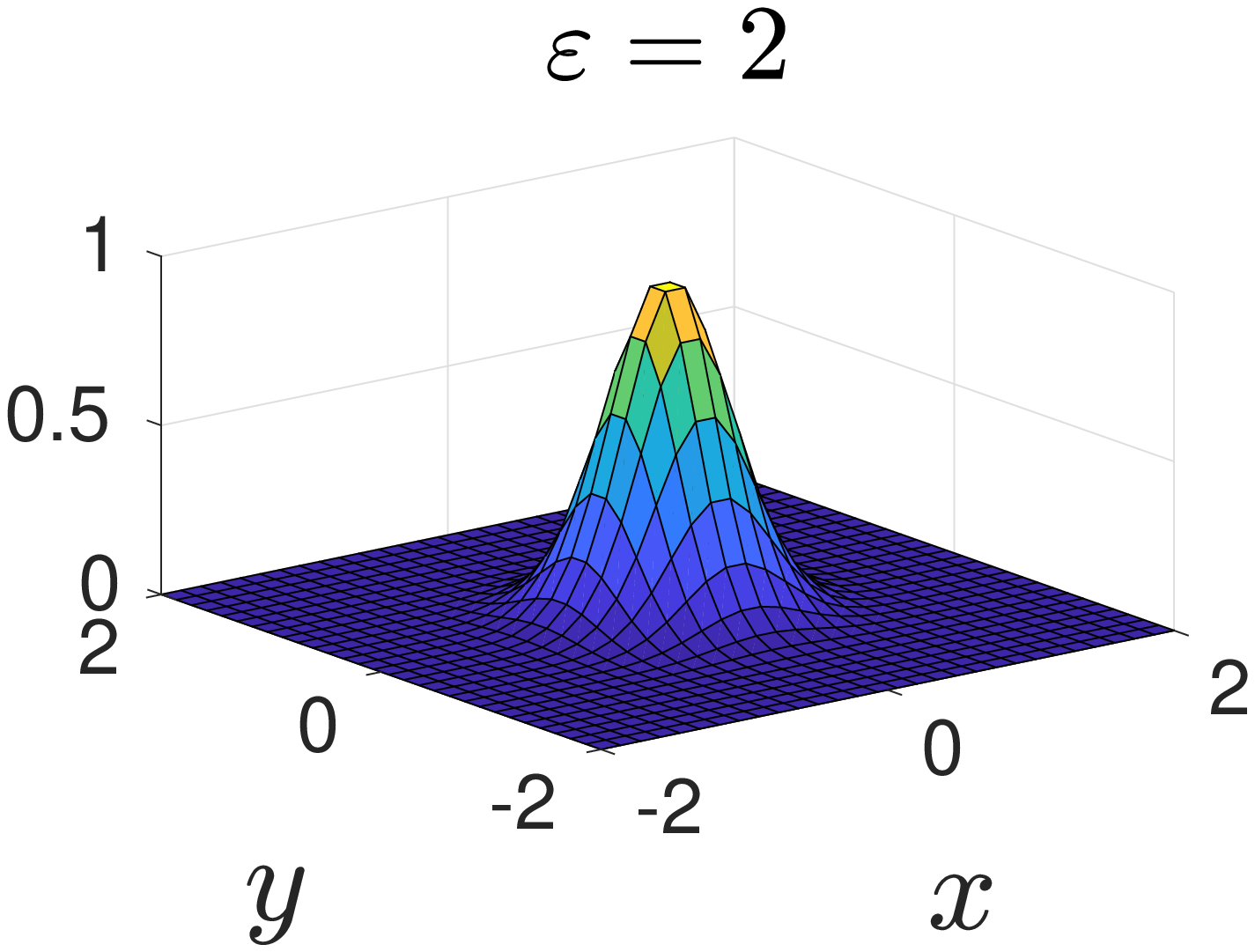}
    \includegraphics[scale=0.225]{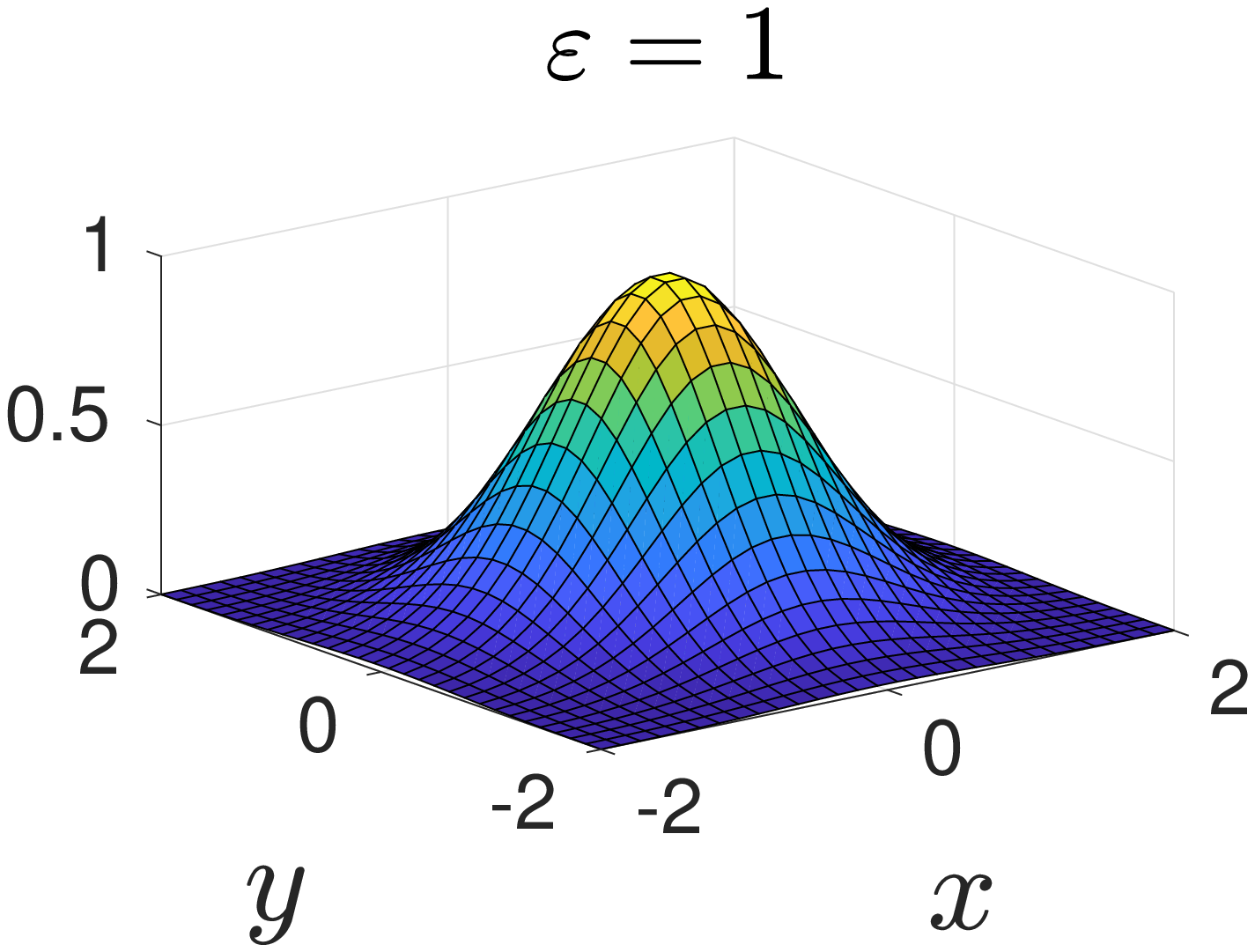}
    \includegraphics[scale=0.225]{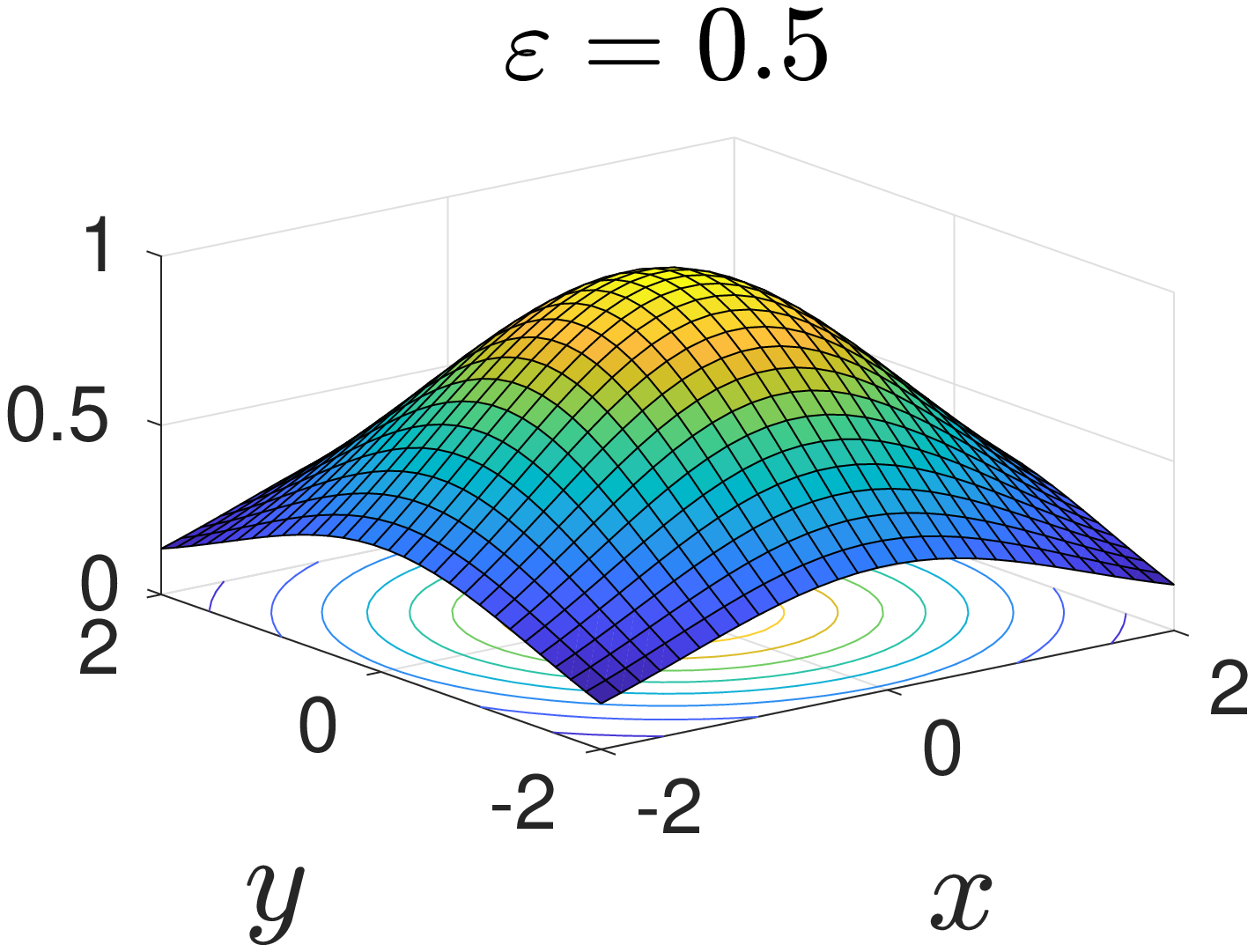}    
    \includegraphics[scale=0.225]{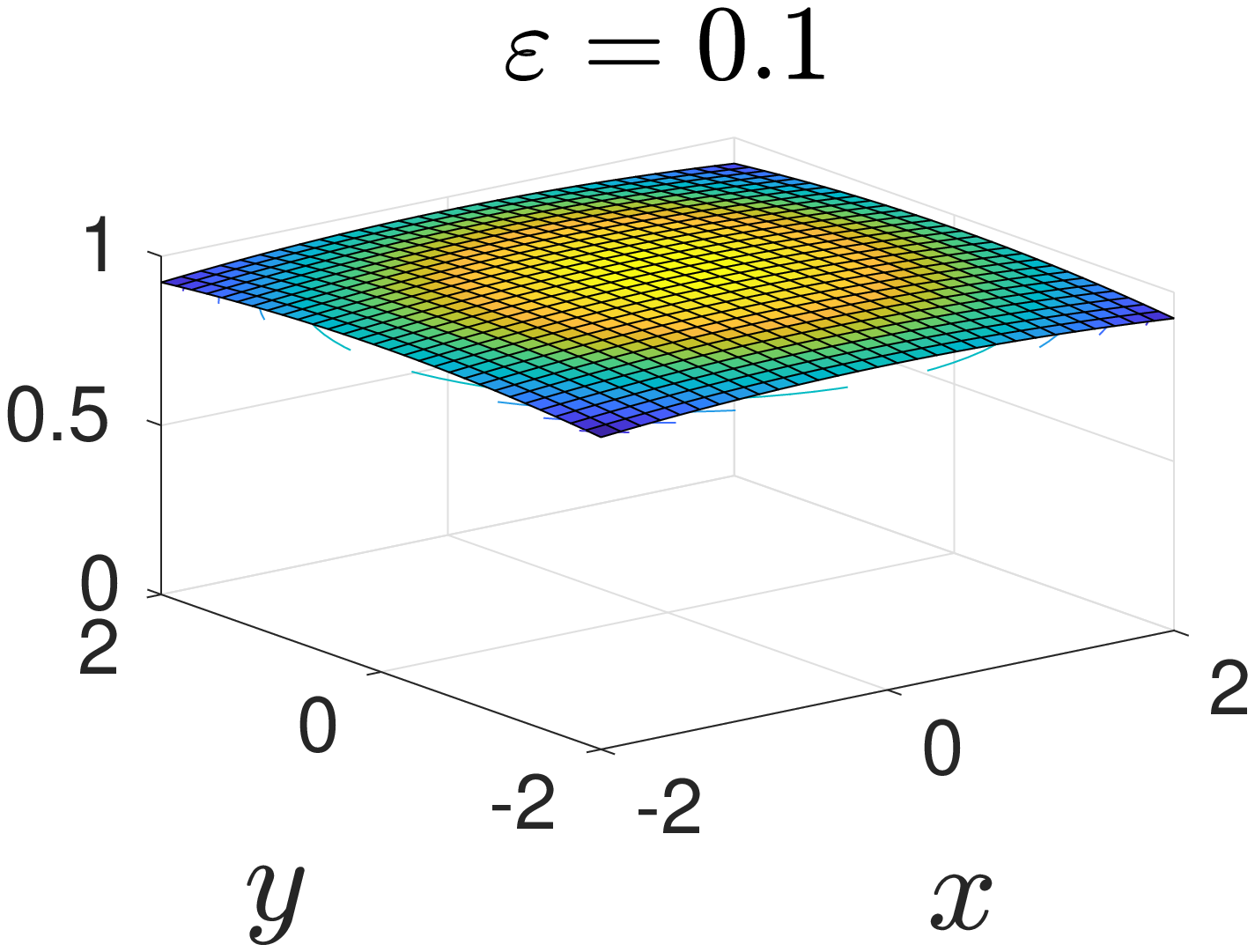}
   \end{center}
  \vspace{-.5cm} 
  \caption{RBF GA when  $\varepsilon=2,1,0.1$.}
  \label{Fig:rbf-ga_near_flat}
\end{figure}

The Fig. \ref{Fig:rbf-ga_near_flat} shows the cases for Gaussians RBFs, $\varphi(r,\varepsilon)=e^{-(r\varepsilon)^2}$ for 
differents $\varepsilon s$. For small and fixed shape parameter, the distance matrix of the RBF makes flatter (called near flat RBFs), 
so the linear combination of these kind of RBFs for interpolation becomes almost linear dependant like the case $\varepsilon = 0.1$ in the 
Fig. \ref{Fig:rbf-ga_near_flat}.

To avoid the ill-conditioning presented in the local interpolation of $u$ and also in the interpolation of the density term $b$, 
we change the basis as it was developed in the RBF-QR method in 2-dimensions \cite{fornberg_larsson_flyer_11}.
The new expansion of the Gaussian RBF takes the form:
\begin{eqnarray}
\phi(\|{\bf x}-{\bf x}_j\|) &=& \sum_{k=0}^\infty \sum_{l=0}^{\lfloor k/2 \rfloor} d_{k,l} \; c_{k,j}({\bf x}_j) \; C_{k,l}({\bf x})
  \label{rbf_expansion_1}\\
   &+& \sum_{k=0}^\infty \sum_{l=1-p}^{\lfloor k/2 \rfloor} d_{k,l} \; s_{k,l}({\bf x}_j) \; S_{k,l}({\bf x}),
  \label{rbf_expansion_2}
\end{eqnarray}
where $p=0$ if $k$ even and $p=1$ if $k$ odd. The scale factors $d_{k,l}$ is $O(\varepsilon^{2k})$ are
\begin{equation}
 d_{k,l} = \frac{\varepsilon^{2k}}{2^{k-2l-1} \left(\frac{k+2l+p}{2}\right)! \left(\frac{k-2l-p}{2}\right)!}, 
\end{equation}
and the coefficients $c_{k,l}, s_{k,l}$ are $O(1)$ given by
\begin{eqnarray}
  c_{k,l}({\bf x}_j) = b_{2l+p} \; t_{k-2l} \; e^{-\varepsilon^2 r^2_j} \; \cos((2l+p)\theta_j) \; _1F_2(\alpha_{k,l},\beta_{k,l},\varepsilon^4 r^2_j),\\
  s_{k,l}({\bf x}_j) = b_{2l+p} \; t_{k-2l} \; e^{-\varepsilon^2 r^2_j} \; \sin((2l+p)\theta_j) \; _1F_2(\alpha_{k,l},\beta_{k,l},\varepsilon^4 r^2_j),
  \label{expansion_coefficients}
\end{eqnarray}
where $b_0=1$, $b_k=2,\forall k>0$, $t_0=1/2$, $t_k=1,\forall k>0$, $_1F_2$ is the hypergeometric function with parameters 
$\alpha=\frac{k-2l+p+1}{2}$ and $\beta = \left[k-2l+1,\frac{k+2l+p+2}{2}\right]$ being $(r_j,\theta_j)$ is the polar coordinates 
location of the node ${\bf x}_j$.

The expansion functions $C_{k,l}$ and $S_{k,l}$ in (\ref{rbf_expansion_1}) and (\ref{rbf_expansion_2}) are given by
\begin{eqnarray}
  C_{k,l}({\bf x}) &=& e^{-\varepsilon^2 r^2} \; r^{2l} \; T_{k-2l}(r) \; \cos((2l+p)\theta),
  \label{expansion_functions_cos}\\
  S_{k,l}({\bf x}) &=& e^{-\varepsilon^2 r^2} \; r^{2l} \; T_{k-2l}(r) \; \sin((2l+p)\theta), \;\; 2l+p\neq 0
  \label{expansion_functions_sin}
\end{eqnarray}
where $\{T_n(r)\}$ are the Chebyshev polynomials.

\begin{figure}[!ht]
    \includegraphics[scale=0.31]{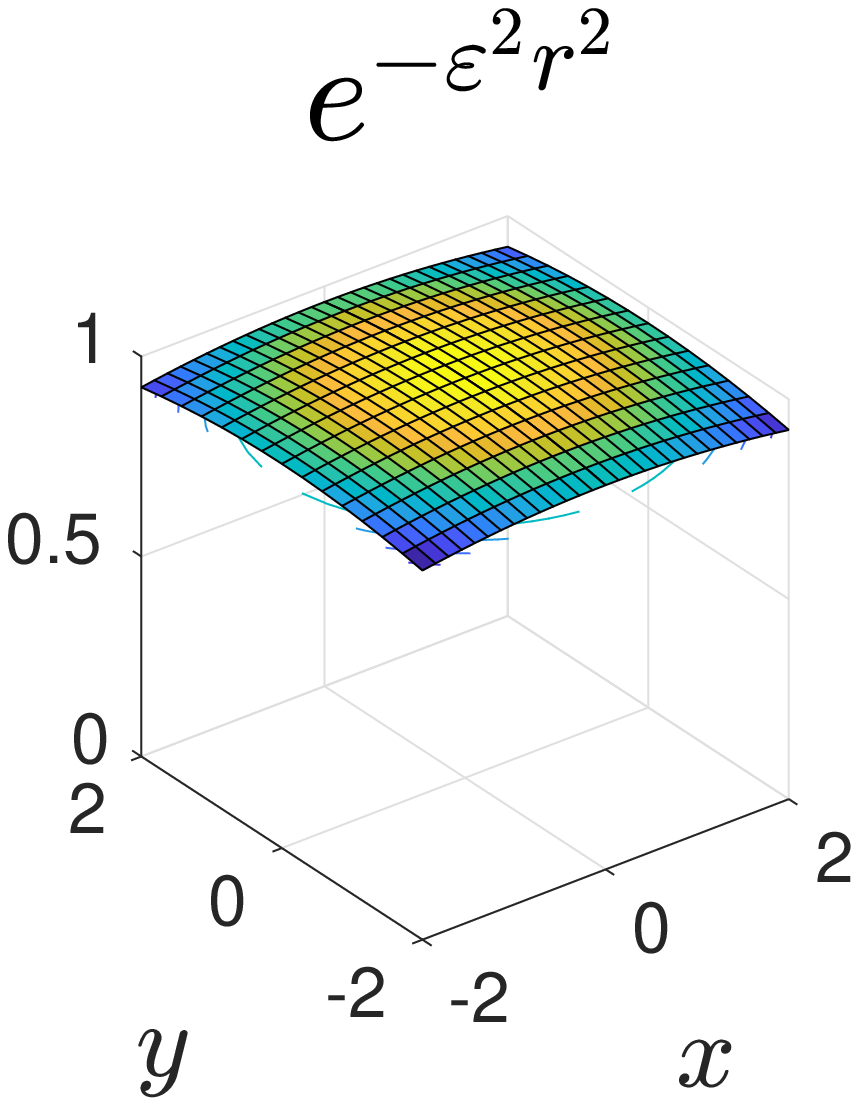}\\
    \includegraphics[scale=0.31]{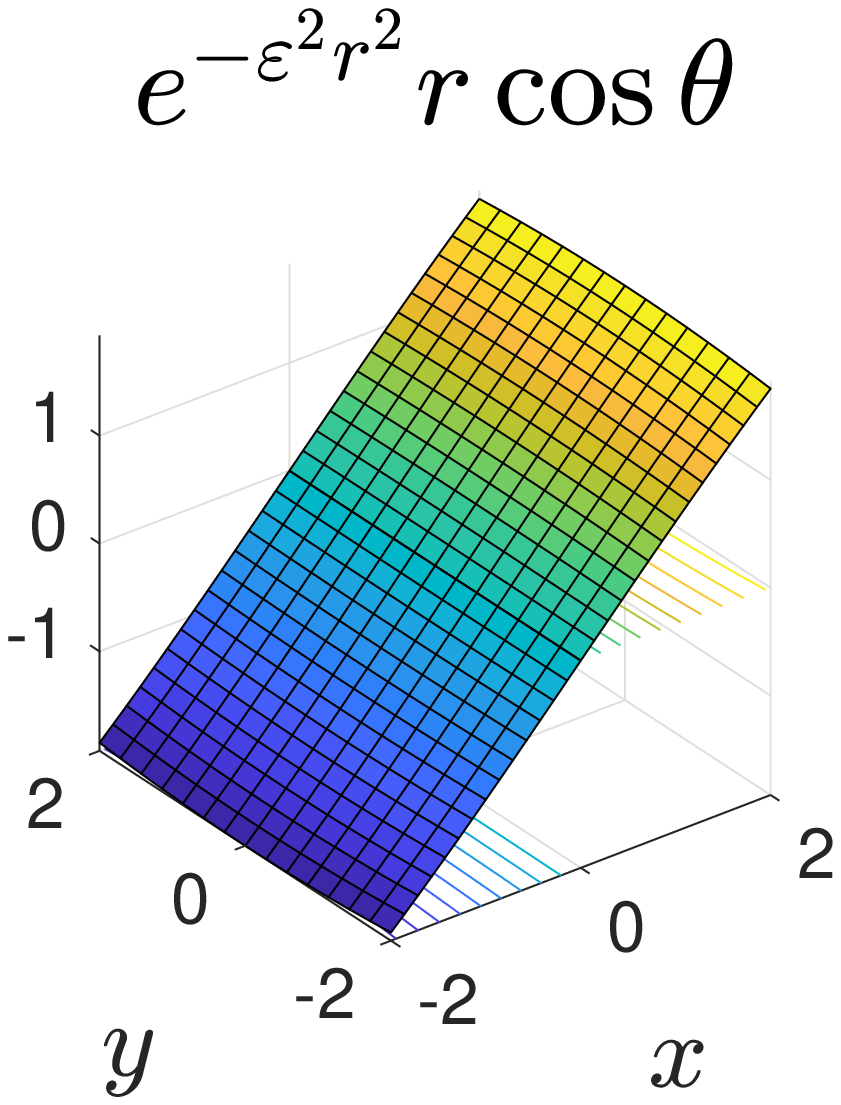}
    \hspace{.5cm} 
    \includegraphics[scale=0.31]{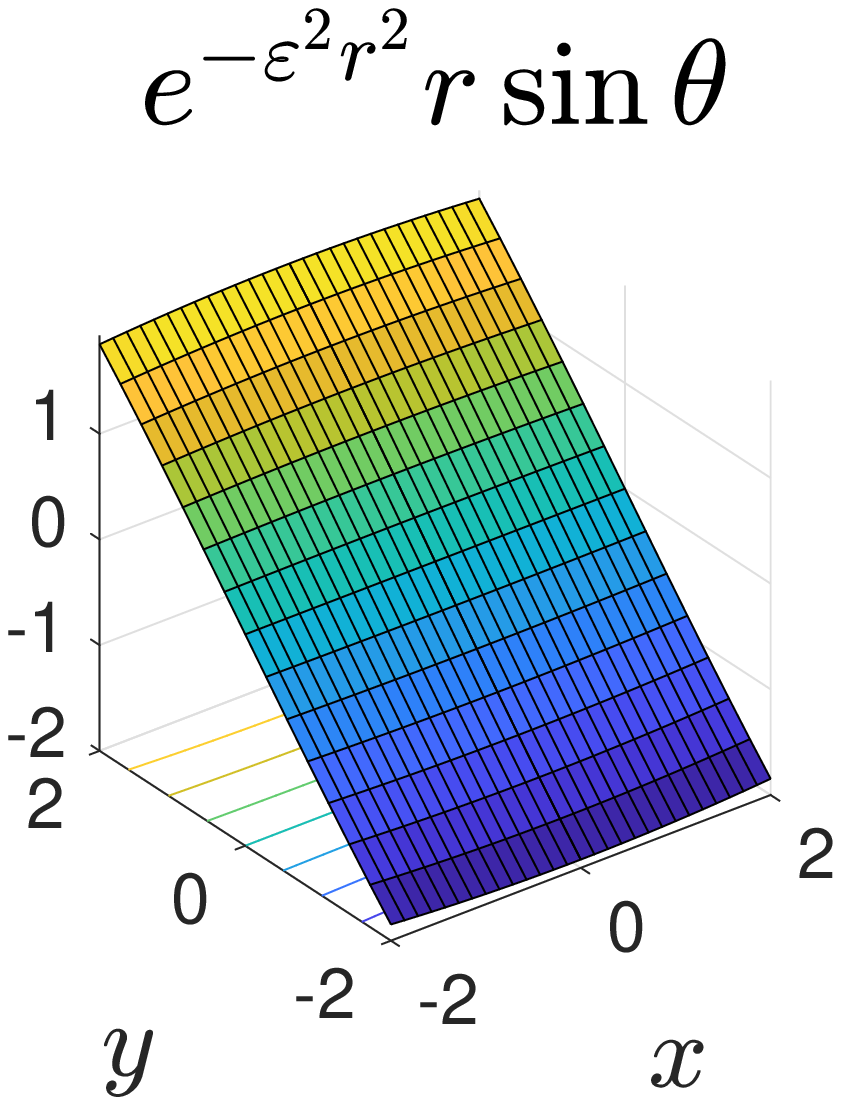}\\
    \includegraphics[scale=0.31]{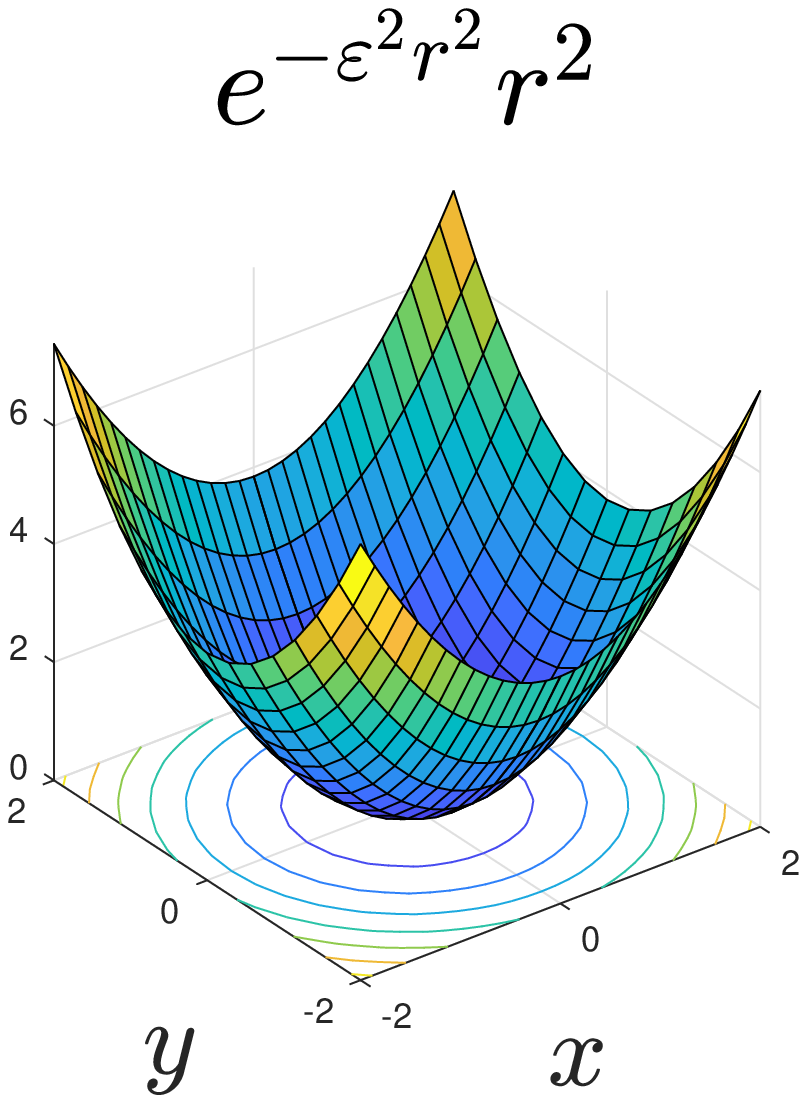}
    \hspace{.5cm} 
    \includegraphics[scale=0.31]{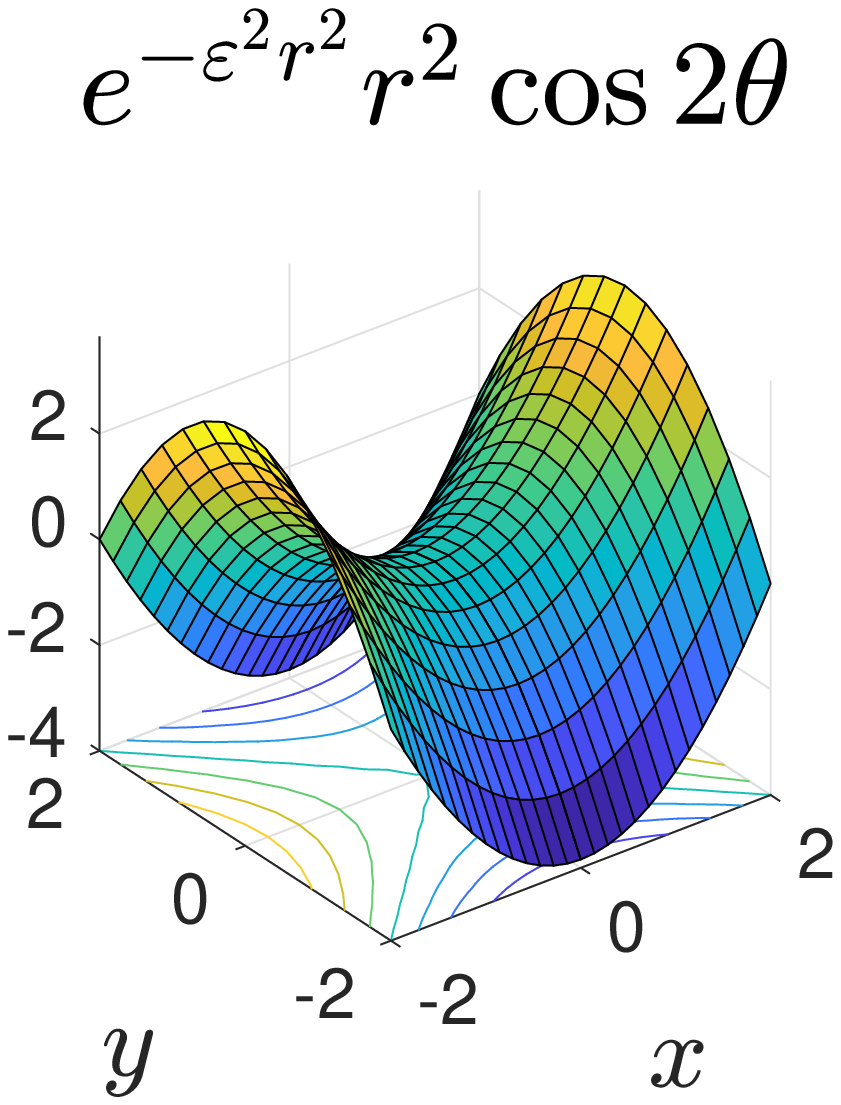}
    \hspace{.5cm} 
    \includegraphics[scale=0.31]{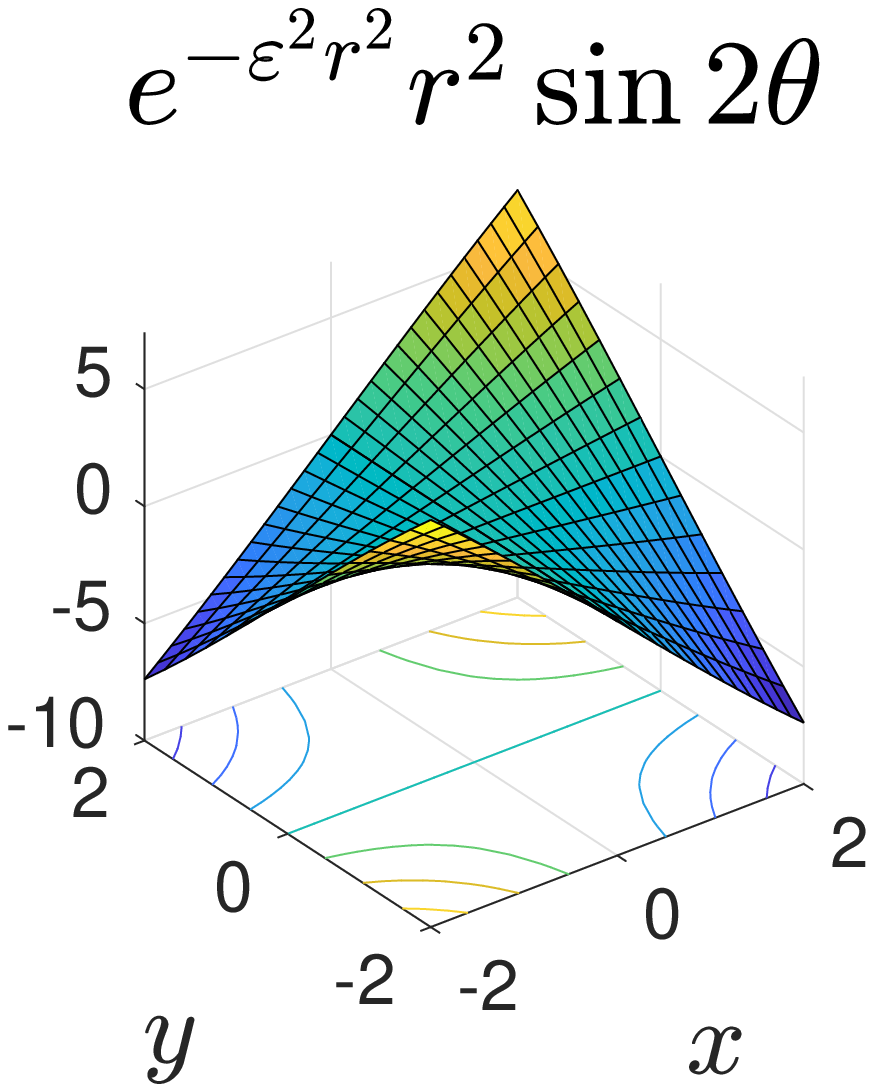}\\
    \includegraphics[scale=0.31]{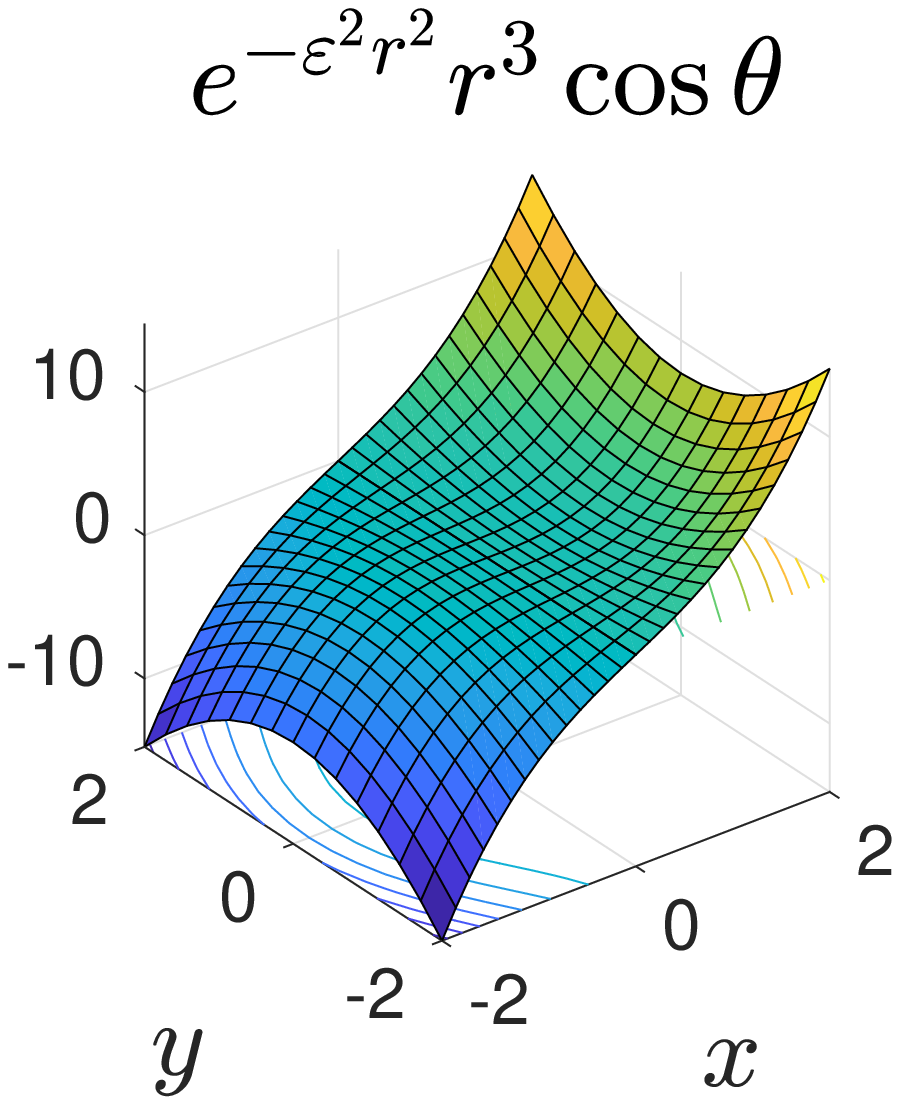}
    \hspace{.5cm} 
    \includegraphics[scale=0.31]{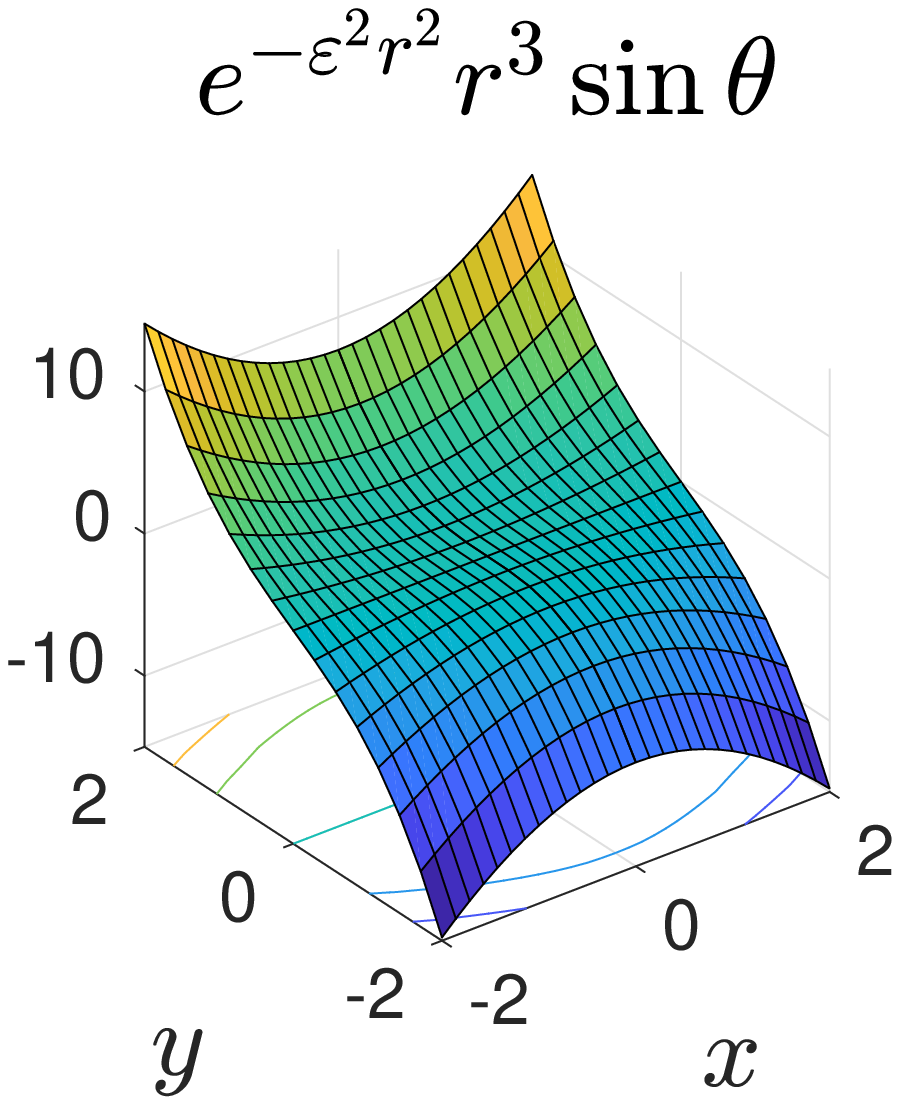}
    \hspace{.5cm} 
    \includegraphics[scale=0.31]{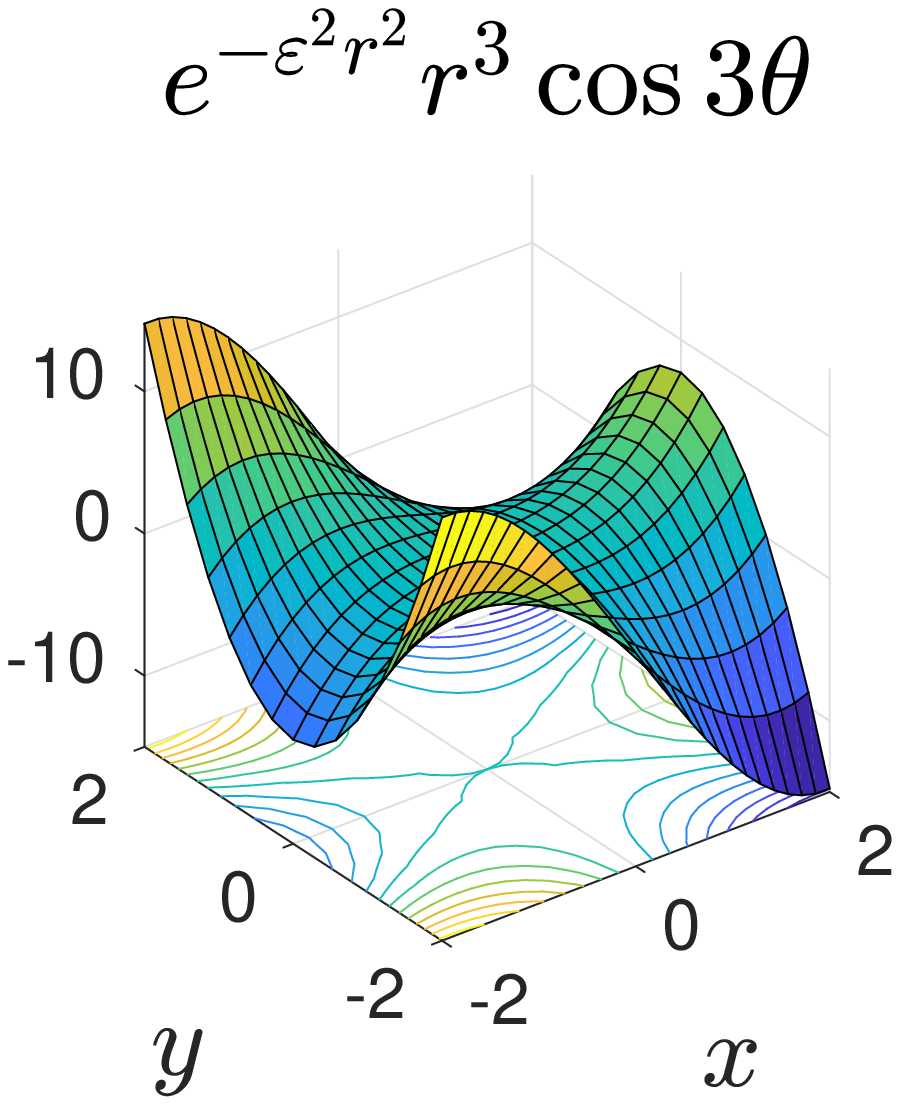}
    \hspace{.5cm} 
    \includegraphics[scale=0.31]{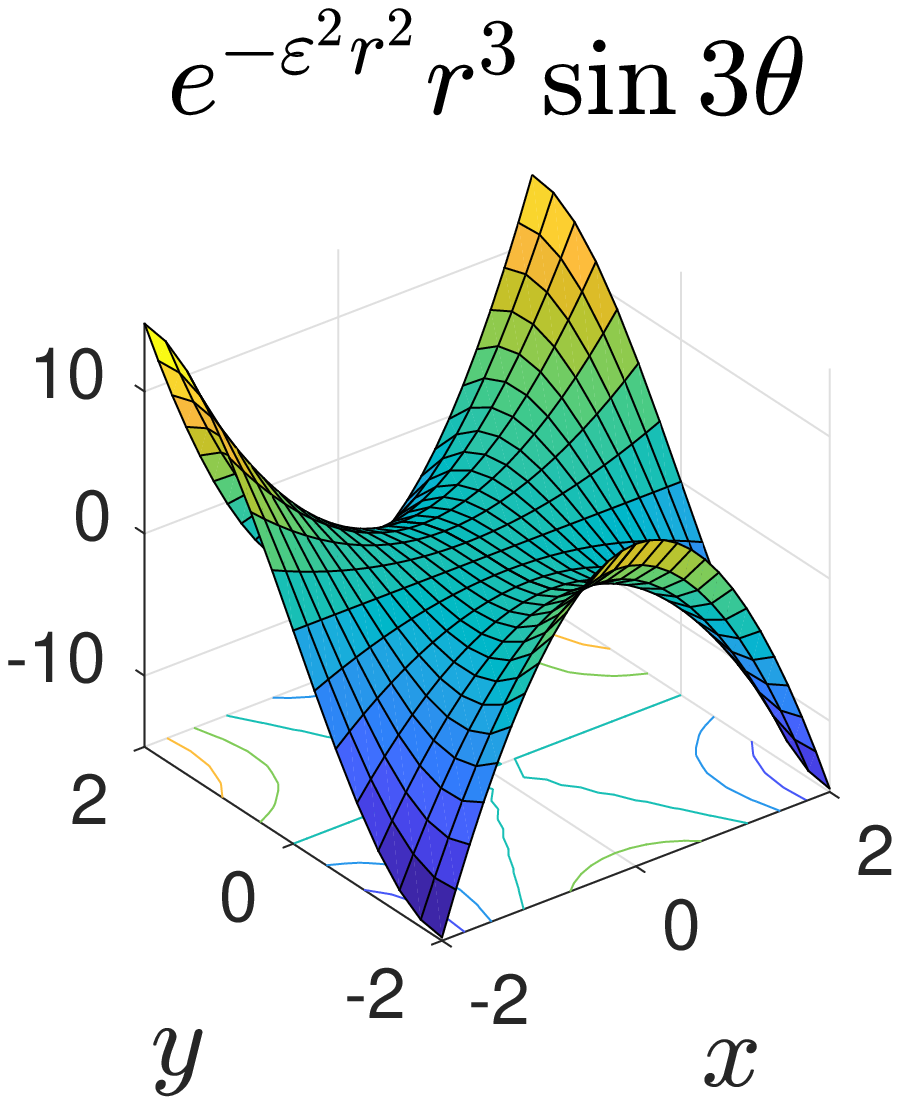}
  \vspace{-.5cm} 
  \caption{Some expansion functions from (\ref{expansion_functions_cos}) and (\ref{expansion_functions_sin}) 
  used to represent RBF GA for $\varepsilon = 0.1$.}
  \label{Fig:expansion_rbf-qr_1}
\end{figure}

Fig. \ref{Fig:expansion_rbf-qr_1} shows the four level of expansion functions used to generated $C_{k,l}({\bf x})$ and $S_{k,l}({\bf x})$ 
for small shape parameter $\varepsilon =0.1$ wich seems to be clearly linear independant. As $\varepsilon \rightarrow 0$, the term 
$e^{-\varepsilon^2 r^2} \rightarrow 1$ and the basis tends to the basis $\{1,r,r^r,r^3,\dots\}$ which is again an ill-conditioned basis because high powers of $r$ 
tend to be nearly dependent. The introduction of the Chebyshev polynomials in Eqs. (\ref{expansion_functions_cos}) and (\ref{expansion_functions_sin}) 
instead of monomials improve this situation. For more details see \cite{fornberg_larsson_flyer_11}.

In matrix form, each Gaussian RBF basis has the infinitely expansion
 \begin{eqnarray}
 \left[
 \begin{array}{c}
 \varphi_1(\mathbf{x})\\
 \varphi_2(\mathbf{x})\\
 \vdots\\
 \varphi_n(\mathbf{x}) 
 \end{array}
 \right]
  = \left[\begin{array}{ccc}
  . & . & \ldots\\
  . &  c_{k,l} & \ldots\\
  . &  s_{k,l} & \ldots\\
  . & . & \ldots\\
 \end{array}
 \right]\hspace{-.2cm}
  \left[\begin{array}{cccc}
  \ddots  &         &         & \\
          & d_{k,l} &         & \\
          &         & \ddots  & 
 \end{array}\right] \hspace{-.2cm}
  \left[\begin{array}{c}
 C_{0,0}(\mathbf{x}) \\
 C_{1,0}(\mathbf{x}) \\
 S_{1,0}(\mathbf{x}) \\
 C_{2,0}(\mathbf{x}) \\
 C_{2,1}(\mathbf{x}) \\
 S_{2,1}(\mathbf{x}) \\
\vdots
\end{array}
\right].
\end{eqnarray}

If we rename the functions $C_{k,l}$ and $S_{k,l}$ as $V_k$, the coefficients $c_{k,l}$ and 
$s_{j,m}$ as $\tilde{c}_k$, $d_{j,m}$ as $\tilde{d}_k$ and truncating for some $m$, we get the new basis 
approximation in matrix form
\begin{eqnarray}
 \left[
 \begin{array}{c}
 \varphi_1(\mathbf{x})\\
 \varphi_2(\mathbf{x})\\
 \vdots\\
 \varphi_n(\mathbf{x}) 
 \end{array}
 \right]
 \approx
 \left[\begin{array}{ccc}
  \tilde{c}_1(\mathbf{x}_1) &  \ldots  & \tilde{c}_m(\mathbf{x}_1)\\
  \vdots &  \ddots  & \vdots\\
  \tilde{c}_1(\mathbf{x}_n) & \ldots & \tilde{c}_m(\mathbf{x}_n) \\
 \end{array}
 \right]\hspace{-.2cm}
  \left[\begin{array}{ccc}
  \tilde{d}_1     &         &         \\
          & \ddots  &         \\
          &         & \tilde{d}_m  \\ 
 \end{array}\right] \hspace{-.2cm}
  \left[\begin{array}{c}
 V_1(\mathbf{x}) \\
 \vdots \\
 V_m(\mathbf{x}) \\
\end{array}
\right].
\end{eqnarray}

QR-factorizing the matrix $\mathbf{C}$, we obtain
\begin{eqnarray}
  \left[
  \begin{array}{c}
  \varphi_1(\mathbf{x})\\
  \varphi_2(\mathbf{x})\\
  \vdots\\
  \varphi_n(\mathbf{x})
  \end{array}
  \right]
  \approx \mathbf{C\;D\;}
  \left[
  \begin{array}{c}
  V_1(\mathbf{x}) \\
  V_2(\mathbf{x}) \\
  \vdots \\
  V_m(\mathbf{x}) \\
  \end{array}
  \right]
  &=& \mathbf{(QR) D\;V(x)}\\
  &=& \mathbf{Q} \left[ \; \mathbf{R}_1 \mathbf{D}_1 \; | \; \mathbf{R}_2 \mathbf{D}_2 \; \right]  \mathbf{V}(\mathbf{x})
\end{eqnarray}
where $\mathbf{C}$ is a $n\times m$ rectangular matrix where the elements of the coefficients are $O(1)$, 
$\mathbf{D}$ is a $m \times m$ diagonal matrix with the scaling coefficients $\tilde{d}_k$ proportional to $\varepsilon^{2k}$ for the integer $k\geq 0$, 
$\mathbf{Q}$ the $n \times n$ orthonormal matrix and $\mathbf{R}$ the $n\times n$ upper-triangular from QR algorithm. 

The matrix $\mathbf{R}$ is partitioned as $\left[ \mathbf{R}_1 \;|\; \mathbf{R}_2 \; \right]$ where $\mathbf{R}_1$ the upper triangular that contains 
the $n$ first columns of $\mathbf{R}$ and $\mathbf{R}_2$ is a matrix block $n\times m$. The scaling matrix $D$ is partitioned correspondingly with 
$\mathbf{D}_1$ a $n \times n$ diagonal block and $\mathbf{R}_2$ of size $(m-n)\times (m-n)$.

The vector function $\mathbf{V}(\mathbf{x})$ has components the functions $V_k\mathbf{(x)}$ combination of monomials, Chebyshev polynomials and trigonometric functions.

The new basis $\{\psi_j(\mathbf{x})\}$ is given by:
\begin{eqnarray}
 \left[
 \begin{array}{c}
 \psi_1(\mathbf{x})\\
 \psi_2(\mathbf{x})\\
 \vdots\\
 \psi_n(\mathbf{x}) 
 \end{array}
 \right]
 = \mathbf{D}_1^{-1} \mathbf{R}_1^{-1} \mathbf{Q}^T \left[
 \begin{array}{c}
 \varphi_1(\mathbf{x})\\
 \varphi_2(\mathbf{x})\\
 \vdots\\
 \varphi_n(\mathbf{x}) 
 \end{array}
 \right]
 \approx \left[ \; \mathbf{I}_n \; | \; \mathbf{\tilde{R}} \; \right] \mathbf{V(x)}
\label{rbf-qr_decomposition}
\end{eqnarray}
with $\mathbf{I}_n$ the identity matrix of size $n\times n$ and 
$\mathbf{\tilde{R}} = \mathbf{D}_1^{-1} \mathbf{R}_1^{-1} \mathbf{R}_2 \mathbf{D}_2$ the correction matrix 
that contains nonnegative powers of $\varepsilon$ due to the ordering of the scaling coefficients that generate 
local interpolation matrices better conditioned.

Now, for introducing the RBF-QR method into the integral formulations presented in Section \ref{sec:formulation}, the field variable $u$ is calculated numerically in the new basis over each integral subregions $\Omega_i$ and its corresponding boundary $\partial \Omega_i$, 

\begin{equation}
 u\left(\mathbf{x}\right) = \sum^{n}_{k=1}\gamma_k \psi_k\left(\mathbf{x}\right)
 \label{local_u_rbfqr}
\end{equation}
where the new coefficients $\{\gamma_k\}_{k=1}^n$ are to be determinated. Also the non-homogeneous term $b$ is interpolated in the new basis as in Eq. (\ref{Eq:int_term_no_hom})

\begin{equation}
 b \approx \sum_{k=1}^{n} \lambda_k  \psi_k\left(\mathbf{x}\right)
 \label{Eq:int_term_no_hom_new_basis}
\end{equation}
with the new coefficients $\{\lambda_k\}_{k=1}^n$ to be determinated.

In the same way that in the integral method presented above, we have an interpolation matrix for $u$ (denoted as $\mathbf{B}$) and another for $b$ (denoted as $\mathbf{\tilde{B}}$). These matrices can be obtained using the new basis.

For the internal stencils it is well known that they are equals and can be computed in a direct way, $(\mathbf{B})_{ik}=\psi_k(\mathbf{x}_i)$ for $i,k=1,\dots,n$, applying the transpose relation to Eq. (\ref{rbf-qr_decomposition}) at each center node $\mathbf{x}_i$ to get the new matrix

\begin{eqnarray}
  \mathbf{B_\psi} = \mathbf{V}
 \left[\begin{array}{c}
 \mathbf{I}_n \\
 \mathbf{\tilde{R}}^T \\
 \end{array}
 \right]. 
 \label{rbf-qr_collocation_matrix}
\end{eqnarray}
where the matrix $\mathbf{V}$ has elements $v_{ik}=V_k(\mathbf{x}_i)$.

It is well known that for positive definite Gaussians RBFs the interpolation matrix is always nonsingular for distinct nodes and $\varepsilon >0$ as was 
demonstrated by C.A. Micchelli in \cite{micchelli_86}. In this case the matrix $\mathbf{B_\psi}$ is nonsingular since the change of basis is well defined Eq. (\ref{rbf-qr_decomposition}) and as it 
was also discussed in \cite{larsson_lehto_heryudono_fornberg_13}.

When we have Neumann boundary conditions, we need to calculate the partial derivative $ \frac{\partial \psi_k}{\partial n}\left(\mathbf{x}\right)$ at some nodes 
$\{\mathbf{x}_i\}_{i=n_i+1}^n$ of the new basis functions. For this, we need to observate that functions  $\{\psi_k\}$ depend linearly on the 
expansion functions $\{V_k\}$. So, from Eq. (\ref{rbf-qr_decomposition}) it is possibly to calculate numerically the action of a boundary linear 
operator $\mathcal B$ on this basis as

\begin{eqnarray}
 \left[
 \begin{array}{c}
 \mathcal B\psi_1(\mathbf{x})\\
 \mathcal B\psi_2(\mathbf{x})\\
 \vdots\\
 \mathcal B\psi_{n}(\mathbf{x}) 
 \end{array}
 \right]
 \approx \left[ \; \mathbf{I}_n \; | \; \mathbf{\tilde{R}} \; \right] 
 \left[\begin{array}{c}
 \mathcal BV_1(\mathbf{x}) \\
 \mathcal BV_2(\mathbf{x}) \\
 \vdots \\
 \mathcal BV_m(\mathbf{x})
\end{array}
\right]=\left[ \; \mathbf{I}_n \; | \; \mathbf{\tilde{R}} \; \right] \mathbf{W(x)}.
\label{rbf-qr_diff_operator}
\end{eqnarray}
where the vector function $\mathbf{W(x)}$ has components $\mathcal B V_k(\mathbf{x})$,  $k=1,\dots,m$.

So, the local matrix interpolation $ \mathbf{B}$ that arises for the boundary stencil is formed with matrix blocks 

\begin{eqnarray}
  \mathbf{B} = 
 \left[\begin{array}{c}
 \mathbf{B_\psi} \\
 \mathbf{B_{\mathcal B\psi}} \\
 \end{array}
 \right]. 
 \label{rbf-qr_matrix}
\end{eqnarray}
where the block $\mathbf{B}_\psi$ has coefficients $(\mathbf{B}_\psi)_{ik}=\psi_k(\mathbf{x}_i)$ for $i=1,\dots,n_i$ and $k=1,\dots,n$ and the other block $\mathbf{B}_{\mathcal B\psi}$ has coefficients $(\mathbf{B}_{\mathcal B\psi})_{ik}=\mathcal B\psi_k(\mathbf{x}_i)$ for $i=n_i+1,\dots,n$ and $k=1,\dots,n$.
The Matlab's implementation used for for calculating these matrix blocks is the algorithm {\tt RBF\_QR\_diffmat\_2D} available from the first author's website in \cite{larsson_lehto_heryudono_fornberg_13}.

As in Eq. (\ref{alfa}) and (\ref{eq_beta_2}), the new local interpolation coefficients are now given by:

\begin{equation}
 \mathbf{\lambda}=\mathbf{B}^{-1}\mathbf{d}
\label{Eq:alfa_rbf-qr}
\end{equation}
and

\begin{equation}
 \beta = \mathbf{\widetilde{B}^{-1}}\left(\mathbf{f}_i + \mathbf{B}_b\mathbf{B}^{-1} \mathbf{d}\right)
\label{Eq:beta_rbf-qr}
\end{equation}
with matrix coefficients $(\mathbf{B}_b)_{kj}=\widetilde{b}\left(\psi_j\left(\mathbf{x}_k\right),\nabla\psi_j\left(\mathbf{x}_k\right)\right)$.

So taking $\mathbf{\xi} = \mathbf{\xi}_i$ as in the local approach in Eq. (\ref{Eq:DRM-meshless-i}) for each interior collocation point $\{\mathbf{\xi}_i\}$, we have

\begin{equation}
 u_i=  \sum^{n}_{k=1} l_{ik} \ \gamma_k  +  \sum^{n}_{k=1} \tilde{l}_{ik} \ \lambda_k + f_i
 \label{new_collocation_u}
\end{equation}
where the integrals 

\begin{eqnarray}
l_{ik}  =  \int_{\Gamma_{i}}Q\left(\mathbf{x},\mathbf{x}_i\right) \psi_k\left(\mathbf{x}\right) d\Gamma_i,\\ 
\label{hij_fornberg}
\tilde{l}_{ik}  =  \int_{\Omega_{i}} G\left(\mathbf{x},\mathbf{x}_i\right) \psi_k\left(\mathbf{x}\right) d\Omega_i
\label{fij_fornberg}
\end{eqnarray}
are calculated with the new basis $\{\psi_k\}$ instead of the Gaussian RBF basis $\{\varphi_j\}$, being $\Omega_{i}$ the local region of integration and $\Gamma_{i}$ the correspondingly boundary for each $\mathbf{\xi}_i$.
These line and volume integrals are calculated numerically using Gauss-Legendre method of cuadrature. 

Similarly to Eq. (\ref{local_int_eq-ii_zi}), here we obtained the following discretized form for the unknown field $u$ for each internal point

\begin{equation}
 u_i = f_i +\left(\mathbf{l}_i^T \mathbf{B}^{-1} +
 \mathbf{\widetilde{l}}_i^T \mathbf{\widetilde{B}}^{-1} \mathbf{B_b}
 \mathbf{B}^{-1}\right)\mathbf{d}.
\label{local_int_eq-ii_rbf-qr}
\end{equation}
with column vectors $\mathbf{l}_i=\left[\ldots,l_{ij},\ldots\right]^T$ and 
$\mathbf{\widetilde{l}}_i=[\ldots,\widetilde{l}_{ij},\ldots]^T$ that is solve in analogous algorithmic procedure as before.

This alternative to the LRDRM was called \emph{Local Integral RBF-QR Method (LIM RBF-QR)} which avoid finding an auxiliary particular solution for each element of the new basis $\{\psi_k\}$ stable for small shape parameters.

\section{Numerical results}
\label{sec:num-results}

In this section we explore with numerical results the stability of the RBF-QR method with integral methods and validate the formulation 
presented above. We consider Boundary Value Problems for different PDEs: three 2D Poisson Equations as follows, with mixed boundary 
conditions, with Dirichlet BC both cases in square domains and a third one over a circular domain with Dirichlet BC. Adiotionally, 
1D and 2D Convection-Diffusion Equations with mixed BCs.
Numerical results are presented for different domains and node sets distributions uniform, halton, quasi-uniform and scattered repel distribution.

The obtained numerical results were compared with the corresponding exact solutions when available. Equations for the errors 
presented in this work are: \emph{Absolute Maximum error ($L_{\infty}-error$  Eq. \ref{Eq:Linf}), L2 porcentual error 
($L_2-error \%$ Eq. \ref{Eq:L2error-porcentual}) and Root Mean Square error ($L_{RMS}-error$ Eq. \ref{Eq:RMS})}.
\begin{equation}
L_\infty-error = \max_{i=1,...,N} \left| u^{i}_{exact}-u^{i}_{apx}\right|
\label{Eq:Linf} 
\end{equation}
\begin{equation}
L_2-error \%=100\% \sqrt{\frac{\sum_{i=1}^{N}\left(u^{i}_{exact}-
u^{i}_{apx}\right)^2}{\sum_{i=1}^{N}\left(u^{i}_{exact}\right)^2}} 
\label{Eq:L2error-porcentual} 
\end{equation}
\begin{equation}
L_{RMS}-error = \sqrt{\frac{\sum_{i=1}^{N}\left(u^{i}_{exact}-
u^{i}_{apx}\right)^2}{N}} 
\label{Eq:RMS} 
\end{equation}
with $u^{i}_{exact}$ as the nodal values of the exact solution and $u^{i}_{apx}$ 
the corresponding values of the approximation. Comparisons with results in 
\cite{caruso_portapila_power_15,ooi_popov_12,ooi_popov_12_ComMech,bayona_flyer_fornberg_barnett_16} 
for the same equations are also reported.

\subsection{Domain discretizations}

The domains $\Omega$ considered in this paper were discretized using uniform and scattered nodes. 
For the 2D scattered nodes we used Halton \cite{halton_60} nodes, quasi-uniform nodes \cite{fornberg_flyer_15} and 
a repel algorithm presented in \cite{bayona_flyer_fornberg_barnett_16}. The construction of the first and second are based on 
deterministic method.

The 2-dimensional Halton nodes were created from the van der Corput sequences taking a number prime as its basis to generate 
well-spaced points from the interval $(0,1)$. To generate Halton points in $(0,1)^2$ each coordiante was generated with a different 
prime number. Then they were transformed linearly or traslated to a rectangular domain in $\mathbb{R}^2$. The boundary nodes 
were also generated as Halton nodes. The Matlab code used was {\tt halton} written by B. Fornberg available at Matlab 2017a version.

The quasi-uniform nodes were created using the Matlab code {\tt node\_placing}, the implementation of a fast generation 
algorithm for 2D meshfree PDE discretizations developed by B. Fornberg and N. Flyer in \cite{fornberg_flyer_15}. This is 
an advancing front type method that creates a node set from a varialbe density function in rectangular and irregular domains. 
Its start at some boundary advancing until the total domain is filled.

The repel algorithm to create scattered nodes was described in \cite{bayona_flyer_fornberg_barnett_16}. The idea is to embed the domain $\Omega$ 
(a circle in our experiment) into a rectangle in 2-D and discretized it using structured nodes discarding nodes lying within certain distance to boundary. 
Applied a random displacement to the interior nodes in the circle and also a displacement in the direction of the repulsion force
$\vec{F}(x,y)=\sum_{i=1}^n \frac{\vec{\mathbf{r}}_i}{\|\vec{\mathbf{r}}_i\|^3_2}$, with $\vec{\mathbf{r}}_i = \left(x-x_i,y-y_i\right)$, 
$n$ the number of closest nodes to the $i$-node $(x_i,y_i)$. The boundary $\partial \Omega$ is dicretized with equispaced nodes and fixed to their position. 

\begin{figure}[!ht]
\centering
  \begin{center}
    \includegraphics[scale=.375]{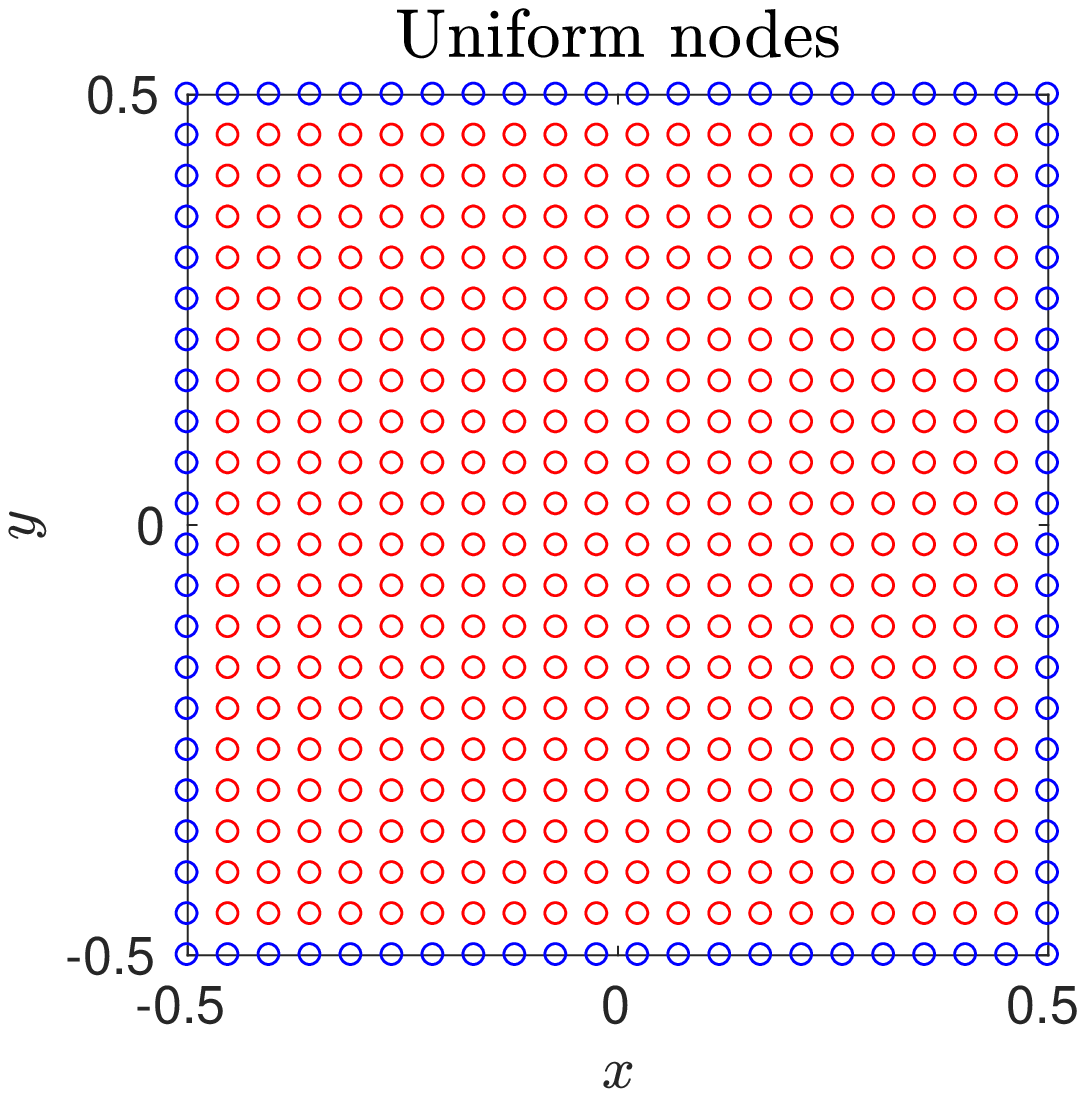}
    \hfill
    \includegraphics[scale=.375]{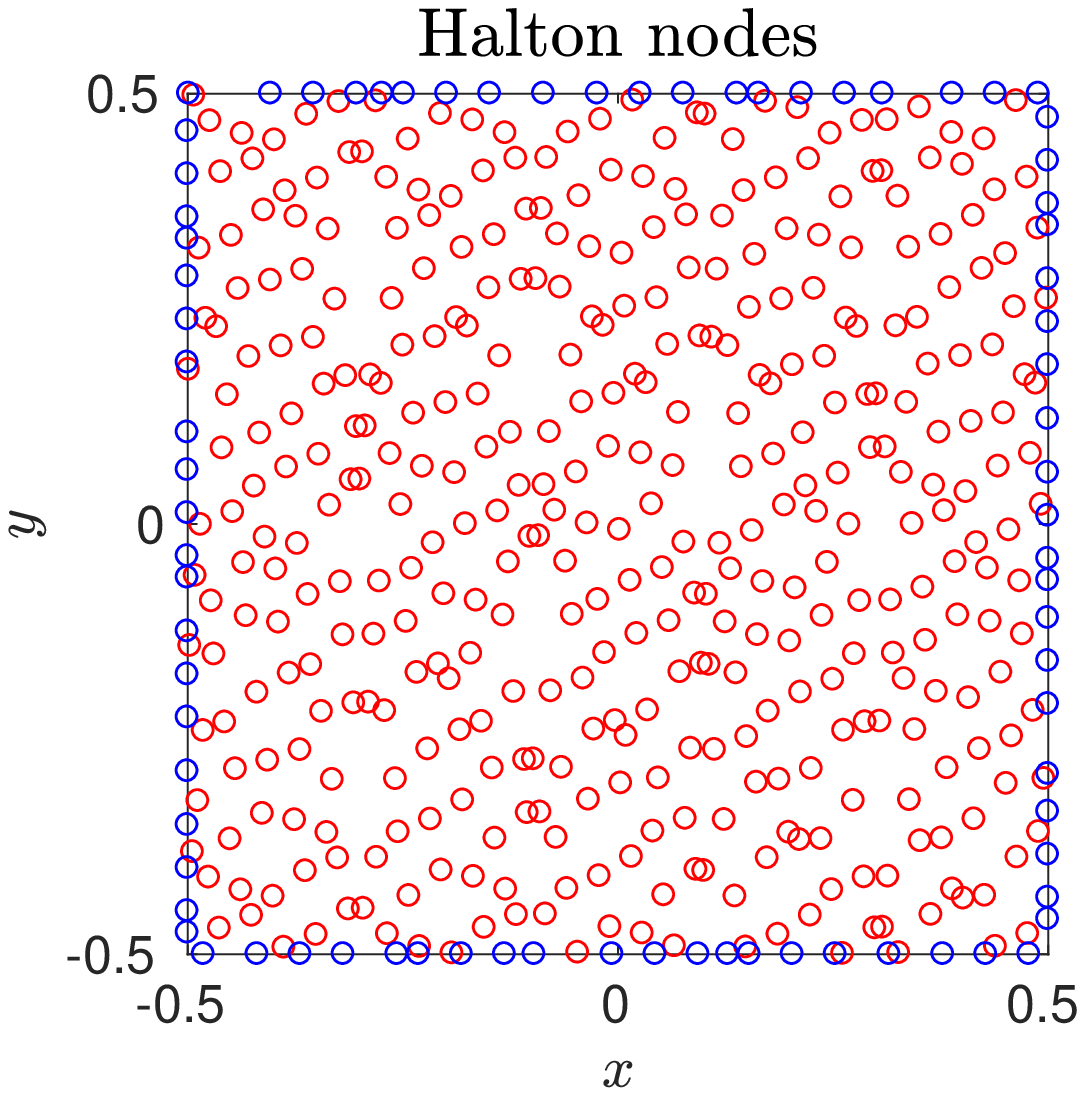}
    \hfill
    \includegraphics[scale=.375]{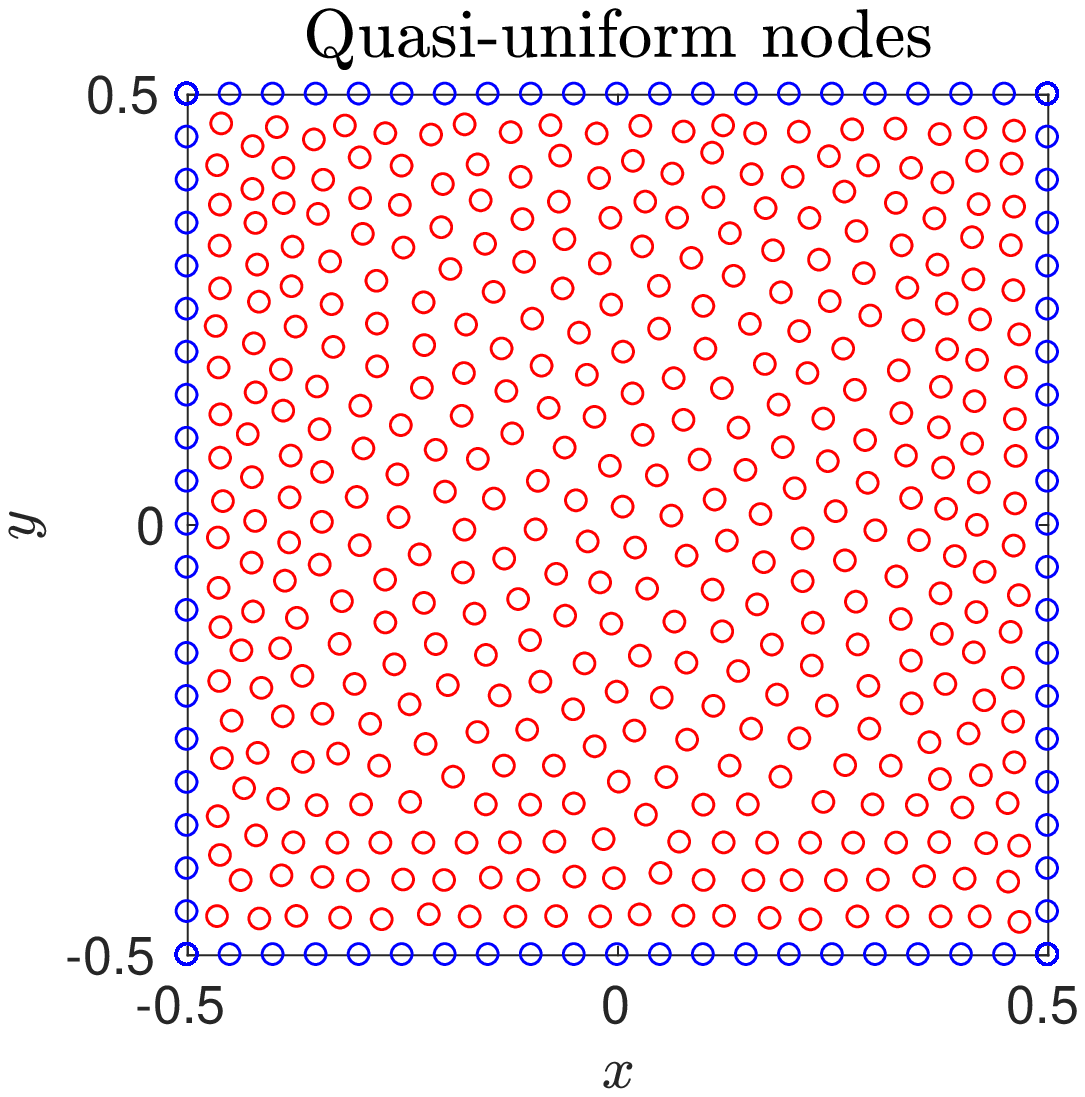}
   \end{center}
     \vspace{-.75cm}
  \caption{Discretized domain using uniform nodes. $N=400$ interior nodes, 84 boundary nodes (left). 
  Halton type nodes centers. $N=400$ interior poitns, 80 boundary nodes (center). Quasi-uniform node set 
  using $N=401$ interior centers, 76 boundary nodes (right)}
  \label{Fig:PEP5_node_sets}
\end{figure}

Uniform, Halton and quasi-uniform node distributions used in this paper for $N=400$ and $\Omega = \left[-\frac{1}{2},\frac{1}{2}\right]\times [-\frac{1}{2},\frac{1}{2}]$ 
are shown in Fig. \ref{Fig:PEP5_node_sets}. The node distributions with more centres are similar, althought denser.

For every domain we considered a band near the boundary whereif the center node is located inside this band, 
the local stencil $\Theta_i$ takes $n_i$ interior nodes and $n_b$ boundary nodes. 

All experiments were with a fixed number of nodes in the stencil except test problem \ref{PEP8} were we considered an 
increasing number of local stencils.

\subsection{Poisson's equation with mixed boundary conditions}
\label{PEP5}

Let us consider the following elliptic PDE, a Poisson's problem whereas the non-homogenous term is a product of trigonometric functions.
The governing equation in the square domain $\Omega=\left[-\frac{1}{2},\frac{1}{2}\right]^2$ is:
\begin{equation}
\Delta u\left(x_1,x_2\right)  =  \frac{5}{4}\sin(\pi x_1)\cos\left(\frac{\pi x_2}{2}\right), 
\quad \quad (x_1,x_2)\in \Omega,
\label{E:NR_Prob_1_E}
\end{equation}
and the mixed boundary conditions are:
\begin{equation}
(BCs)\left\{
\begin{array}{lcl}
u\left(-0.5,x_2\right) & = & -\frac{\sqrt{2}}{2}\cos(\pi x_2),\\
u\left(0.5,x_2\right) & = & \frac{\sqrt{2}}{2}\cos(\pi x_2),\\
\dfrac{\partial u}{\partial x_2}\left(x_1,-0.5\right) & = &  \pi\sin\left(\frac{\pi}{2}x_1\right),\\
\dfrac{\partial u}{\partial x_2}\left(x_1,0.5\right)  & = & -\pi\sin\left(\frac{\pi}{2}x_1\right).\\
\end{array}
\right. 
\label{E:NR_Prob_1_BV}
\end{equation}
The analytical solution to this problem is given by $u\left(x_1,x_2\right)  =  \sin(\pi x_1)\cos\left(\frac{\pi x_2}{2}\right).$

For this first example we considered uniform, Halton and quasi-uniform node set distribution, starting from a coarse distribution with $N=400$ 
integration subregions, up to a denser number of subregions, $N=2500$.

The integral equation is applied at only one source point per subdomain located at its centre with each stencil subdomain having also the same point 
for the local interpolation. In the numerical results in this work, $N$ define the total number of collocation points interior of the problem domain 
which coincide with the number of integration subregions or Green's elements.
The number of the stencil size $\Theta_i$ for the local interpolation is fixed at $n=25$ points and the interpolation RBFs are Gaussians.

\begin{figure}[!ht]
\centering
  \begin{center}
    \includegraphics[scale=.75]{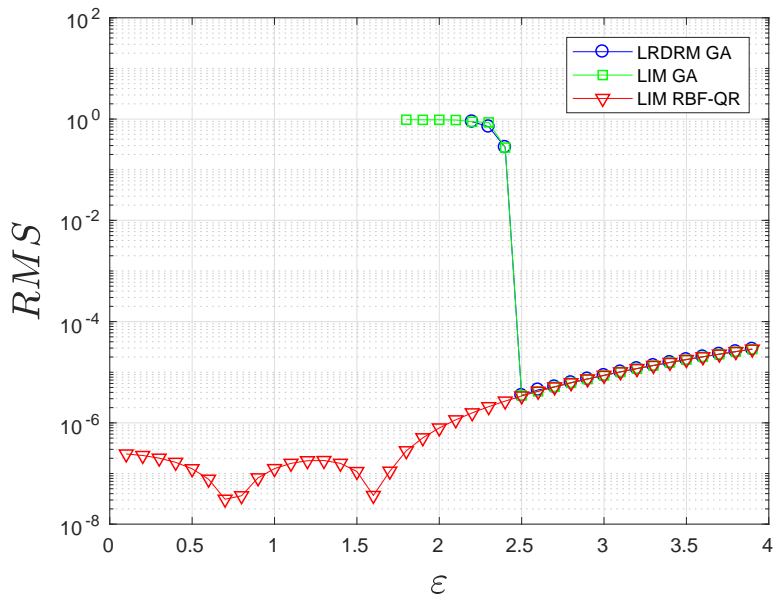}
  \end{center}
  \vspace{-.75cm}   
\caption{Comparison between LRDRM, LIM, LIM RBF-QR with Gaussians for $N=1600$ interior points}
  \label{Fig:PEP5_Errors_uni40x40}
\end{figure}

The obtained RMS error for $u$ is reported in Fig. \ref{Fig:PEP5_Errors_uni40x40}. This Figure  shows the RMS error as a function of 
the shape parameter $\varepsilon$ for three numerical methods. The fisrt one is the Localized Regular Dual Reciprocity Method (LRDRM) developed 
in \cite{caruso_portapila_power_15} using RBF GA for the local interpolation of $u$ and the non-homogeneous term $b$. The second one is a modification 
of this method where we replaced the use of the Reciprocity method implemented before by evaluating domain integrals with Gauss-Legendre quadrature 
and using RBF GA basis for the local interpolation. This way of calculating the integral equations was called Localized Integral Method (LIM GA). 
And finally, we do a further improvement for near flat RBFs using RBF-QR for the local interpolation for $u$ and also for the term $b$ with Gaussians RBFs.

For large $\varepsilon$, the three methods give the same results. The figure shows that the RMS of the LRDRM desestabilizes when the shape parameter 
is $\varepsilon = 2.5$ for the uniform node set distribution formed by $N=1600$ interior nodes arriving to $3.4971E-06$ for the LRDRM and $3.4848E-06$ for 
LIM GA. Using the RBF-QR method locally for small shape parameters estabilizes the RMS method significantly. The RMS in this case is $3.1063E-08$, 
2 orders of magnitude better. The uniform discretizations for $N=400,900,2500$ have the same behavior arriving to $1,8359E-07$, $5,8726E-08$ and $2.0682E-08$ 
respectively. The total comparison is shown in Table \ref{Table:PEP5_comparison_RMS} for the uniform distribution. 

\begin{table}[!ht]
\centering
\begin{scriptsize}
\begin{tabular}{|c|c c|c c|c c|c c|}
\hline
  $N$ & \multicolumn{2}{|c|}{LRDRM GA} & \multicolumn{2}{|c|}{LIM GA} & \multicolumn{2}{|c|}{LIM RBF-QR (1st dip)} & \multicolumn{2}{|c|}{LIM RBF-QR (2nd dip)}\\
      & $\epsilon$ & $RMS$ & $\epsilon$ & $RMS$ & $\epsilon$ & $RMS$ & $\epsilon$ & $RMS$\\
\hline
 $400$  & 1.60 & 1.1460E-06 & 1.60 & 1.1208E-06 & 1.60 & 1.1174E-06 & 0.80 & 1.8359E-07 \\ 
 $900$  & 1.90 & 1.6175E-06 & 1.90 & 1.6082E-06 & 1.60 & 1.5320E-07 & 0.80 & 5.8726E-08 \\
 $1600$ & 2.50 & 3.4971E-06 & 2.50 & 3.4849E-06 & 1.60 & 3.6854E-08 & 0.70 & 3.1063E-08 \\
 $2500$ & 2.90 & 4.1908E-06 & 2.90 & 4.3092E-06 & 1.60 & 2.3318E-08 & 0.70 & 2.0682E-08 \\
\hline
\end{tabular}
\caption{Poisson PDE with mixed BC - RMS - Uniform distribution}
\label{Table:PEP5_comparison_RMS}
\end{scriptsize}
\end{table}

\begin{figure}[!ht]
\centering
  \begin{center}
    \includegraphics[scale=.75]{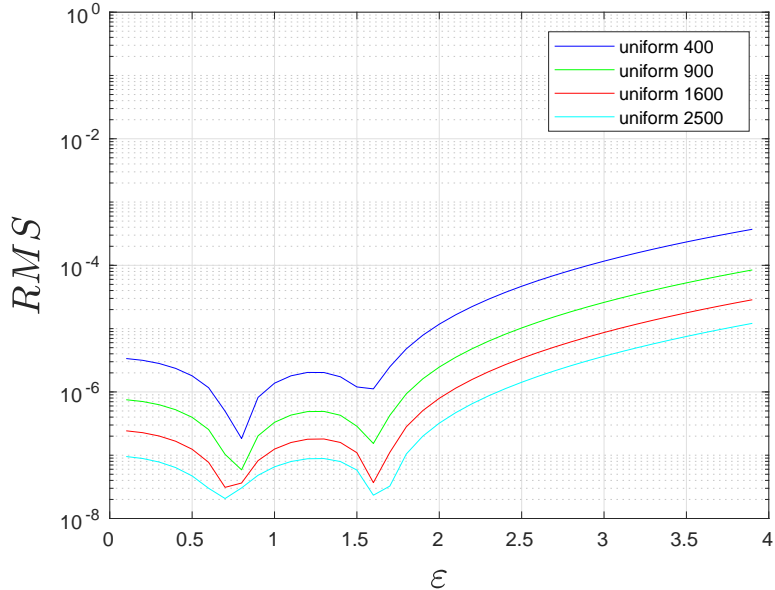}
  \end{center}
  \vspace{-.75cm}   
\caption{Runge phenomenon in the RMS error comparison between different uniform node sets for the Local Integral RBF-QR Method for $\varepsilon \in [0.1,4.0]$.}
  \label{Fig:PEP5_LRDint-Sint_rbf-qr_uniform_RMS}
\end{figure}

Fig. \ref{Fig:PEP5_LRDint-Sint_rbf-qr_uniform_RMS} shows the RMS versus the shape parameter for the different uniform distributions 
when $\varepsilon \rightarrow 0$. It is known that decreasing the shape parameter $\varepsilon$ produce more flat RBFs Gaussians which 
allow more accuracy in the numerical solution of the two systems for the local interpolation matrix $\mathbf{A}$ and $\mathbf{\tilde{A}}$. 
This ill conditioning dominates the error of the Local Integral Method. It is observes also that the error due to Runge phenomenon emerges. 
Such is the case in Fig. \ref{Fig:PEP5_LRDint-Sint_rbf-qr_uniform_RMS}, where we show for different node sets. This error reaches 
very low levels but for values near $\varepsilon=0.8$ and $\varepsilon=1.6$ increases. When increasing the number of total nodes, the diference 
between the first dip and the second becomes smaller (column 3rd and 4rd from Table). This is due to the Runge phenomenon studied in \cite{fornberg_zuev_07} 
since no particular mechanism was taken for controlling that. Note that for the uniform nodes distribution the location of these \emph{curve dips} 
not depends on the number of points $N$ since they appear for the same $\varepsilon$s. An advanced strategy for dealing with the Runge phenomenon 
is to use \emph{spatially varying} shape parameters like in \cite{fornberg_zuev_07}.\\

\begin{figure}[!ht]
 \centering
   \begin{center}
     \includegraphics[scale=.45]{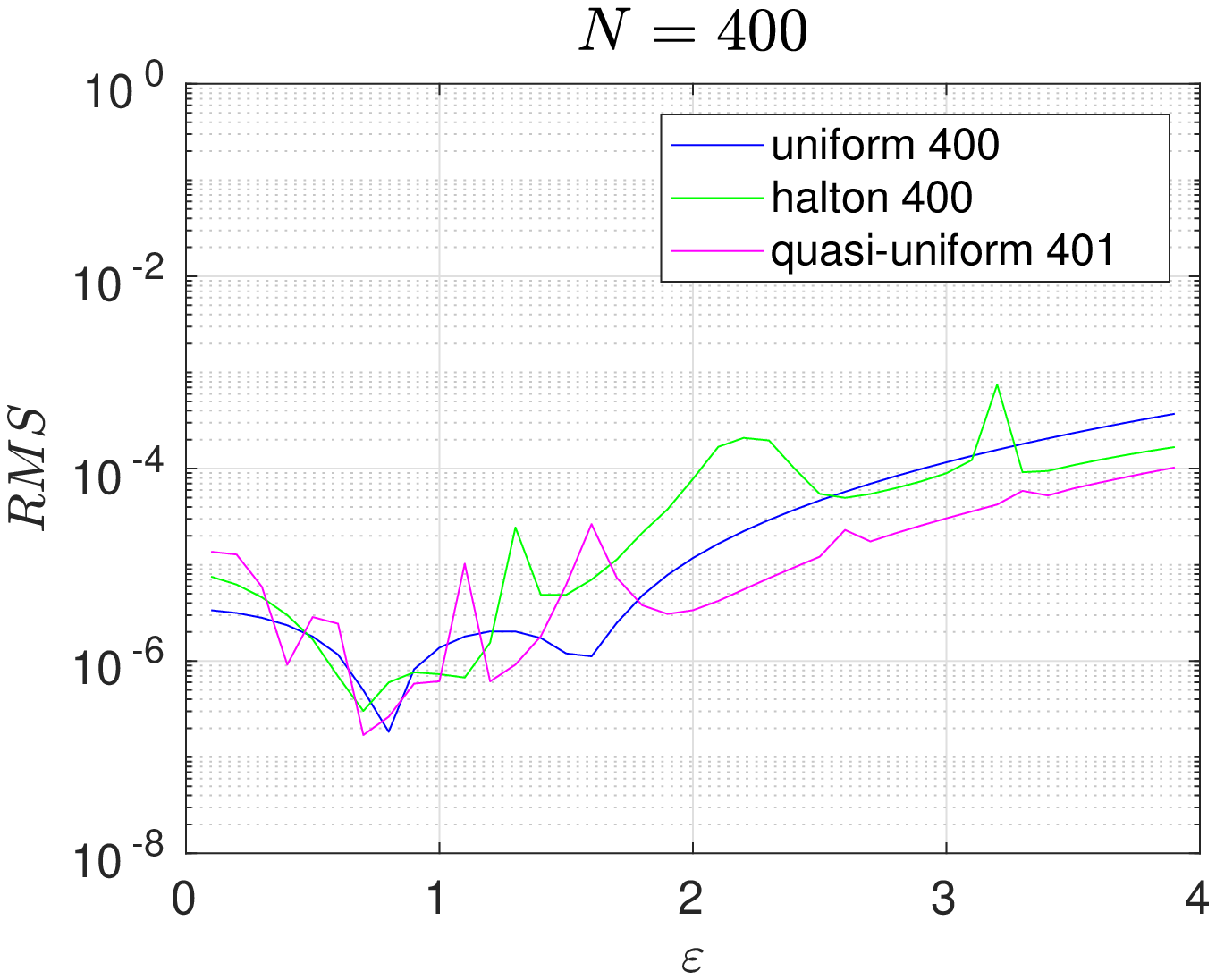}
     \hfill
     \includegraphics[scale=.45]{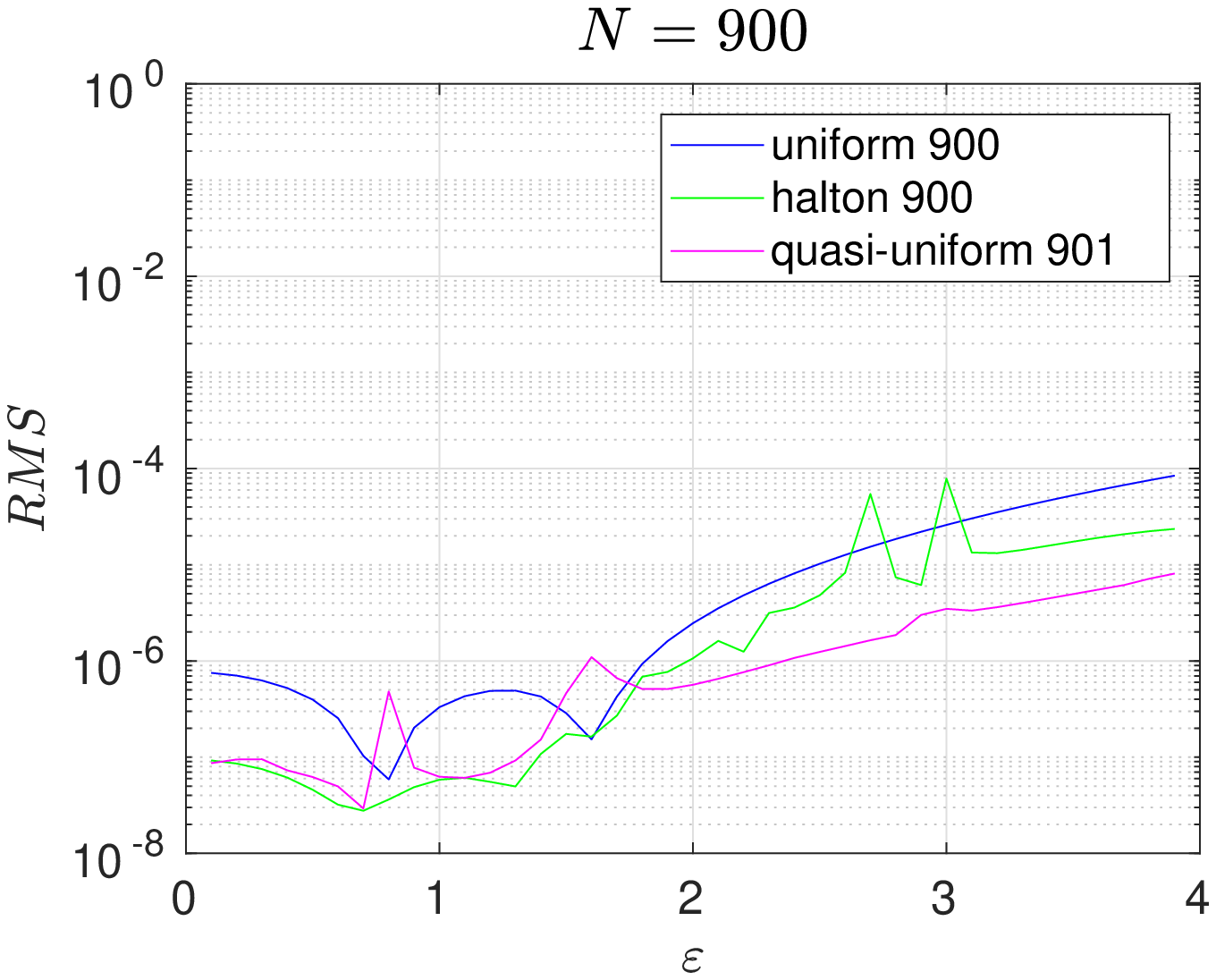}\\
     \vspace{.5cm}
     \includegraphics[scale=.45]{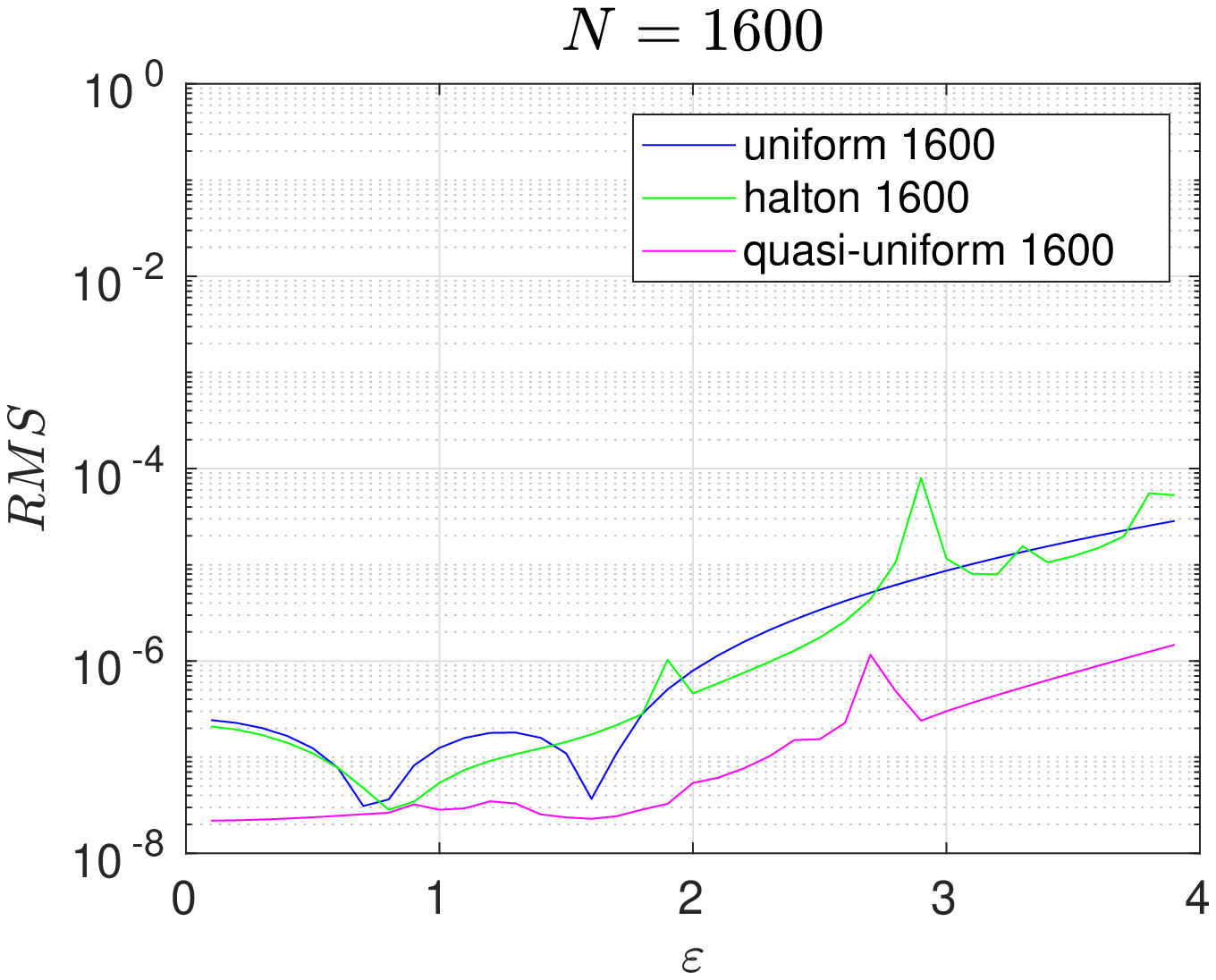}
     \hfill
     \includegraphics[scale=.45]{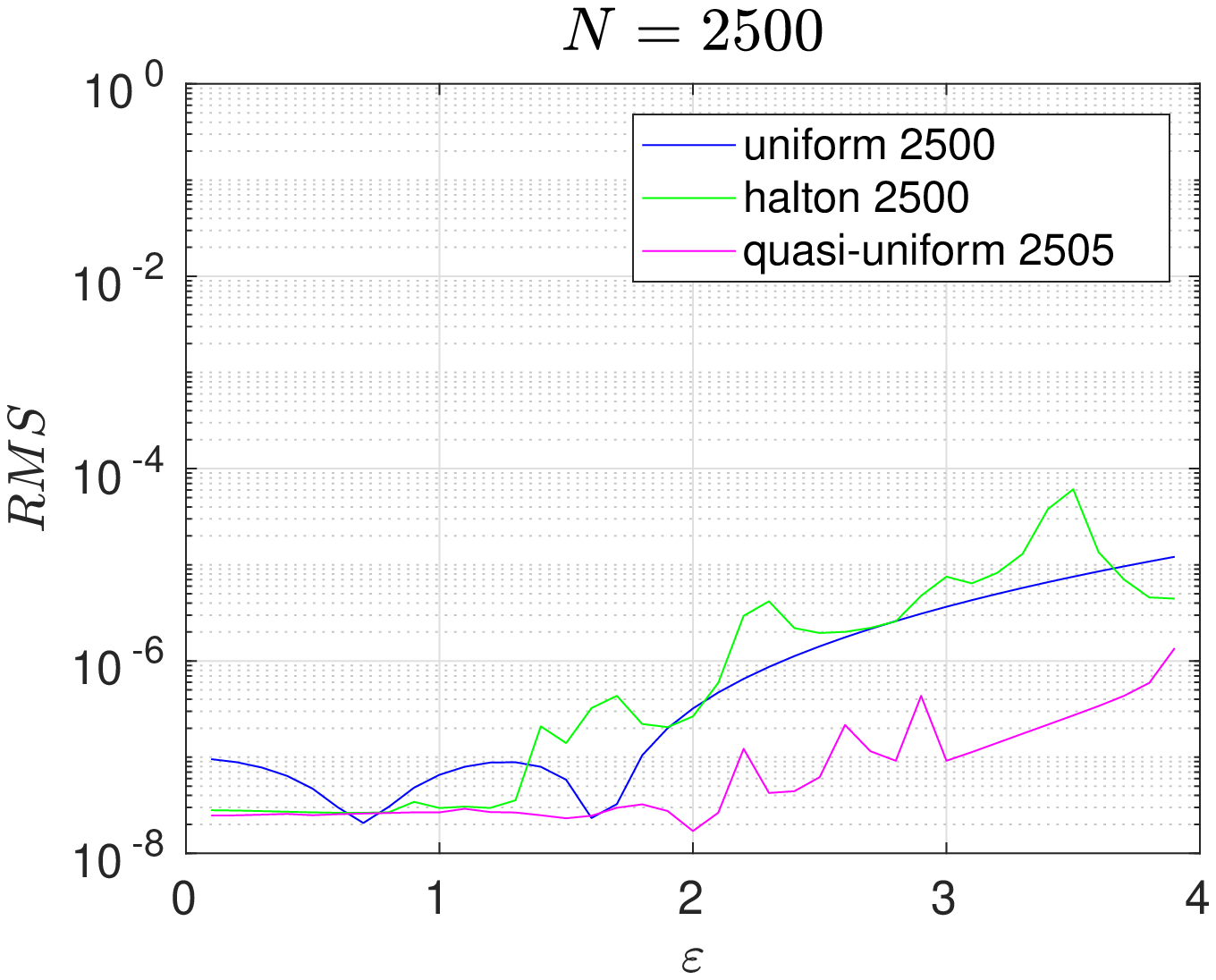}
   \end{center}
   \vspace{-.75cm}
   \caption{Comparison for uniform, halton and quasi-uniform node sets distributions for $N=400$ (left) and $N=2500$ (right) interior nodes.}
     \label{Fig:PEP5_comparison_distributions}
\end{figure}

Fig. \ref{Fig:PEP5_comparison_distributions} shows numerical experiments for the RMS error of the LIM RBF-QR for the uniform, Halton and quasi-uniform 
node sets distributions, with total number of interior points $N=400,900,1600,2500$ in each subfigure. In all cases, we see that the best accuracy of the 
numerical solution is obtained for the smaller values of $\varepsilon$ and also, we we increase the number of total collocation points. The method LIM RBF-QR 
with the uniform distribution shows the Runge phenomenon explained above as in Fig. \ref{Fig:PEP5_LRDint-Sint_rbf-qr_uniform_RMS}. For cases of halton points 
and quasi-uniform poitns, we see that the major change in the trend of $\varepsilon \rightarrow 0$ when the number of collocation points is increase from 
$N=900$ to $N=2500$. For $N=1600$, we see that the halton distribution follows the same pattern of the Runga phenomenon as in the uniform case, while the 
quasi-uniform distribution stabilizes for $\varepsilon <2$ and continues decreasing as $\varepsilon \rightarrow 0$. In the case $N=2500$, we see that the halton 
distribution stabilizes around $\varepsilon < 1.3$ while the quasi-uniform distribution around $\varepsilon <1.7$. This results are the expected for small
$\varepsilon$'s since the RBF-QR was originally implemented for small values.

All these results are due to the severe dependence of the conditioning of the interpolation matrix for the local interpolation for $u$ and $b$ 
of the shape parameter $\varepsilon$. As it was study in \cite{fornberg_zuev_07}, for a constant shape parameter, the interpolation matrix $\mathbf{A}$ 
of a 2-D non-periodic distribution of points for IQ, MQ or GA RBF is:
\begin{equation}
cond(\mathbf{A}) = O\left(\frac{1}{\varepsilon^{2[\sqrt{8n-7}-1]}}\right). 
\label{Eq:cond_inter_matrix}
\end{equation}
where $[.]$ denotes the integer part and $n$ the points of the local interpolation. So for the Halton and quasi-uniform distributions we have 
$cond(A)=O(\varepsilon^{-12})$. Solving the local linear systems of equations for the ill-conditioned interpolation matrix (\ref{u_def_bound}) and 
(\ref{Eq:int_term_no_hom}) with a direct method (RBF-Direct) give worst results that making the change of basis proposed in the RBF-Qr method.

\subsection{Poisson's equation with Dirichlet boundary conditions}
\label{PEP7}

Let's consider the following Poisson's equation defined in the domain $\Omega=[1,2]^2$:
\begin{eqnarray}
 \Delta u\left(x,y\right) = &-& \frac{751}{144} \pi^2 \sin\left(\frac{\pi}{6} x\right)\sin\left(\frac{7}{4} \pi x\right) \sin\left(\frac{3}{4} \pi y\right) \sin\left(\frac{5}{4} \pi y\right) \nonumber\\
 &+& \frac{7}{12} \pi^2 \cos\left(\frac{\pi}{6} x\right) \cos\left(\frac{7}{4} \pi x\right) \sin\left(\frac{3}{4} \pi y\right) \sin\left(\frac{5}{4} \pi y\right) \nonumber \\
 &+& \frac{15}{8} \pi^2 \sin\left(\frac{\pi}{6} x\right) \sin\left(\frac{7}{4} \pi x\right) \cos\left(\frac{3}{4} \pi y\right) \cos\left(\frac{5}{4} \pi y\right)
 \label{Poisson_equation}
\end{eqnarray}
with Dirichlet boundary conditions:
\begin{equation}
 (BCs)\left\{
 \begin{array}{lcl}
   u\left(x,1\right) &= & -\frac{1}{2}\sin\left(\frac{\pi}{6}x\right)\sin\left(\frac{7\pi}{4}x\right),\\
   u\left(x,2\right) &= & - \sin\left(\frac{\pi}{6}x\right)\sin\left(\frac{7}{4}y\right),\\
   u\left(1,y\right) &= & -\frac{1}{2\sqrt{2}}\sin\left(\frac{3\pi}{4}y\right)\sin\left(\frac{5\pi}{4}y\right),\\
   u\left(2,y\right) &= &  \frac{3}{\sqrt{2}}\sin\left(\frac{3\pi}{4}y\right)\sin\left(\frac{5\pi}{4}y\right).\\
 \end{array}
\right.
\label{boundaries_cond_possion}
\end{equation}
The analytical solution to this problem is: $ u\left(x,y\right) = \sin\left(\frac{\pi}{6} x\right) \sin\left(\frac{7}{4}\pi x\right) \sin\left(\frac{3}{4}\pi y\right) \sin\left(\frac{5}{4}\pi y\right)$.

As in the numerical experiment before, the node set distributions implemented were uniform, Halton and quasi-uniform for 
$N=400,900,1600,2500$ interior points, $N_b=84,124,164,204$ boundary points and $n=25$ points for the local stencils $\Theta_i$.

The objetive in solving this PDE was to compare numerical results of three different local integral methods reported in the literature, 
the results obtained with the LRDRM in \cite{caruso_portapila_power_15}, those found by Ooi and Popov \cite{ooi_popov_12} using the 
Radial Basis Integral Equation Method (RIBEM) and finally the LIM RBF-QR presented in this paper.

\begin{figure}[!ht]
 \centering
   \begin{center}
     \includegraphics[scale=.65]{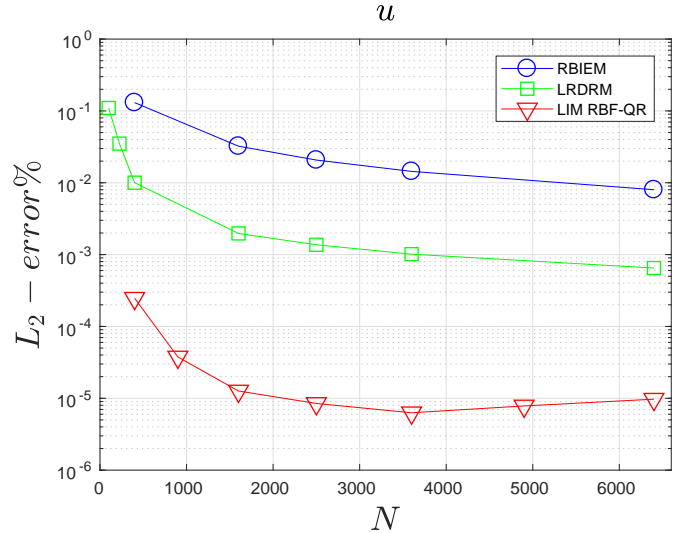}
   \end{center}
  \vspace{-.5cm}
 \caption{Comparison of the $L_2-error\%$ between RBIEM, LRDRM and LIM RBF-QR for uniform node set. 
 $N$ number of interior nodes of the domain.}
\label{Fig:PEP7_Error_L2p_vs_N_uniform} 
\end{figure}

\begin{figure}[!ht]
 \centering
   \begin{center}
     \includegraphics[scale=.31]{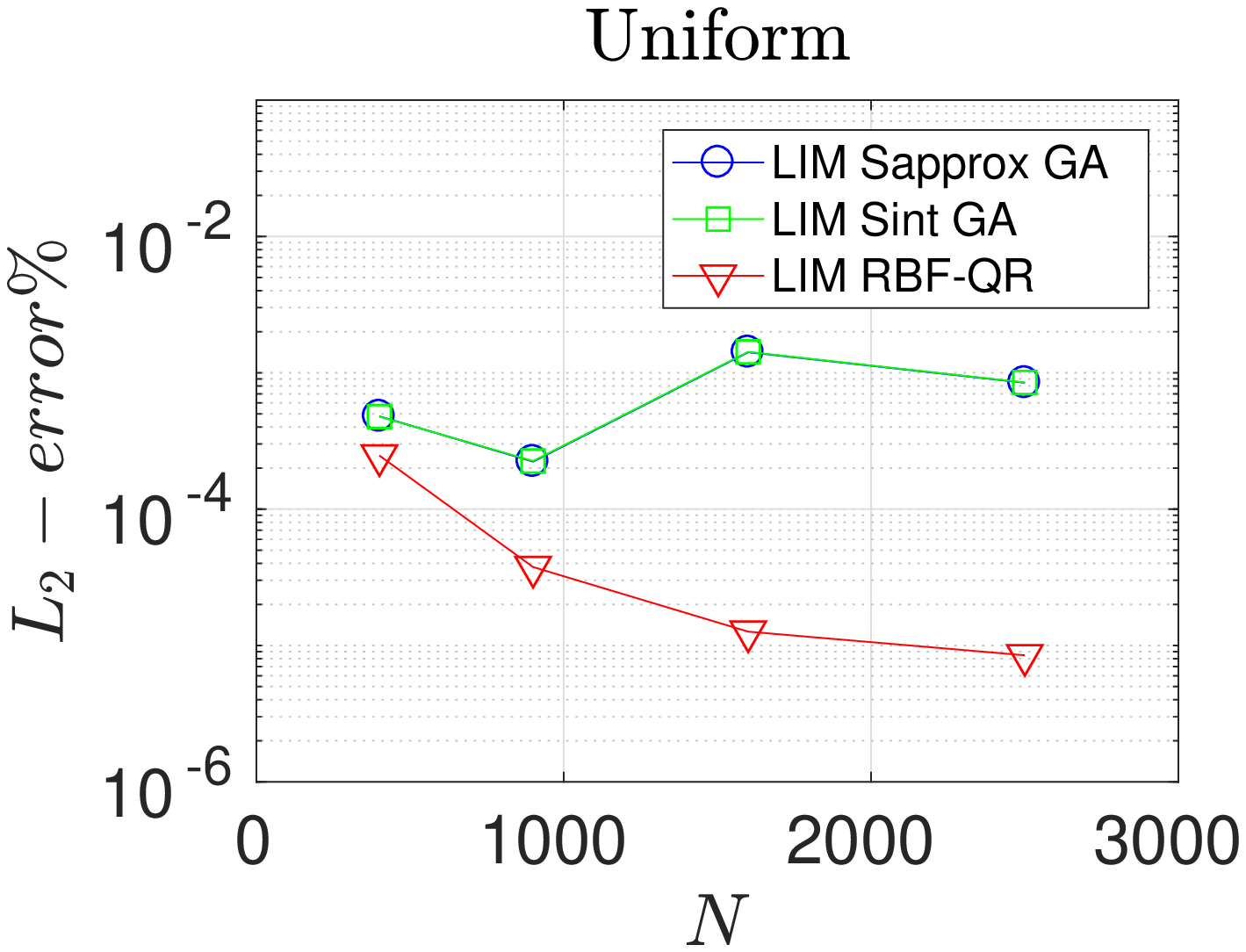}
     \includegraphics[scale=.31]{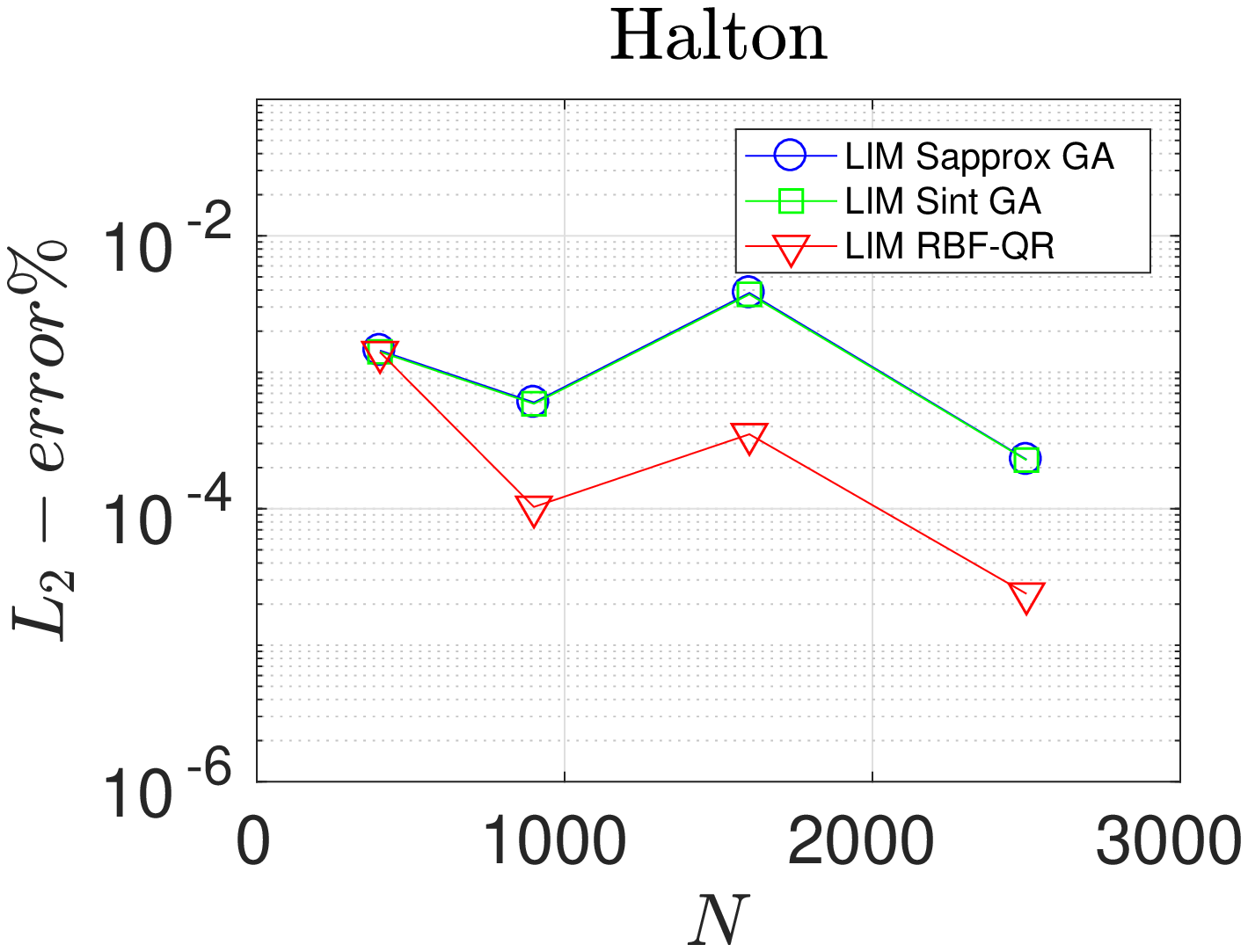}
     \includegraphics[scale=.31]{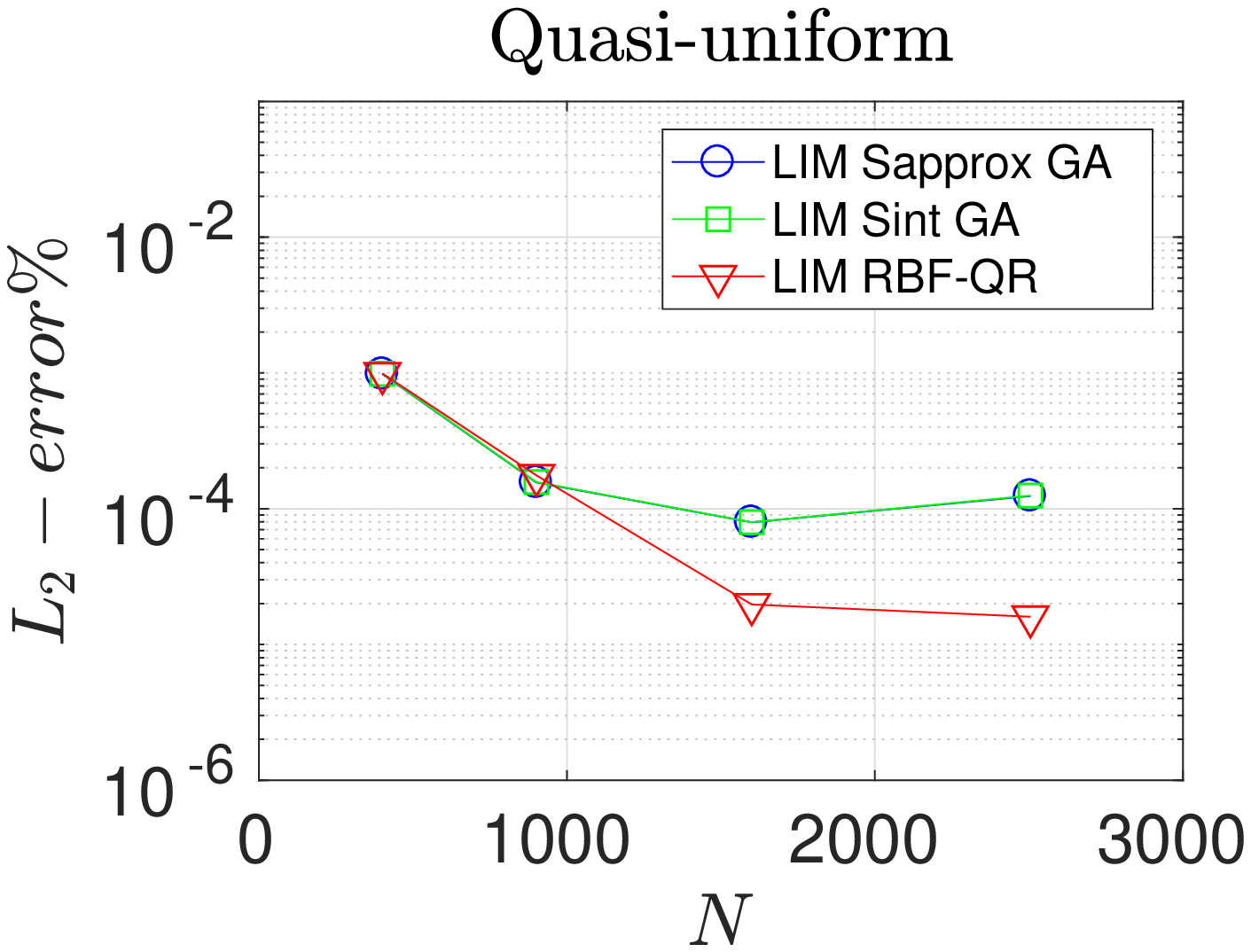}
   \end{center}
  \vspace{-.75cm}
 \caption{$L_2-error\%$ versus $N$ for uniform, Halton and quasi-uniform distributions}
 \label{Fig:PEP7_Error_L2p_vs_N_all_distribution} 
\end{figure}

For the uniform case, the convergence analysis through the L2\% error norm is presented in Fig. \ref{Fig:PEP7_Error_L2p_vs_N_uniform} and 
\ref{Fig:PEP7_Error_L2p_vs_N_all_distribution} for $u$. 
Comparison between the numerical results obtained with the LRDRM and the RBIEM, \cite{ooi_popov_12}, are also 
presented in Fig. \ref{Fig:PEP7_Error_L2p_vs_N_uniform}, where as before it can be seen that the LIM RBF-QR results 
are more accurate than those obtained with the LRDRM and RBIEM by one or two order of magnitude or more. 
For $N=6400$ the RBIEM achieved the $L_2-error\%$ of 8.000E-03 and the LRDRM 6.5021E-04. The best results for this 
PDE were obtained with LIM RBF-QR were the errors varies from 2.4698E-04, 3.75833E-05, 1.2626E-05, 8.4696E-06, 6.2864E-06 
for $N=400,900,1600,2500,3600$ respectively.

In Fig. \ref{Fig:PEP7_Error_L2p_vs_N_all_distribution} we compare the $L_2-error\%$ versus the number of interior points for the different types 
of distribution node sets. As expected, all the schemes show that the stability with RBF-QR for the local interpolation matrix of the unknown field $u$ and 
the non-homogeneous term $b$ gave the best numerical results, that we call LIM RBF-QR. As in this case $b$ depends just on the variable $\mathbf{x}$, 
it was integrated in two different ways. The first one was to considered the integral over the domain $\Omega_i$ like in Eq. (\ref{int_formulation_split}) 
being $b_2=0$ (that was called LIM Sint GA). And the second one, was to considered $b$ as a linear combination of Gaussians RBFs Eq. (\ref{Eq:int_term_no_hom})) 
and then interpolated (called LIM Sapprox GA). This numerical modification gave no numerical difference in orders of magnitude of the error.

This two alternatives of integrating $b$ had a better behavior for the quasi-uniform distribution arriving to values of the $L_2-error\%$ as low as 
7.9410E-05 and 7.9549E-05 respectively for $N=1600$. In the uniform case, the $L_2-error\%$ achieve 2.2318E-04 and 2.2530E-04 for $N=900$.
For the LIM RBF-QR the best behavior is obtained with uniform and quasi-uniform points achieving 8.4696E-06 in the uniform case for $N=2500$ and 
1.5995E-05 in the quasi-uniform case for $N=2505$.

The Halton case does not show uniform convergence for all the methods.

\subsection{Poisson's equation over the unit disk}
\label{PEP8}

The final Poisson's equation with Dirichlet boundary conditions is defined over the circular domain 
$\Omega=\{(x,y)/ x^2+y^2 \leq 1\}$:
\begin{equation}
 \left\{
 \begin{array}{lcl}
  \Delta u &=& -200 \sin[10(x+y)] \hspace{.5cm}  (x,y)\in \Omega,\\
  u &=& \sin[10(x+y)] \hspace{.5cm} (x,y)\in \partial\Omega \quad (BCs).
 \end{array}
\right.
\label{PDE_Bayona_case1}
\end{equation}

\begin{figure}[!ht]
  \begin{center}
  \includegraphics[scale=.45]{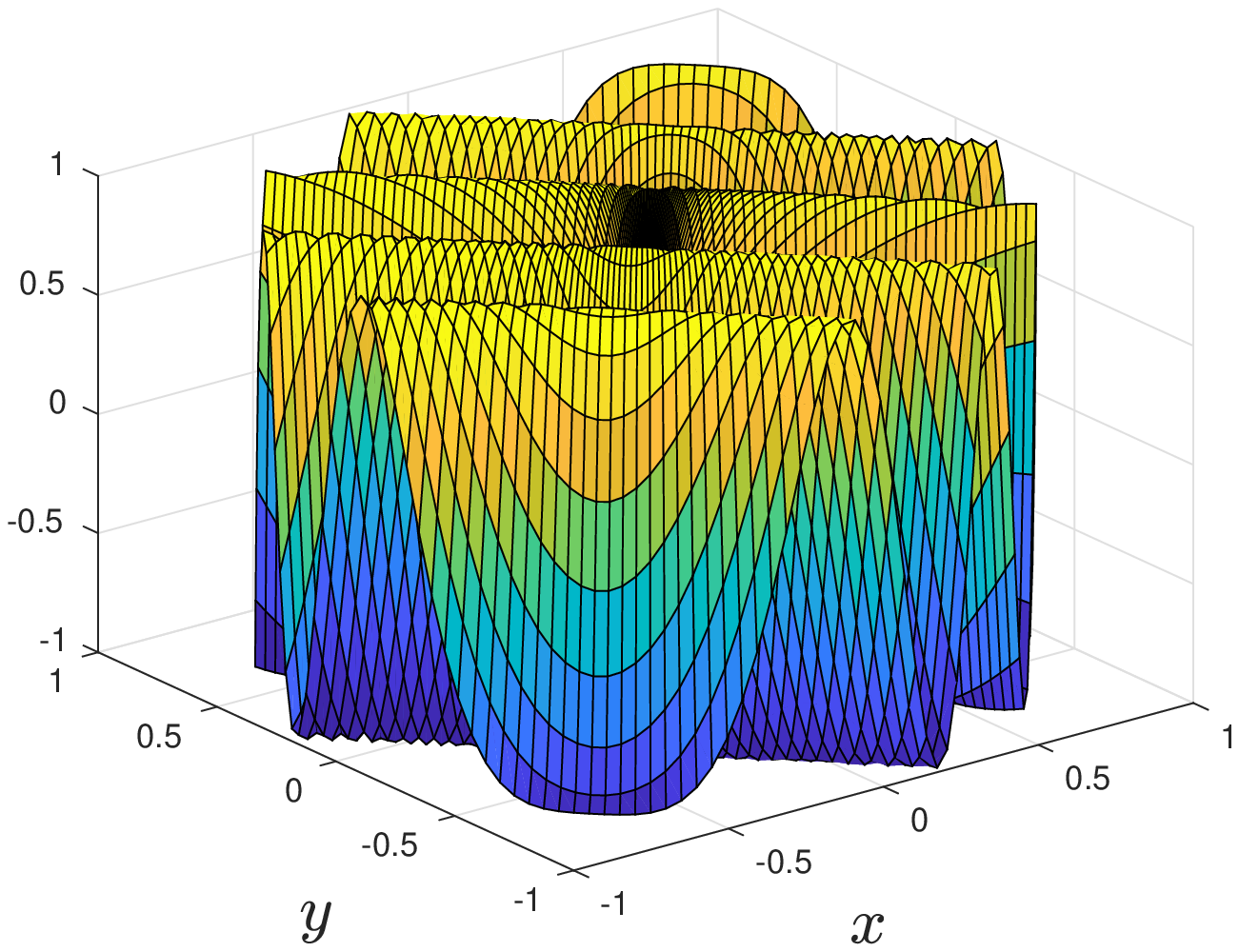}
  \hfill
  \includegraphics[scale=.45]{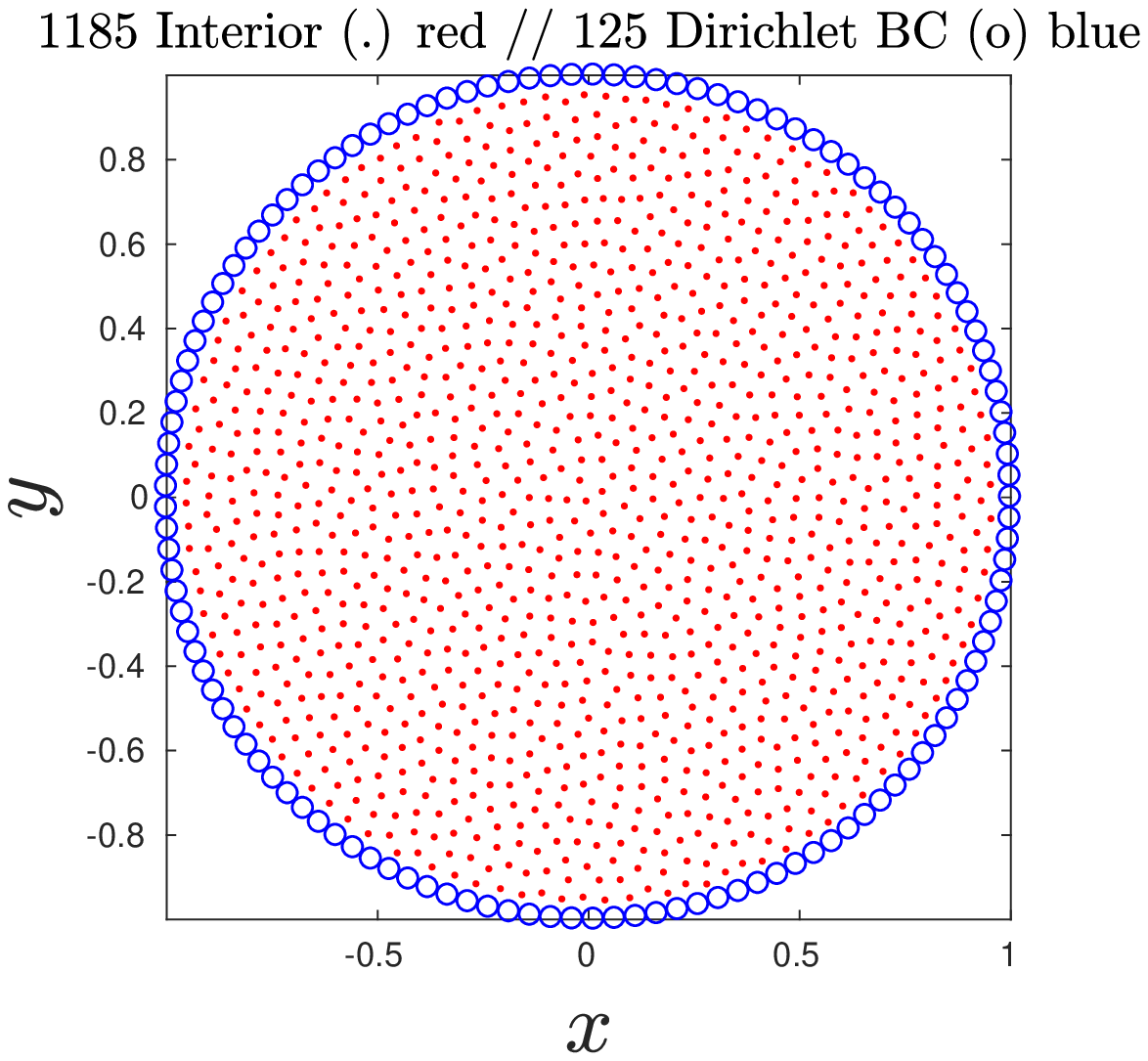}
  \end{center}
  \vspace{-.5cm}
  \caption{Exact solution (left) and quasi-uniform point distribution (right)}
  \label{Fig:bayona_et_al_exact_domain}
 \end{figure}

The exact solution to this problem is plotted in Fig. \ref{Fig:bayona_et_al_exact_domain} and is given by: 
$u\left(x,y \right) = \sin[10(x+y)].$

The results are compared with those found by Bayona et al \cite{bayona_flyer_fornberg_barnett_16} using RBF-Generatad 
Finite Difference method with polynomial augmentation. This method is a combination of polyharmonic splines (PHS) with multivalue polynomials 
for solving PDEs that has emerge in the last decade as a powerful and flexible numverical approach. Int that paper, several strategies 
were used to avoid the accuracy and stability problems of using one-side stencils near boundaries for elliptic PDEs. Besides, our objetive in 
this experiment is to compare integral and difference numerical method that use near flat RBF without any treatment at the boundaries and 
also using the RBF-QR method to stabilize the local interpolarion errors.

The present scattered nodes distribution to experiment was formed by $N=1185,4880$ interior nodes and $N_b=125,251$ boundary nodes respectely.
Fig. \ref{Fig:bayona_et_al_exact_domain} (right) shows the nodes distribution for $h=0.025$ ($N=1185$ interior nodes), where the structure 
for $h=0.01$ ($N=4880$ interior nodes) is similar, but denser. The left figure shows the oscillatory behavior of the analytical solution over the 
unit disk.

\begin{figure}[!ht]
 \centering
   \begin{center}
   \includegraphics[scale=.5]{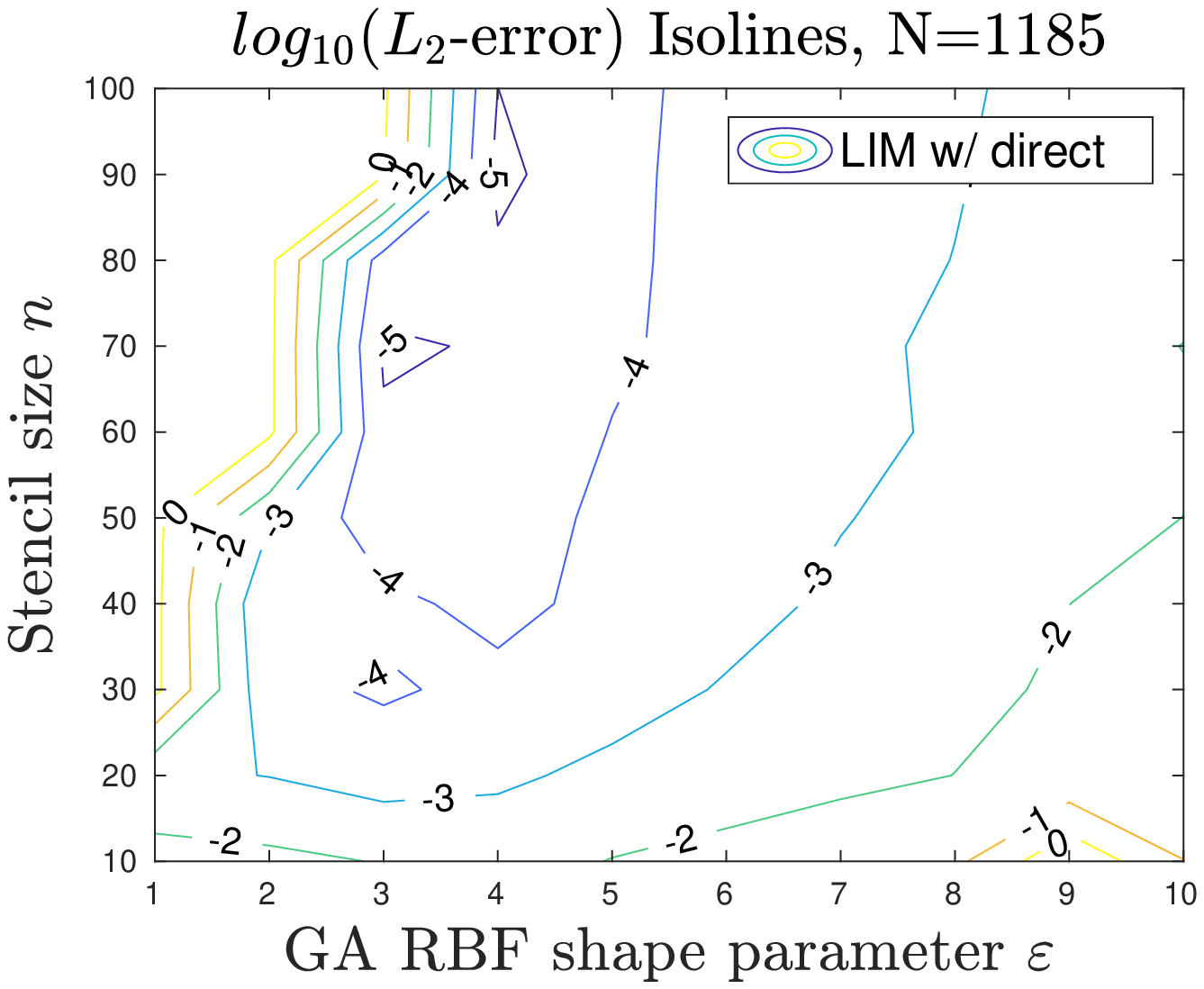}
   \hfill
   \includegraphics[scale=.5]{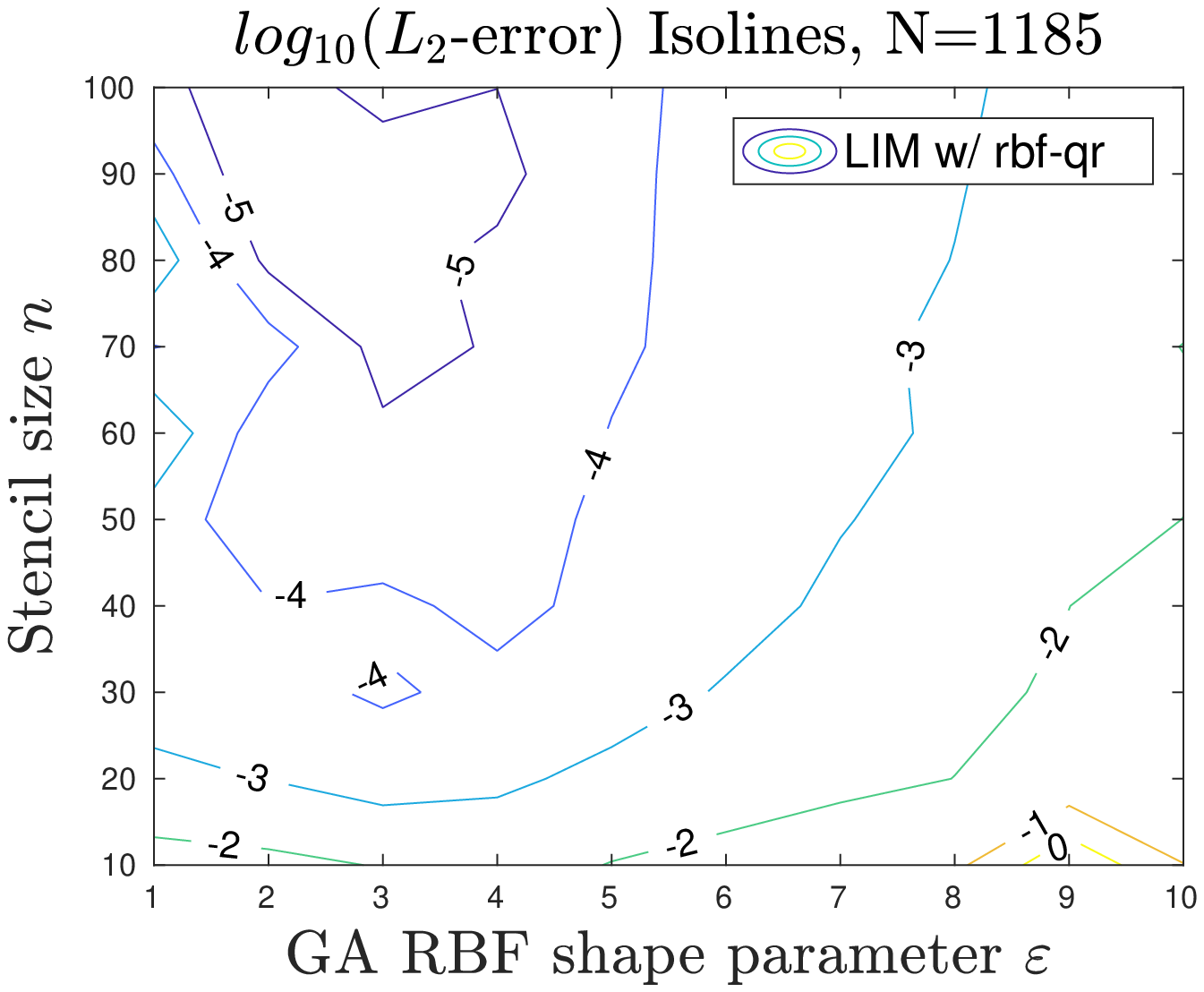}\\
   \vspace{.5cm}
   \includegraphics[scale=.5]{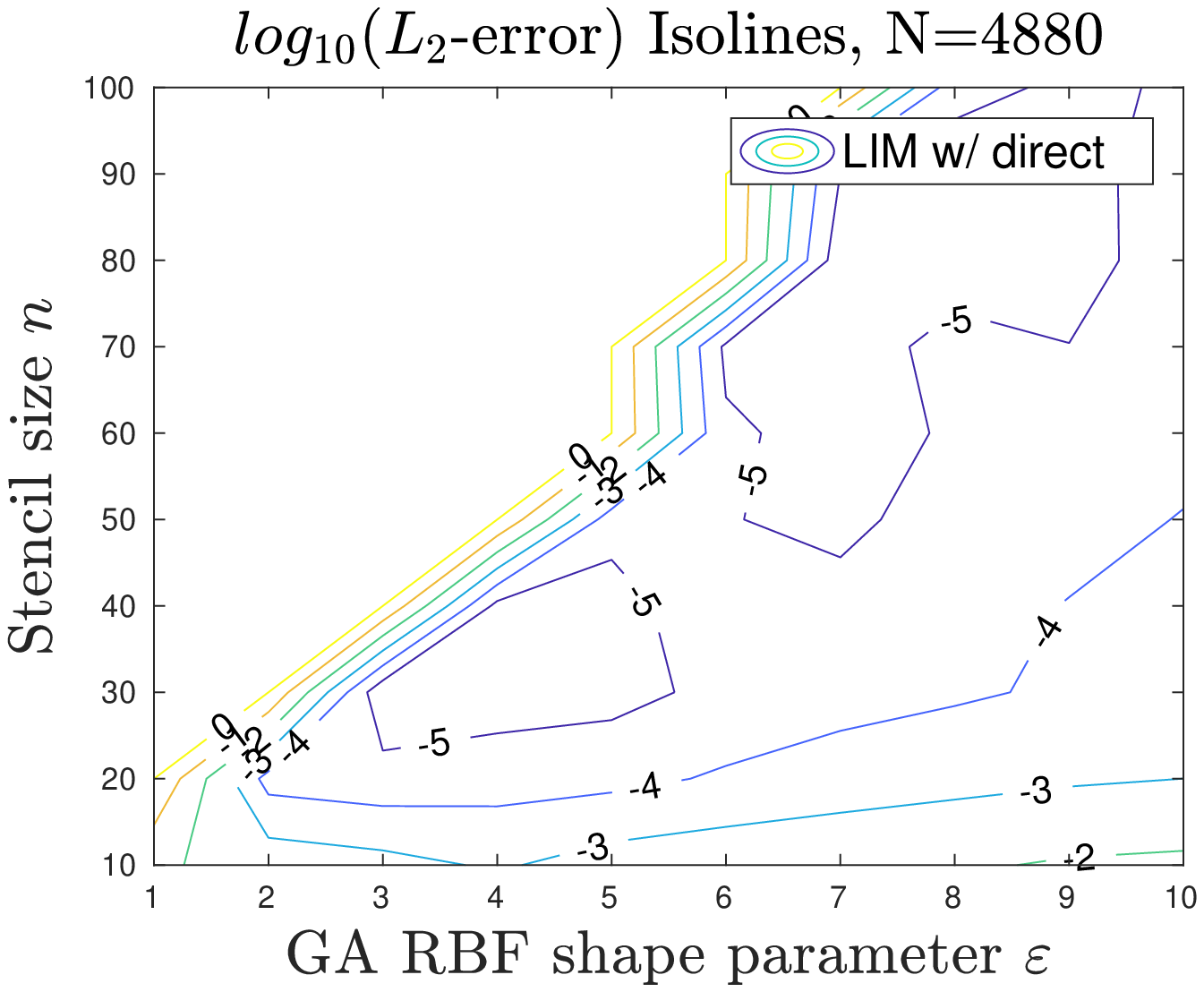}
     \hfill
   \includegraphics[scale=.5]{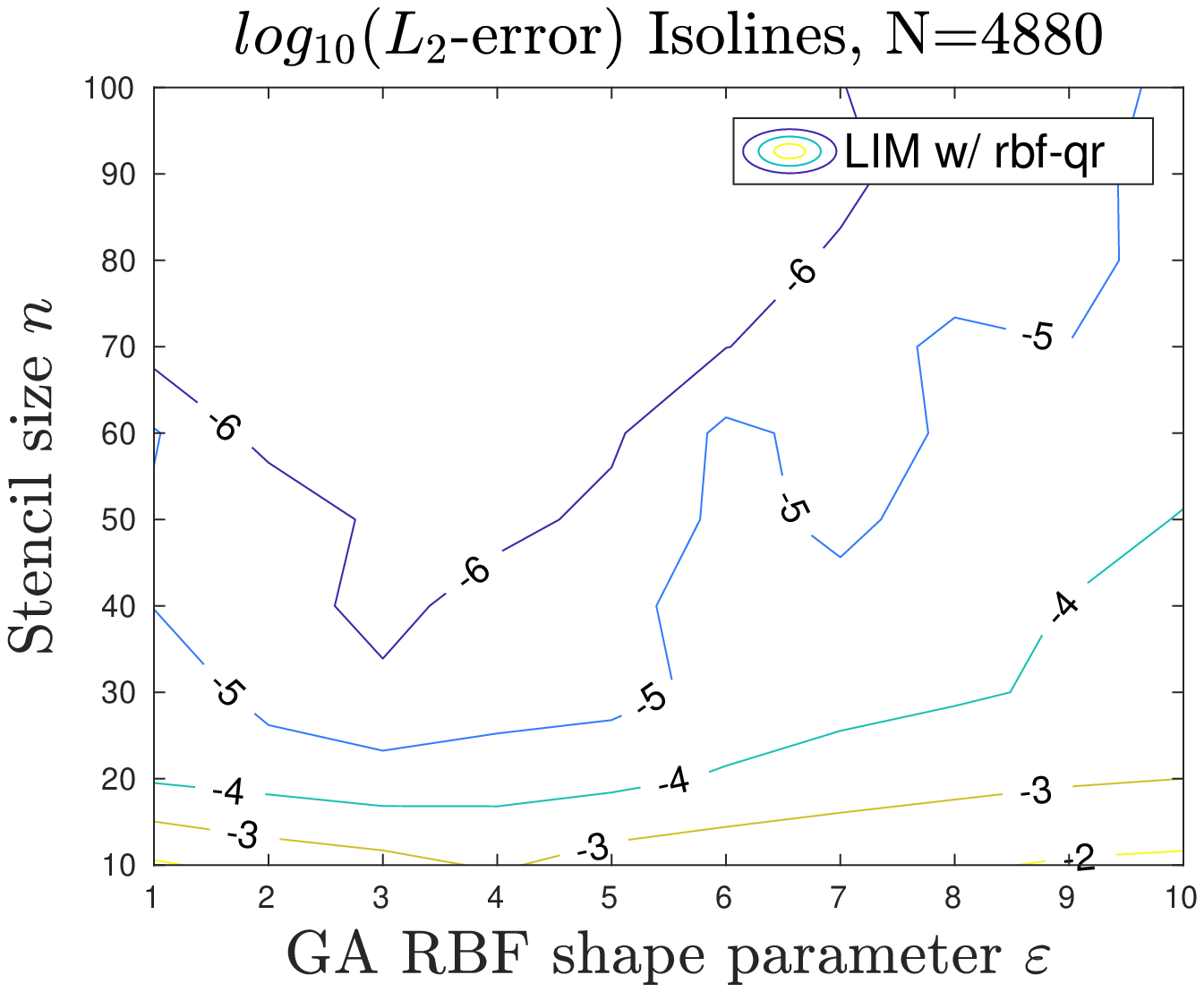}\\
   \end{center}
   \vspace{-.5cm}
   \caption{Accuracy ($log_{10}(L_2-error)$ isolines) when using RBF GA as a function of the stencil size and shape parameter for 
   quasi-uniform distribution for $N=1185$ and $N=4880$ with and without RBF-QR.}
   \label{Fig:bayona_isolines}
\end{figure}

Fig. \ref{Fig:bayona_isolines} shows the isolines order $L_2-error$ as a function of the stencil size and the shape parameter for the different 
node sets (row subplots) and local integral numerical methods with/without RBF-QR (column subplots). Without any special treatment of the boundary, 
the interior nodes are used as collocation nodes obtaning a one side stencils near the boundary. The number of points for local stencils varies 
from $n=10$ to $n=100$. 

The first observation is that the error decreases as the node distribution is refined from $N=1185$ to $N=4880$ and also when the shape parameter 
$\varepsilon$ tends to zero in the four subplots. The introduction of the RBF-QR method in the local interpolation of the integral methods 
produces the same effects than in the Poisson equations before. Decreasing the shape parameter of the Gaussian RBF for the local interpolation, 
improved the accuracy of the LIM until the breakdown error ocurrs. 

In the first column of subplots we show a wide region of inestability for $N=1185$ that began for $\varepsilon <4$ and $n>30$ to the left and up. 
The best order error is $10^{-5}$ obtained in two peaks bewteen $3<\varepsilon <3.6$ and $65<n<70$ and $4<\varepsilon <4.5$ and $85<n<100$.
Also for $N=4880$ this inestability's region began from $\varepsilon <8$ and $n>20$ to left and up. This is because the near flat Gaussians RBFs 
produce ill-conditioned matrix for the local interpolation matrices that increases their size so the direct local solver increases the local error. 
Is is expectable that  increasing the number of points in the local stencils, the accuracy of the LRDIM decreases. The best order error is $10^{-5}$ 
from $3<\varepsilon <5.6$ and $20<n<40$ and there is another region of the same order for $6<\varepsilon <9.8$ and $45<n<100$.

The second column shows the effecs of the RBF-QR. The inestability region expands from $1.5<\varepsilon <4.2$ and $63<n<100$ for order $10^{-5}$ at 
the case $N=1185$, and from $1<\varepsilon <7$ and $35<n<100$ for order $10^{-6}$ in the case $N=4880$.

\subsection{One-dimensional Convection-Diffusion equation}
\label{subsect conv-diff}
To test the performance of the proposed LIM RBF-QR on differential problems that comes from applications, 
we consider a steady state convection-diffusion equation with a variable velocity field that 
has been used before as test example of different implementations in the literature \cite{popov_bui_10,portapila_power_08}
\begin{equation}
 D\frac{\partial^{2} u\left(\mathbf{x}\right)}{\partial x_i^2} - V_{x_1}\frac{\partial u\left(\mathbf{x}\right)}{\partial x_1} - 
 k u\left(\mathbf{x}\right) = 0
\label{Eq:eq_conv_diff_1D}
\end{equation}
with convective velocity field 
\begin{equation}
V_{x_1}=\ln \frac{U_{1}}{U_{0}}+k\left( x_1-\frac{1}{2}\right) ;\ V_{x_2}=0
\end{equation}
corresponding to the flow of a hypothetical compressible fluid with a density variation inversely
proportional to the velocity field.  

The analytical solution of the above boundary value problem for a diffusion coefficient $D=1$ is given by
\begin{equation}
 u\left(\mathbf{x}\right)=U_0\exp\left\{\frac{k}{2}x^2_1+\left(\ln\frac{U_1}{U_0}-\frac{k}{2}\right)x_1\right\}
 \label{solution}
\end{equation}
showing the formation of shock structures at each side of the problem domain. To analyse the performance of the 
numerical scheme different values of the decay parameter $k$ in the convective velocity are considered, where 
larger values of $k$ correspond to stronger shock structures at the problem boundaries.

For the numerical solution, this 1-D problem is considered as a 2-D one in a rectangular domain 
$\Omega = [0, 1] \times [-0.1, 0.1]$ subject to the following boundary conditions
\begin{equation}
 \left\{
 \begin{array}{lc}
   u\left(0,x_2\right)=U_0, & u\left(1,x_2\right)=U_1,\\
   \frac{\partial u}{\partial x_2}\left(x_1,-0.1\right)= 0, &\frac{\partial u}{\partial x_2}\left(x_1,0.1\right) = 0,\\
 \end{array}
\right.
\label{boundaries_cond_conv_diff_1D}
\end{equation}
where the domain $\Omega$ is subdivided into subdomains $\Omega_i$ that are used to construct the interpolation 
stencils $\Theta_i$. The distribution node sets considered were uniform, Halton and quasi-uniform. 
See Table \ref{Table:CDRE_discretizations}.

\begin{table}[!ht]
\centering
\begin{small}
\begin{tabular}{|c c c|c c c|c c c|}
\hline
  \multicolumn{3}{|c|}{Uniform} & \multicolumn{3}{|c|}{Halton} & \multicolumn{3}{|c|}{Quasi-uniform} \\
  $N$ & $N_{Dir}$ & $N_{Neu}$ & $N$ & $N_{Dir}$ & $N_{Neu}$ & $N$ & $N_{Dir}$ & $N_{Neu}$\\
\hline
   500 & 100 & 24 &   500 & 100 & 20 &   500 & 100 & 22 \\
  1125 & 150 & 34 &  1125 & 150 & 30 &  1127 & 150 & 32 \\ 
  2000 & 200 & 44 &  2000 & 200 & 40 &  1981 & 198 & 42 \\ 
  3125 & 250 & 54 &  3125 & 250 & 50 &  3125 & 250 & 52 \\
  4500 & 300 & 64 &  4500 & 300 & 60 &  4501 & 300 & 62 \\ 
  6125 & 350 & 74 &  6125 & 350 & 70 &  6158 & 352 & 74 \\
  8000 & 400 & 84 &  8000 & 400 & 80 &  7987 & 404 & 84 \\
\hline
\end{tabular}
\caption{Number of discretizations}
\label{Table:CDRE_discretizations}
\end{small}
\end{table}

\begin{figure}[!ht]
\centering
  \begin{center}
    \includegraphics[scale=.5]{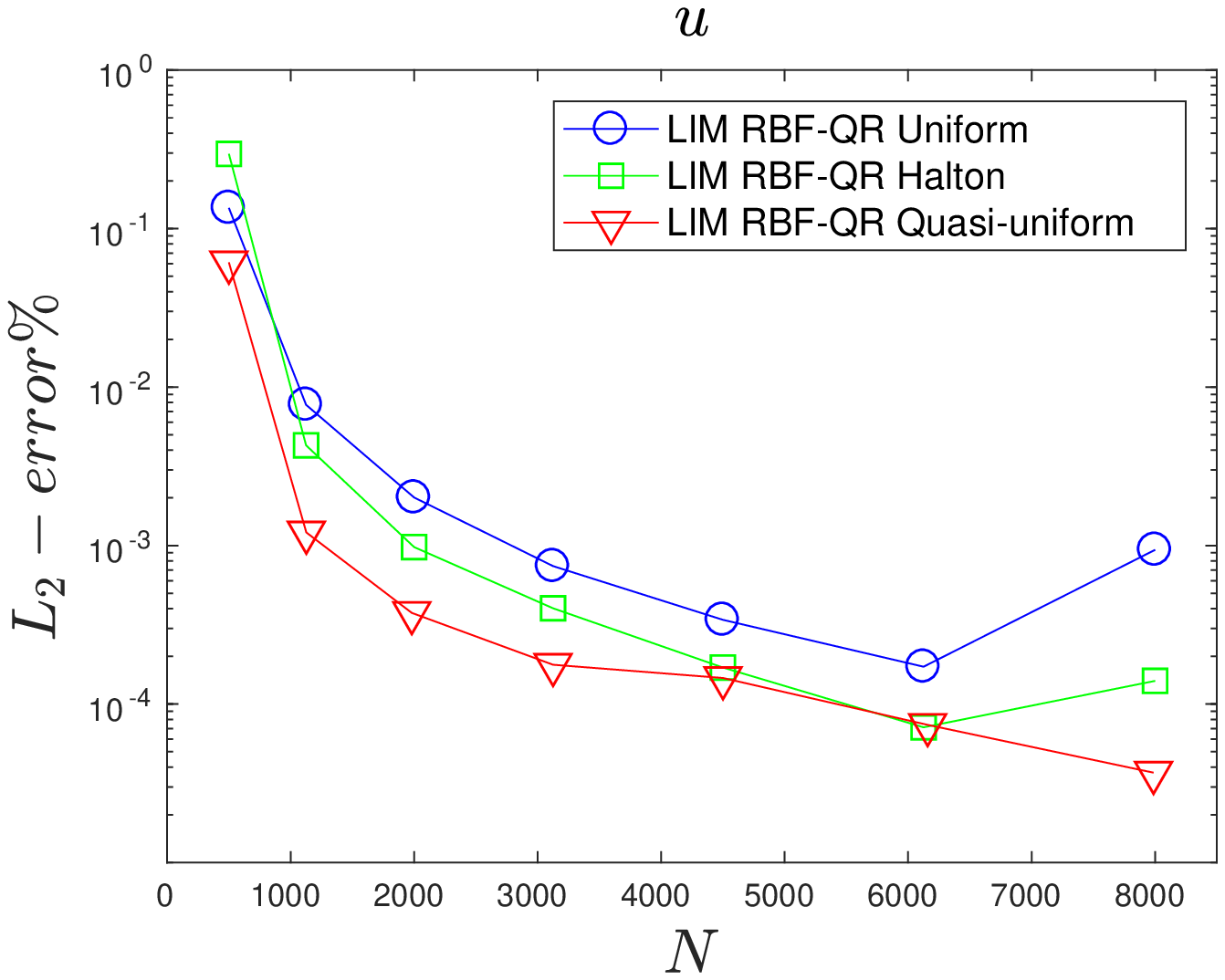}
    \hfill
    \includegraphics[scale=.5]{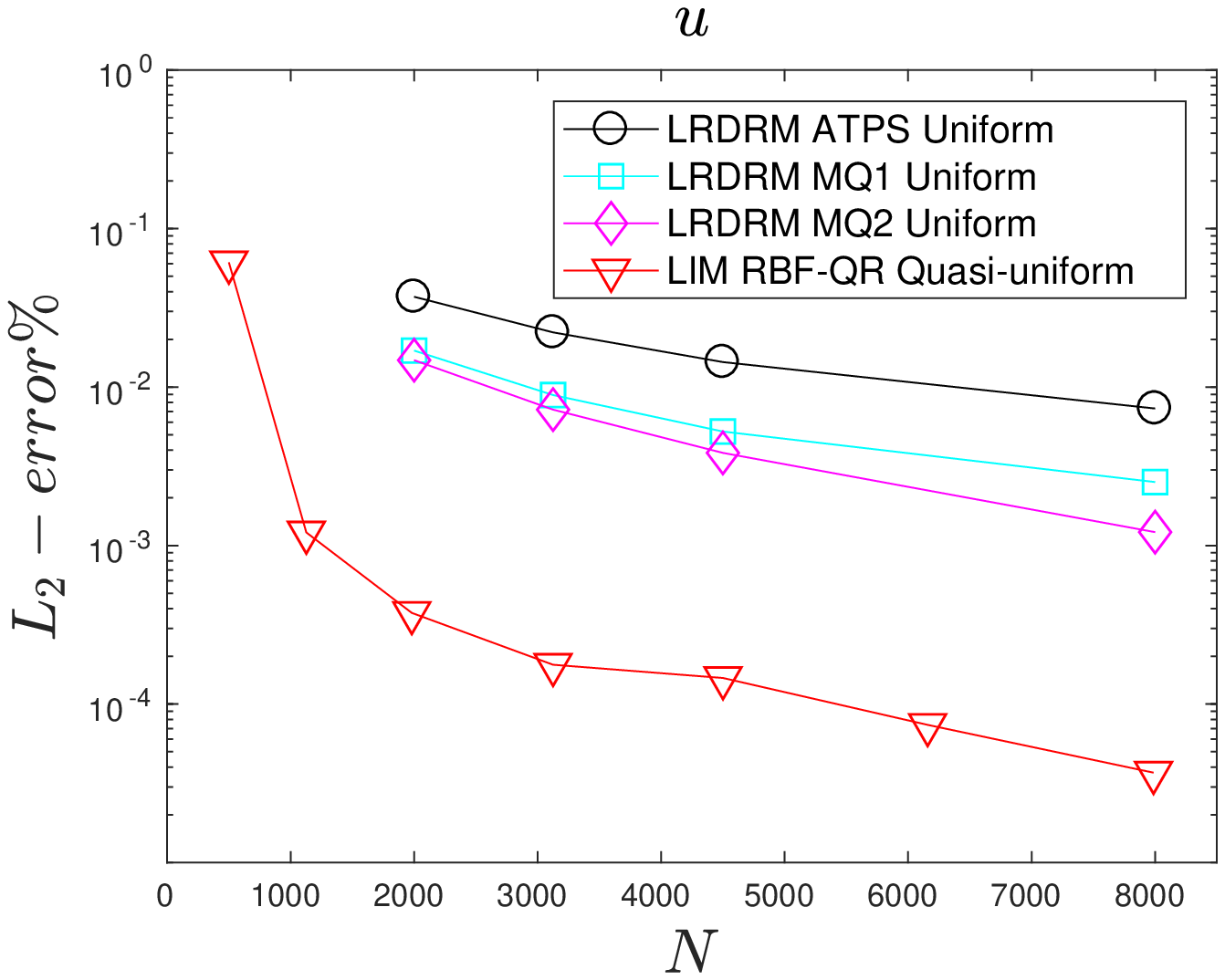}
  \end{center}
  \vspace{-.75cm}   
\caption{Comparison between different methods $k=40$}
  \label{Fig:CDRE_1D_comparison}
\end{figure} 

The difference between the LIM RBF-QR and LRDRM presented in \cite{caruso_portapila_power_15} are shown in Fig. \ref{Fig:CDRE_1D_comparison} (left). It can be seen that for $k=40$, the best $L_2-Error\%$ is obtained for the LIM RBF-QR over the quasi-uniform points, that is using locally RBF-QR method to achieve low errors with small shape parameters. In this case the integral equation is applied at only one collocation point per subdomain located at its centre with each stencil subdomain having $n=25$ points every stencil. With this 25-points stencils for the quasi-uniform distribution the values of the $L_2-Error\%$ is as low as $3.6878E-05$ for $u$. 
In Table \ref{Table:CDRE_1D_comparison} we show the comparison for the different distributions and different values of $k$. It is observed that the results for the LRDRM where obtained with $N=20480$ uniform interior nodes, but in the case $k=40$ and $k=100$ for the LIM RBF-QR we need less points to achieve one order of magnitude less. For the case $k=200$ we achieved the same order that the LRDRM with $N=20480$ interior points but with $N=8000,6125,7987$ for the distributions uniform, halton and quasi-uniform respectely. The LIM RBF-QR for $N=7987$ quasi-uniform nodes achieved $1.6580E-02$ lower that $3.4485E-02$ for $N=20480$ for uniform points.

Also, if we compare the different distributions we observed that the best result was obtained for the quasi-uniform points. See Fig \ref{Fig:CDRE_1D_comparison} (right). For $k=40$ with the LIM RBF-QR the uniform and Halton nodes set distribution the $L_2-Error\%$ versus the number of interior points decreases from $N=500$ to $N=6125$ achieving the best errors values of $1.7186E-04$ and $7.1334E-05$ for $N=6125$. While for the quasi-uniform, the error continue decreasing until $3.6878E-05$ for $N=7987$ points.

\begin{table}[!ht]
\centering
\begin{scriptsize}
\begin{tabular}{|c|c c|c c|c c|c c|}
\hline
       & \multicolumn{2}{|c|}{LRDRM}  & \multicolumn{6}{|c|}{LIM RBF-QR} \\
\hline       
       & \multicolumn{2}{|c|}{{\bf Uniform}} & \multicolumn{2}{|c|}{{\bf Uniform}} & \multicolumn{2}{|c|}{{\bf Halton}} & \multicolumn{2}{|c|}{{\bf Quasi-uniform}} \\
 $k$   & $N$ & $L_2-error\%$ & $N$ & $L_2-error\%$ & $N$ & $L_2-error\%$ & $N$ & $L_2-error\%$\\
\hline
 $40$  & 20480 & 6.83E-04 &  6125 & 1.7186E-04 &  6125 & 7.1334E-05 & 7987 & 3.6878E-05 \\ 
 $100$ & 20480 & 4.32E-03 &  6125 & 1.1654E-04 &  6125 & 2.8154E-03 & 7987 & 5.8630E-04 \\
 $200$ & 20480 & 3.44E-02 &  8000 & 6.2399E-02 &  6125 & 8.7114E-02 & 7987 & 1.6580E-02 \\
\hline
\end{tabular}
\caption{\textbf{CDRE 1D} - For different $k$ - $L_2-error\%$}
  \label{Table:CDRE_1D_comparison}
\end{scriptsize}
\end{table}

\subsection{Thermal boundary layer in a two dimensional channel}

As a final example, let us consider a 2D steady state flow entering a parallel channel with different walls temperatures. The governing equation is
\begin{equation}
 \Delta T\left(x_1,x_2\right)-Pe V\left(x_2\right) \frac{\partial T}{\partial x_1} = 0
\label{thermal_bound_layer_2D}
\end{equation}
where the parabolic velocity distribution is $V(x_2)=4 x_2(x_2-1)$ and $Pe$ is the P\'eclet number. The computational domain is 
taken to be $\Omega = [0,1]^2$ and the following Dirichlet and Neumann boundary conditions are imposed:
\begin{equation}
 \left\{
 \begin{array}{lcl}
   T\left(x_1,0\right) &=& 1 \hspace{1cm} 0\leq x_1 \leq 1\\
   T\left(x_1,1\right) &=& 0 \hspace{1cm} 0\leq x_1 \leq 1\\
   T\left(0,x_2\right) &=& 0 \hspace{1cm} 0\leq x_2 \leq 1\\
   \frac{\partial T}{\partial x_1} \left(1,x_2\right) &=& 0 \hspace{1cm}  0\leq x_2 \leq 1\\
 \end{array}
 \right.
\label{boundaries_cond_conv_diff_2D}
\end{equation}
There is no analytical solution for this PDE.

\begin{figure}[!ht]
\centering
  \begin{center}
    \includegraphics[scale=.325]{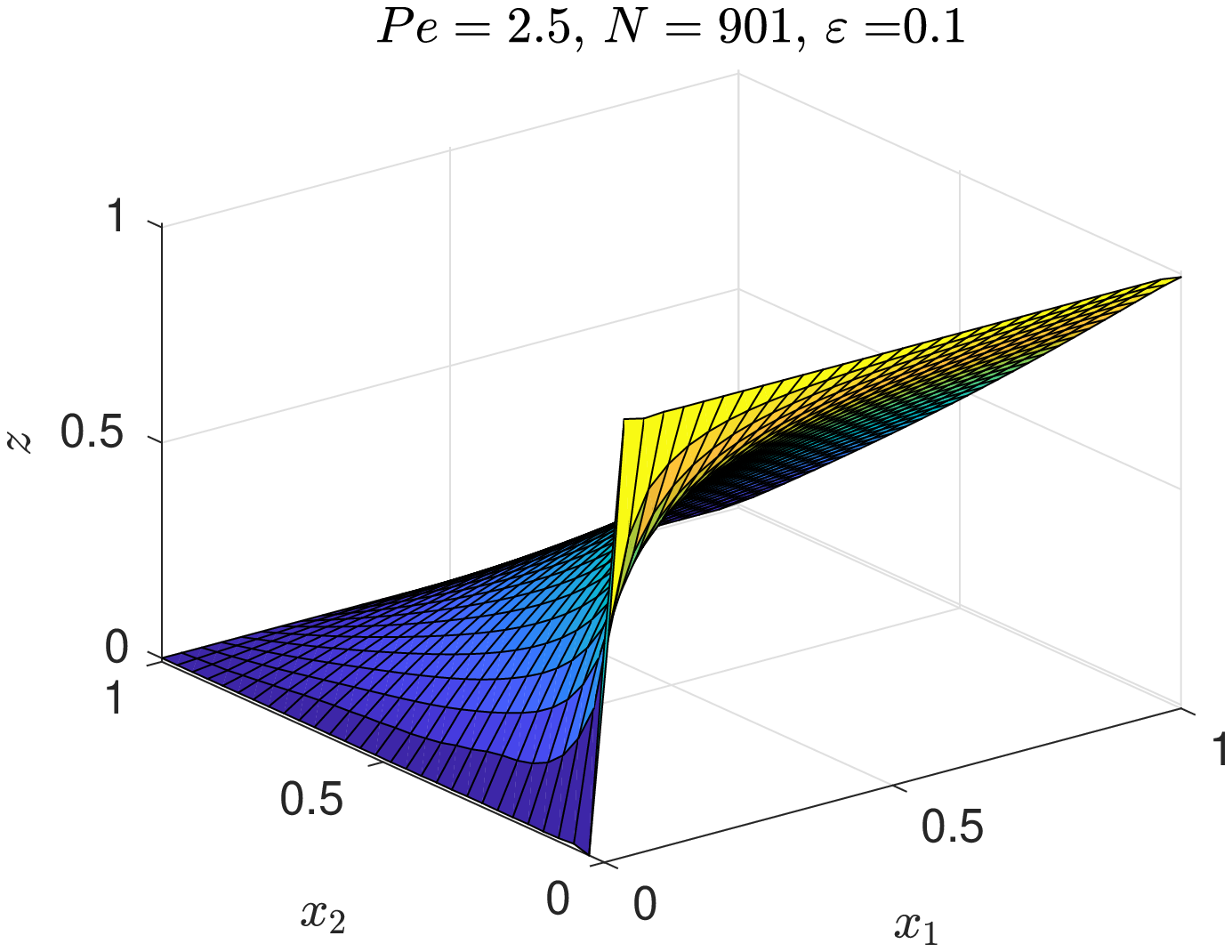}
    \includegraphics[scale=.325]{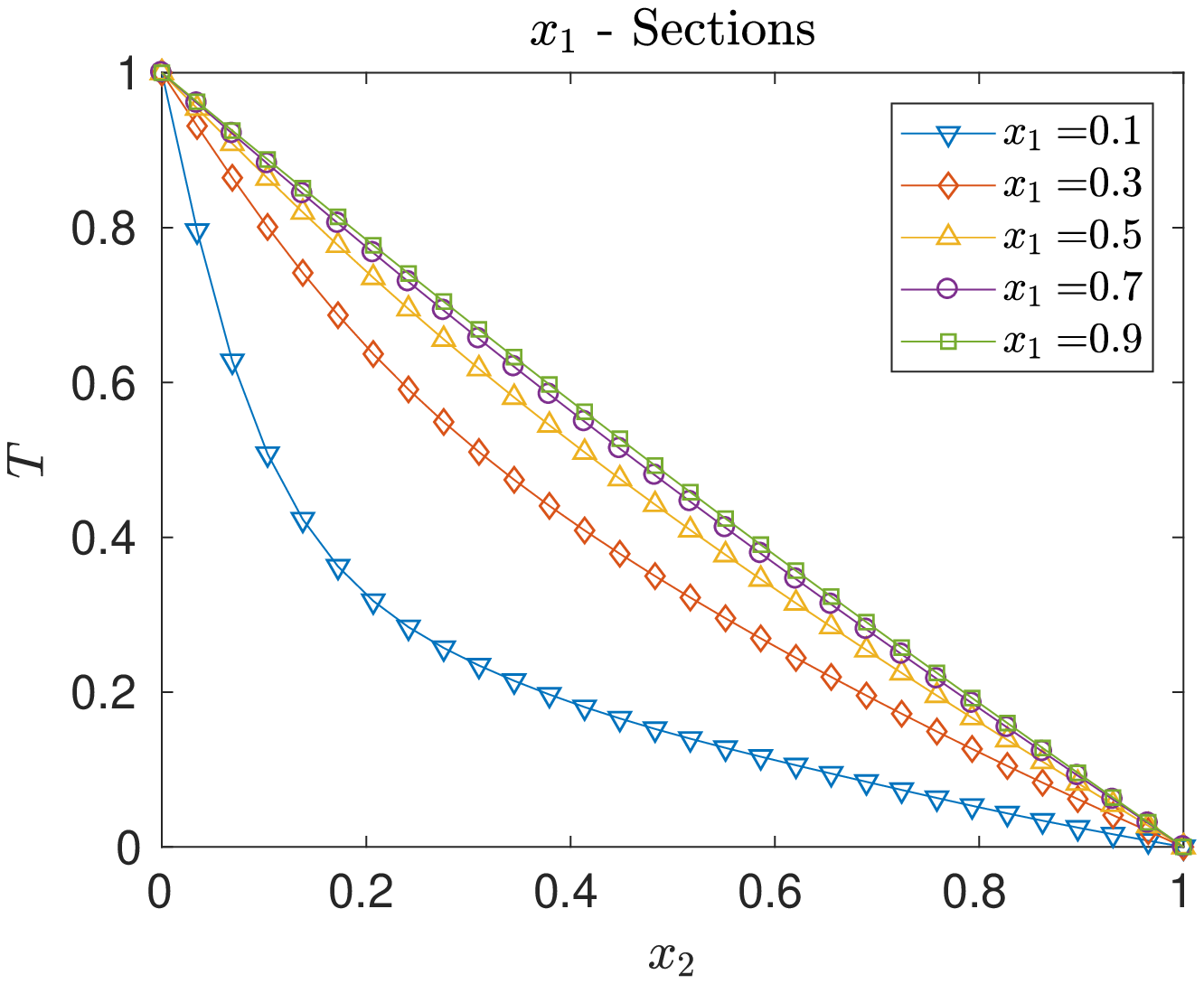}
    \includegraphics[scale=.325]{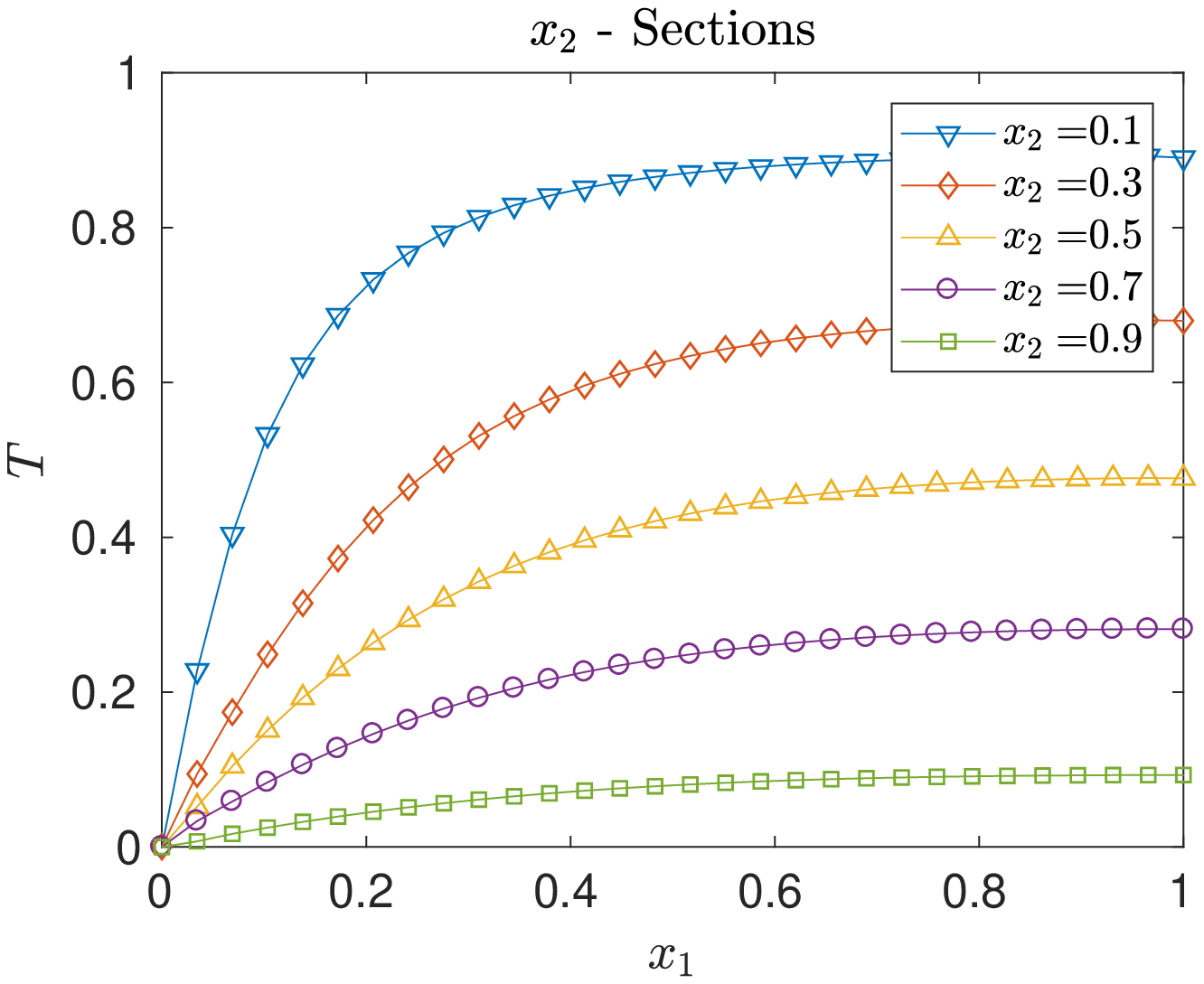}\\
   \end{center}
   \vspace{-.5cm}
  \caption{LIM RBF-QR - Quasi-uniform point distribution - $Pe=2.5$ - $N=901$ - $\varepsilon = 0.1$}
  \label{Fig:TBL_P1_ErrorDomain}
\end{figure}

\begin{figure}[!ht]
\centering
  \begin{center}
    \includegraphics[scale=.325]{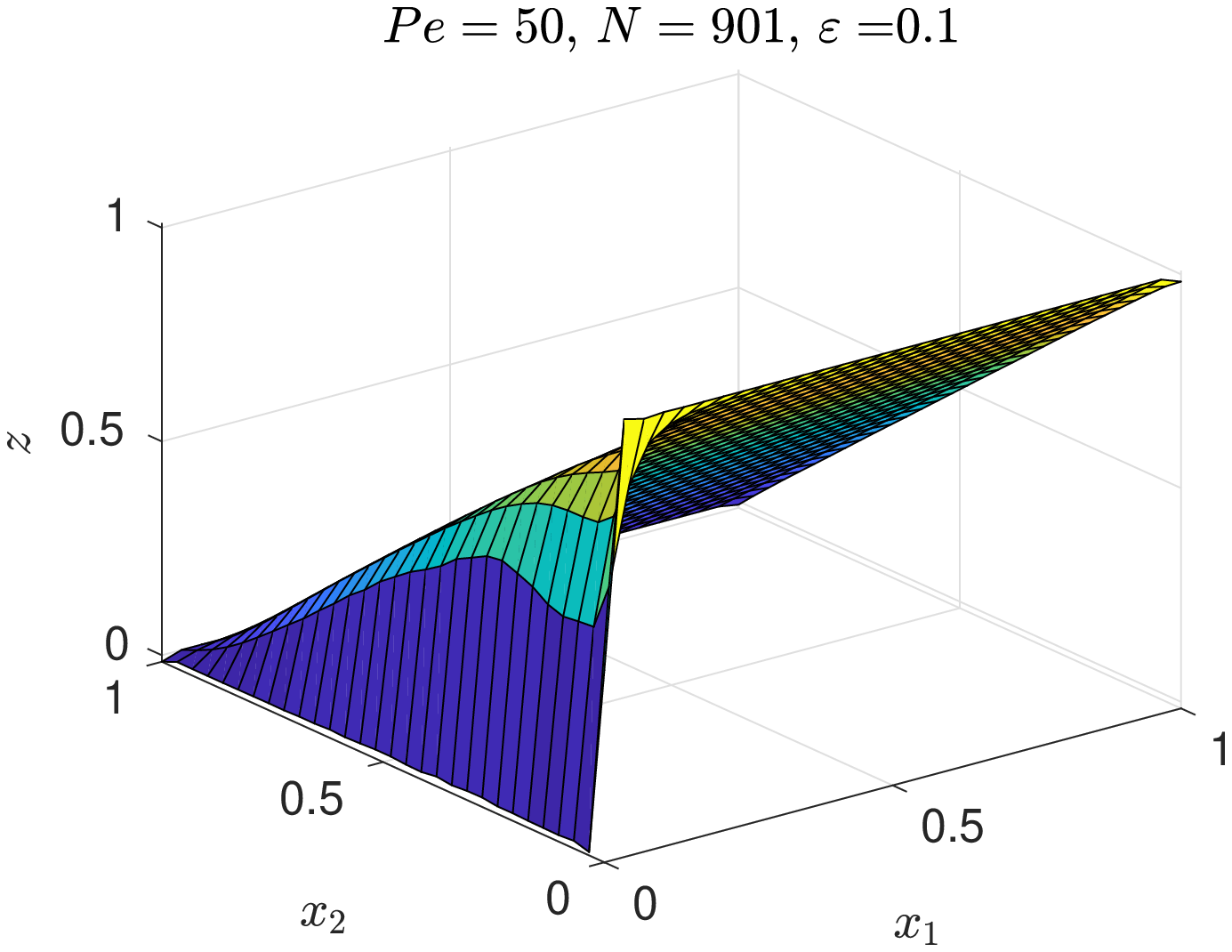}
    \includegraphics[scale=.325]{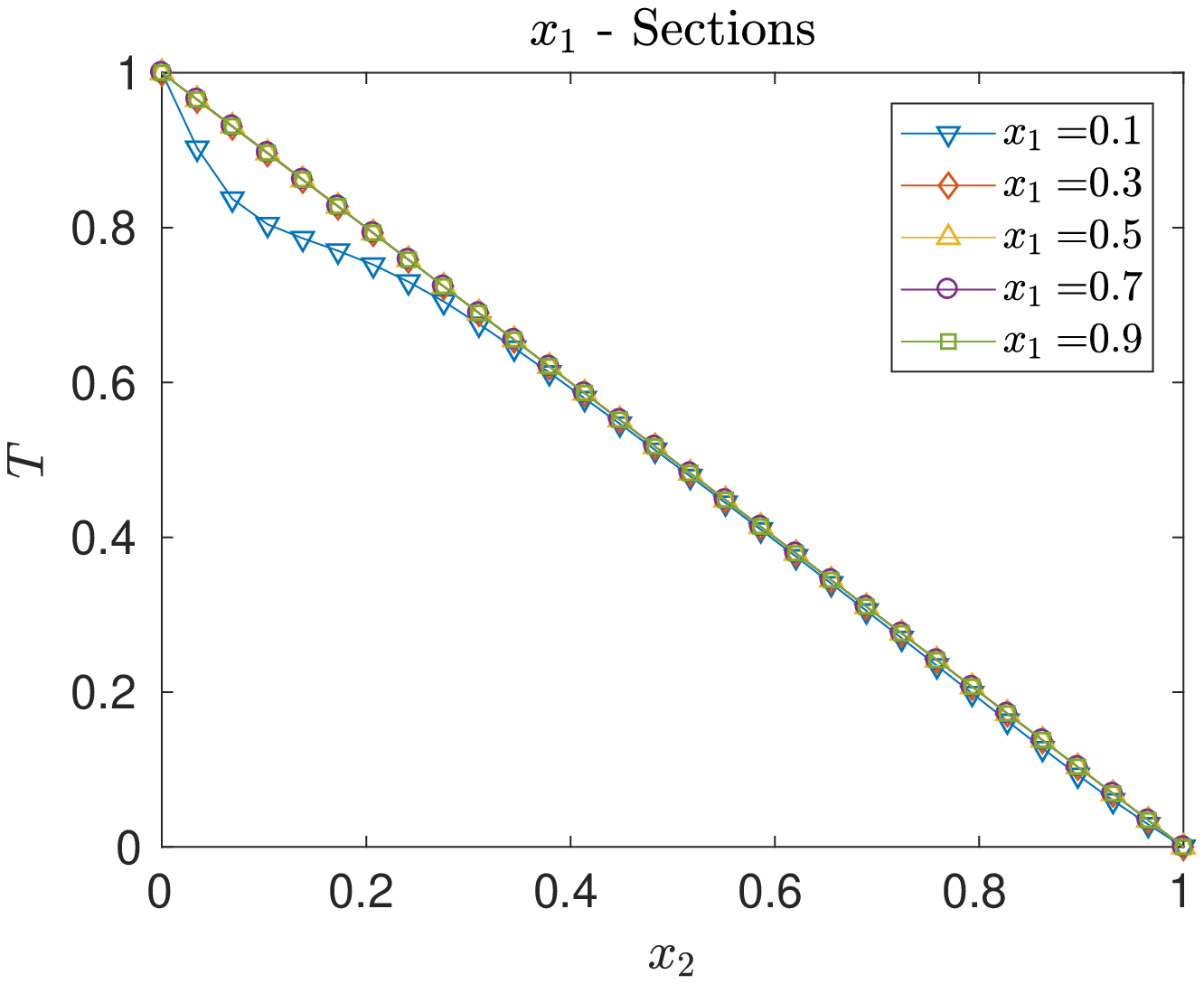}
    \includegraphics[scale=.325]{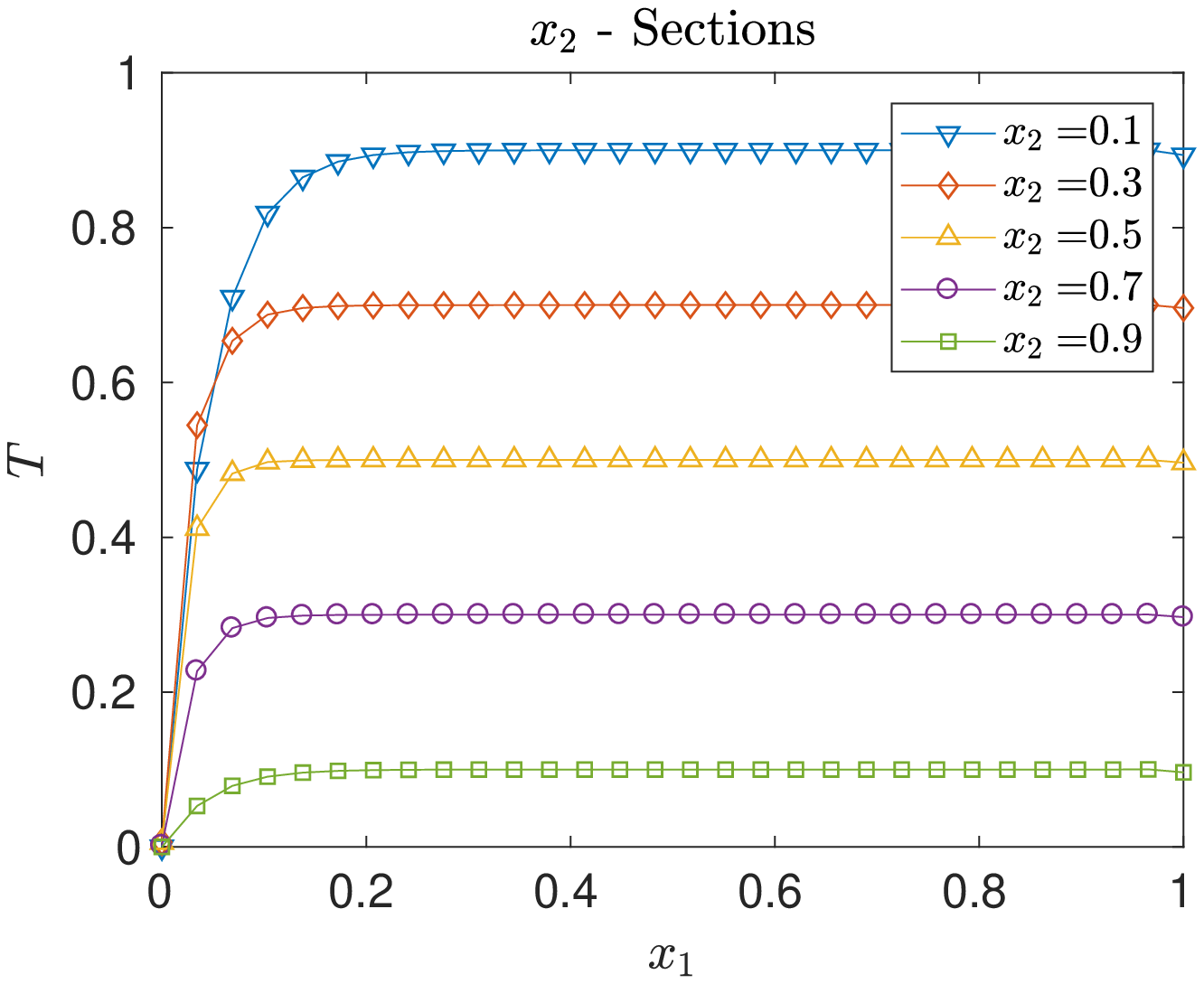}\\
   \end{center}
   \vspace{-.5cm}
  \caption{LIM RBF-QR - Quasi-uniform point distribution - $Pe=50$ - $N=901$ - $\varepsilon = 0.1$}
    \label{Fig:TBL_P2_ErrorDomain}
\end{figure}

\begin{figure}[!ht]
\centering
  \begin{center}
    \includegraphics[scale=.325]{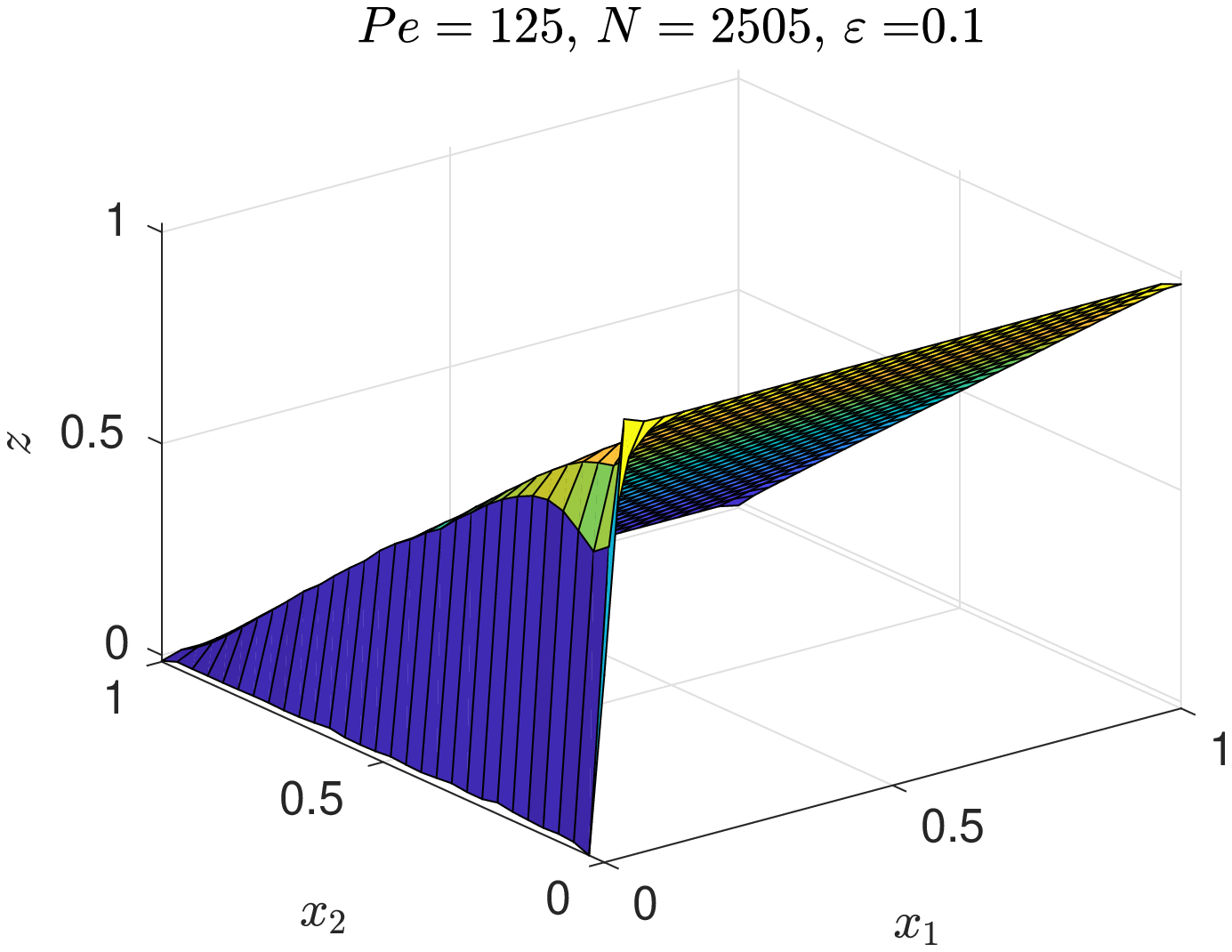}
    \includegraphics[scale=.325]{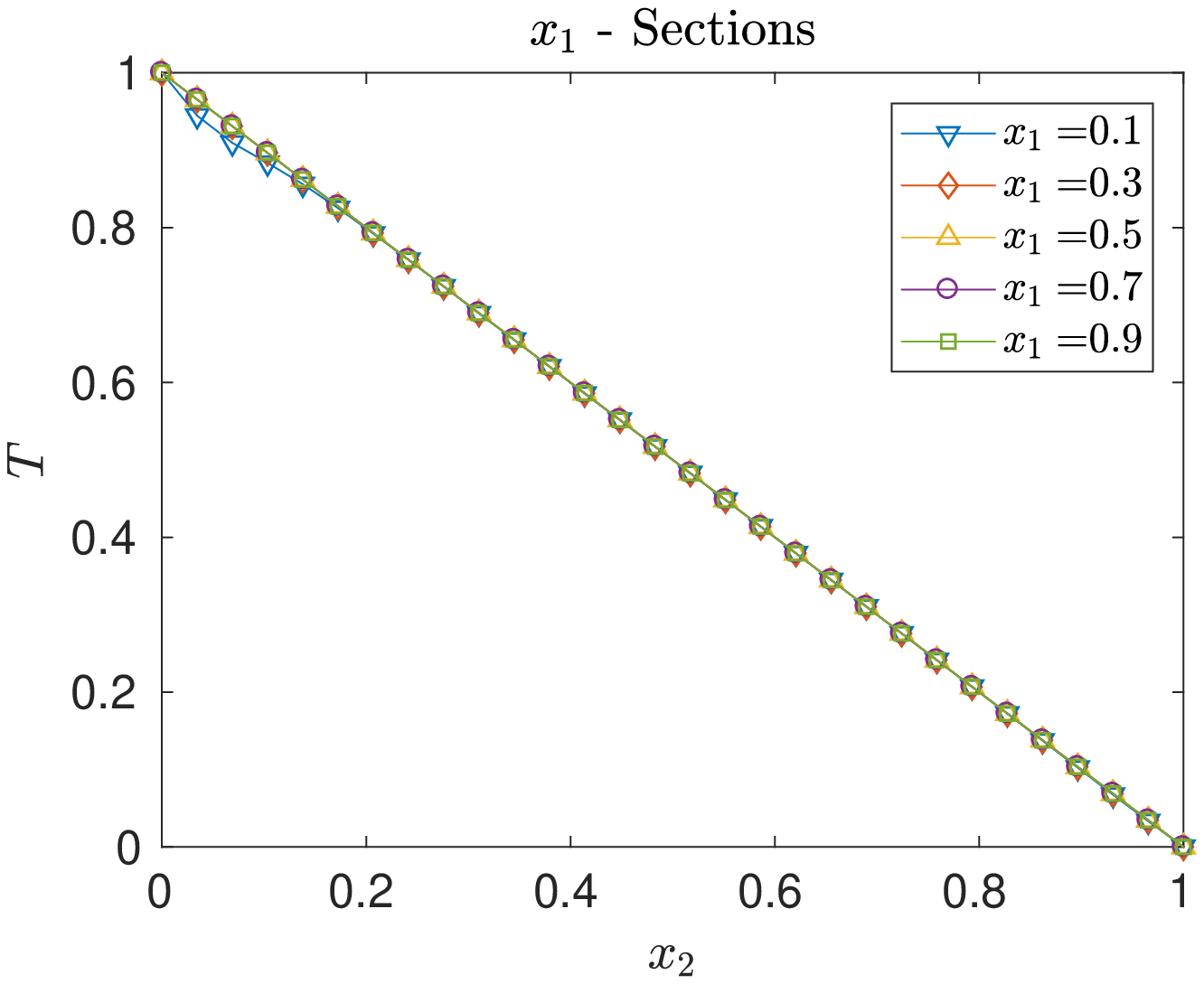}
    \includegraphics[scale=.325]{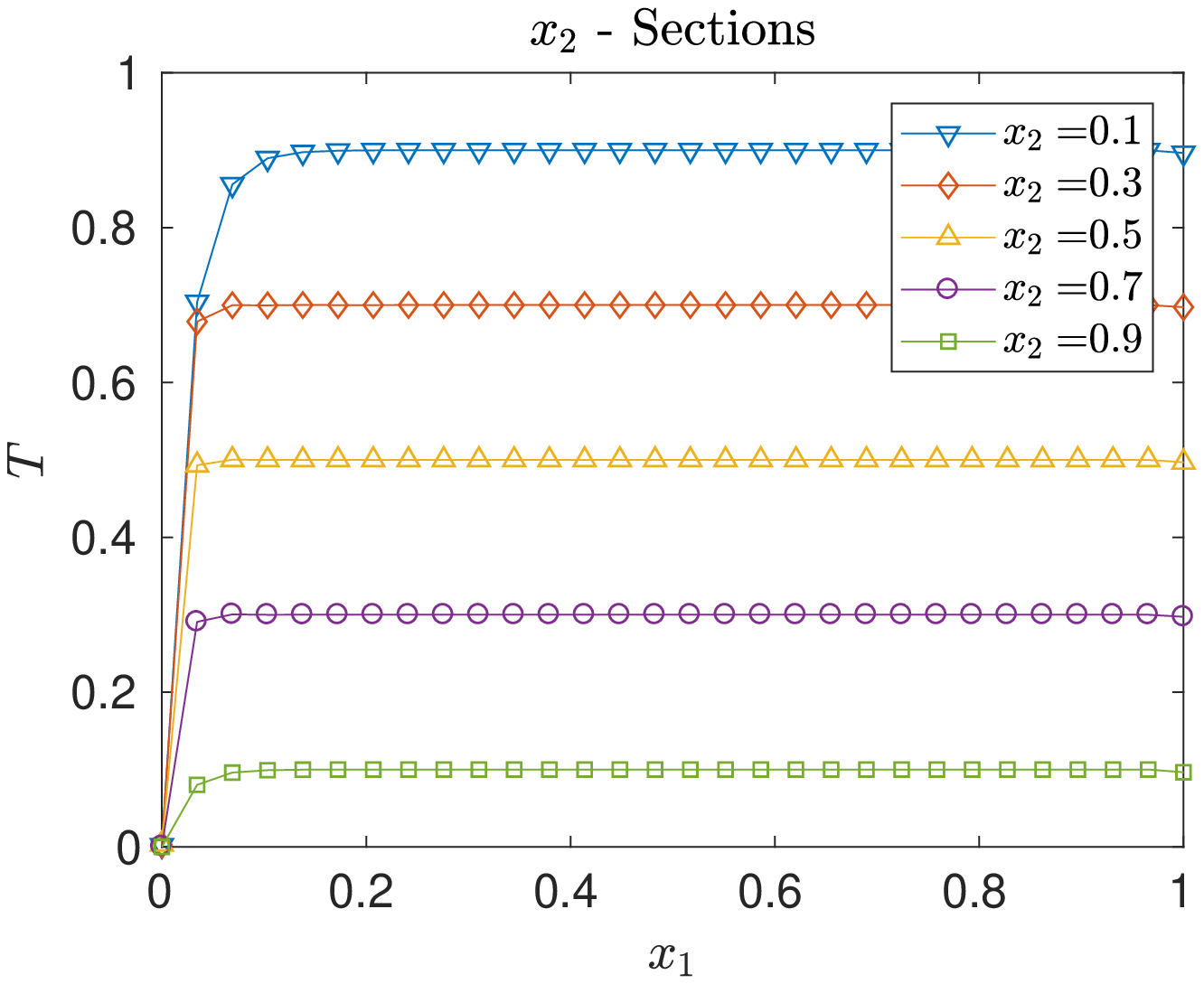}\\
   \end{center}
   \vspace{-.5cm}
  \caption{LIM RBF-QR - Quasi-uniform point distribution - $Pe=125$ - $N=2505$ - $\varepsilon = 0.1$}
  \label{Fig:TBL_P3_ErrorDomain}
\end{figure}

The objetive in this example was to obtain numerical solutions in low values of the shape parameter $\varepsilon$ for different 
quasi-uniform node sets and three different values of global Peclet number $Pe$. Fig. \ref{Fig:TBL_P1_ErrorDomain} show results 
for $Pe=0.25$, Fig.\ref{Fig:TBL_P2_ErrorDomain} for $Pe=50$ and Fig. \ref{Fig:TBL_P3_ErrorDomain} for $Pe=125$. In all Figures 
we show from left to right the approximated solution, the $x_1$ and $x_2$ sections in that order. For the reconstruction of the 
solution we use Biharmonic spline interpolation that is an interpolation of the irregular 2D data points that we obtain from the 
integral method. This interpolating surface is a linear combination of Green functions centered at each data point. For more 
references see \cite{sandwell_87,deng_tang_11}.

For all the values of the Peclet numbers studied with LIM RBF-QR, the obtained numerical results were numerically stable without 
oscillations for low levels of $\varepsilon$. In the Figures we show that to achieve $\varepsilon =0.1$ we need $N=901$ quasi-uniform 
points for the PDE with $Pe=0.25$ and $Pe=50$. For the global Peclet $Pe=125$ we need $N=2505$ interior points.

When we considered the LIM without RBF-QR and with GA RBFs for local interpolations for $u$ and $b$, we obtained that for shape 
parameters $\varepsilon<1$, the numerical approximation presented several oscillations in all cases. This is because the condition 
number of the RBF local interpolation matrix varies between $10^{17}$ and $10^{21}$.

\section{Conclusions}
\label{sec:conclusion}

In this paper, a method based on a local integral approach considering local RBF interpolation has been presented with the improvement 
of the numerical technique RBF-QR to achieve good results for low range of the shape parameter. This method was called the LIM RBF-QR. 
The robustness of this numerical method has been assessed for several elliptic PDEs with Dirichlet and Neumann BC over different domains and with 
scattered distributions as Halton or quasi-uniform points. For a Poisson equation with mixed BC over a square using LIM RBF-QR 
we improve the $RMS$ error one order of magnitude of the LRDRM with Gaussians RBF for $N=400$ uniform points. We also improved two order for 
$N=900, 1600,2500$ with $\varepsilon<1$. For a Poisson problem with Dirichlet BC we improve with $\varepsilon=1.4$ by three orders the results presented 
for the RBIEM with TPS and one order for LRDRM with MQ2 for $N=400$ interior points. For $N=900,1600,2500,3600,4900,6400$ the improvemnet was by three and two 
orders respectively using $\varepsilon=1.4,1.3,1.1,0.2,0.1,0.1$ respectively.
For the Poisson problem over the unit disk with the RBF-QR scheme, a larger region of convergence with $N=1185$ centres is observed for $L_2$-error of orders $10^{-3}$, 
$10^{-4}$ and $10^{-5}$ (with smaller values ​​of the shape parameter and greater number of points per stencil) to the obtained with Direct LIM with Gaussians RBF. 
In the case of $N=4880$ interior centres, the region of convergence obtained of orders $10^{-3}$, $10^{-4}$ and $10^{-5}$ are larger and also we achieved to a region 
of order $10^{-6}$ for small $\varepsilon$ and stencil size bigger than $n=40$ points. All these results improved the numerical $L_2$-errors presented using 
RBF-Generated Finite Difference method with polynomial augmentation for the same problem.
For the Convection-Diffusion PDE, the $L_2\%$ error for the discretization tested was improved by two orders of magnitude. Even better order of magnitude 
results are obtained for LIM RBF-QR with a discretization of $N=1127$ quasi-uniform points, with respect to results of LRDRM for $N=8000$ uniform points with 
TPS, MQ1 and MQ2. 
For the Thermal Boundary Layer PDE (numerically unstable from $Pe=2$ and without analytical solution), using LIM RBF-QR with $\varepsilon=0.1$ we obtained good results 
for the reconstruction with biharmonic splines of the numerical approximation for P\'eclet $2.5$, $50$ and $125$ values.


\section{Acknowledgments}
\label{sec:acknowledgments}

The authors would like to thank Elisabeth Larsson from Uppsala University for the collaboration on the RBF-QR imlementation methos, 
and Bengt Fornberg from University of Colorado at Boulder and Natasha Flyer from the National Center for Atmospheric Research (NCAR) 
at Boulder, CO, US for their valuable discussions and suggestions provided. 
We also acknowledge the contributions of Professor Henry Power (1950-2017) from University of Nottingham in memoriam.

\bibliographystyle{abbrv}
\bibliographystyle{PCP-article-template-v47-simple}
\bibliography{PCP-article-references}
\end{document}